\documentclass[11pt]{article}
\usepackage{amsfonts,latexsym,amssymb,amsthm,amsmath,graphicx,cases,comment}
\usepackage{amsmath}
\usepackage{paralist}
\usepackage{graphics} 
\usepackage{epsfig} 
\usepackage{epstopdf}
\usepackage{ulem}
\usepackage{todonotes}
\usepackage[titletoc,title]{appendix}
\usepackage{empheq}
\usepackage{cancel}
\usepackage{tikz}
\usetikzlibrary{arrows,shapes}
\usetikzlibrary{decorations.pathmorphing,decorations.pathreplacing}
\usetikzlibrary{calc,patterns,angles,quotes}

\usepackage[colorlinks = true, citecolor = blue]{hyperref}

\allowdisplaybreaks
\newtheorem{theorem}{Theorem}[section]
\newtheorem{thm}{Theorem}

\newtheorem{prop}{Proposition}

\theoremstyle{definition}
\newtheorem{definition}[theorem]{Definition}
\newtheorem{rmk}{Remark}

\newcommand{\be}{\begin{equation}}
\newcommand{\ee}{\end{equation}}
\newcommand{\bsubeq}{\begin{subequations}}
	\newcommand{\esubeq}{\end{subequations}}
\renewcommand{\div}{\text{div}}
\newcommand{\ds}{\displaystyle}

\newcommand{\calL}{{\mathcal{L}}}

\newcommand{\calN}{{\mathcal{N}}}

\newcommand{\calF}{{\mathcal{F}}}
\newcommand{\calB}{{\mathcal{B}}}
\newcommand{\calD}{{\mathcal{D}}}

\newcommand{\calA}{{\mathcal{A}}}

\newcommand{\calC}{{\mathcal{C}}}

\newcommand{\calM}{{\mathcal{M}}}

\newcommand{\calV}{{\mathcal{V}}}
\newcommand{\BA}{\mathbb{A}}

\newcommand{\BR}{\mathbb{R}}

\newcommand{\wti}{\widetilde}
\newcommand{\what}{\widehat}

\newcommand{\bpm}{\begin{pmatrix}}
	\newcommand{\epm}{\end{pmatrix}}

\newcommand{\bbm}{\begin{bmatrix}}
	\newcommand{\ebm}{\end{bmatrix}}

\numberwithin{equation}{section}
\numberwithin{thm}{section}
\numberwithin{rmk}{section}
\numberwithin{prop}{section}

\newcommand{\lso}{\mathbf{L}^{q}_{\sigma}(\Omega)}
\newcommand{\lqo}{L^{q}(\Omega)}
\newcommand{\lo}[1]{\mathbf{L}^{#1}_{\sigma}(\Omega)}

\newcommand{\Bso}{ \mathbf{B}^{2-\rfrac{2}{p}}_{q,p}(\Omega)}

\newcommand{\Bt}{\widetilde{\mathbf{B}}^{2-\rfrac{2}{p}}_{q,p}}
\newcommand{\Bto}{\widetilde{\mathbf{B}}^{2-\rfrac{2}{p}}_{q,p}(\Omega)}
\newcommand{\xtpq}{\mathbf{X}^T_{p,q}}
\newcommand{\xtpqs}{\mathbf{X}^T_{p,q,\sigma}}
\newcommand{\xipqs}{\mathbf{X}^{\infty}_{p,q,\sigma}}
\newcommand{\ytpq}{Y^T_{p,q}}
\newcommand{\yipq}{Y^{\infty}_{p,q}}
\newcommand{\xipq}{X^{\infty}_{p,q}}

\newcommand\rfrac[2]{{}^{#1}\!/_{#2}}
\newcommand{\Fs}{\mathbf{F}_{\sigma}}
\newcommand{\lqaq}{\big( \lso, \calD(A_q) \big)_{1-\frac{1}{p},p}}

\newcommand{\norm}[1]{\left\lVert#1\right\rVert}
\newcommand{\abs}[1]{\left\lvert#1\right\rvert}
\newcommand{\lsoo}{\mathbf{L}^{q'}_{\sigma}(\Omega)}
\newcommand{\lqoo}{L^{q'}(\Omega)}

\newcommand{\lplq}{L^p \big( 0,\infty; L^q(\Omega) \big)}
\newcommand{\lplqs}{L^p \big( 0,\infty; \lso \big)}

\newcommand{\SqsO}{\mathbf{W}^q_{\sigma}(\Omega)}

\newcommand{\VqpO}{\mathbf{V}^{q,p}(\Omega)}

\begin{document}
	
	\title{Uniform stabilization of Boussinesq systems in critical $\mathbf{L}^q$-based Sobolev and Besov spaces by finite dimensional interior localized feedback controls \thanks{The research of I.L. and R.T. was partially supported by the National Science Foundation under Grant DMS-1713506. The research of B.P. was supported by the ERC advanced grant 668998 (OCLOC) under the EU's H2020 research program.}}	
		
	\author{Irena Lasiecka \thanks{Department of Mathematical Sciences, University of Memphis, Memphis, TN 38152 USA and IBS, Polish Academy of Sciences, Warsaw, Poland.}
		\and Buddhika Priyasad \thanks{Institute for Mathematics and Scientific Computing, University of Graz, Heinrichstrasse 36, A-8010 Graz, Austria. (b.sembukutti-liyanage@uni-graz.at).}
	 \and Roberto Triggiani \thanks{Department of Mathematical Sciences, University of Memphis, Memphis, TN 38152 USA.}}


	\date{}
	\maketitle	
	\begin{abstract}
		\noindent We consider the d-dimensional Boussinesq system defined on a sufficiently smooth bounded domain, with homogeneous boundary conditions, and subject to external sources, assumed to cause instability. The initial conditions for both fluid and heat equations are taken of low regularity. We then seek to uniformly stabilize such Boussinesq system in the vicinity of an unstable equilibrium pair, in the critical setting of correspondingly low regularity spaces, by means of explicitly constructed, feedback controls, which are localized on an arbitrarily small interior subdomain. In addition, they will be minimal in number, and of reduced dimension: more precisely, they will be of dimension $(d-1)$ for the fluid component and of dimension $1$ for the heat component. The resulting space of well-posedness and stabilization is a suitable, tight Besov space for the fluid velocity component (close to $\mathbf{L}^3(\Omega$) for $ d = 3 $) and the space $L^q(\Omega$) for the thermal component, $ q > d $. Thus, this paper may be viewed as an extension of \cite{LPT.1}, where the same interior localized uniform stabilization outcome was achieved by use of finite dimensional feedback controls for the Navier-Stokes equations, in the same Besov setting.
	\end{abstract}
	
	\section{Introduction}
	
	\subsection{Controlled dynamic Boussinesq equations}\label{subsec_1.1}
	
	\noindent In this paper, we consider the following Boussinesq approximation equations in a bounded connected region $\Omega$ in $\BR^d$ with sufficiently smooth boundary $\Gamma = \partial \Omega$. More specific requirements will be given below.  Let $Q \equiv (0,T) \times \Omega$ and $\Sigma \equiv (0,T) \times \partial \Omega $ where $T > 0$. Further, let $\omega$ be an arbitrary small open smooth subdomain of the region $\Omega$, $\omega \subset \Omega$, thus of positive measure. Let $m$ denote the characteristic function of $\omega$: $ m(\omega) \equiv 1, \ m(\Omega / \omega) \equiv 0$. \\
	
	\noindent \textbf{Notation:} Vector-valued functions and corresponding function spaces will be boldfaced.
	Thus, for instance, for the vector valued ($d$-valued) velocity field or external force, we shall write
	say $\mathbf{y,f} \in \mathbf{L}^q(\Omega) $ rather than $ y,f \in (L^q(\Omega))^d$.
	\medskip

	\noindent The Boussinesq system under the action of two localized interior controls $m(x)\mathbf{u}(t,x)$ and $m(x)v(t,x)$ supported on $Q_{\omega} \equiv (0, \infty) \times \omega$ is then	
	\begin{subequations}\label{1.1}
		\begin{align}
		\mathbf{y}_t - \nu \Delta \mathbf{y} + (\mathbf{y} \cdot \nabla) \mathbf{y} - \gamma (\theta - \bar{\theta}) \mathbf{e}_d+ \nabla \pi &= m(x)\mathbf{u}(t,x) + \mathbf{f}(x)   &\text{ in } Q \label{1.1a}\\
		\theta_t - \kappa \Delta \theta + \mathbf{y} \cdot \nabla \theta &= m(x)v(t,x) + g(x) &\text{ in } Q \label{1.1b}\\
		\begin{picture}(0,0)
		\put(-175,-1){ $\left\{\rule{0pt}{40pt}\right.$}\end{picture}
		\text{div }\mathbf{y} &= 0   &\text{ in } Q \label{1.1c}\\
		\mathbf{y} = 0, \ \theta &= 0 &\text{ on } \Sigma \label{1.1d}\\
		\mathbf{y}(0,x) = \mathbf{y}_0, \quad \theta(0,x) & = \theta_0 &\text{ on } \Omega. \label{1.1e}
		\end{align}
	\end{subequations}	
	
	\noindent In the Boussinesq approximation system, $\mathbf{y} = \{y_1, \dots, y_d\}$ represents the fluid velocity, $\theta$ the scalar temperature of the fluid, $\nu$ the kinematic viscosity coefficient, $\kappa$ the thermal conductivity. The scalar function $\pi$ is the unknown pressure. The term $\mathbf{e}_d$ denotes the vector $(0, \dots, 0, 1)$. Moreover $\ds \gamma = \rfrac{\bar{g}}{\bar{\theta}}$ where $\bar{g}$ is the acceleration due to gravity and $\bar{\theta}$ is the reference temperature. The $d$-vector valued function $\mathbf{f}(x)$ and scalar function $g(x)$ correspond to an external force acting on the Navier-Stokes equations and a heat source density acting on the heat equation, respectively. They are given along with the I.C.s $\mathbf{y}_0$ and $\theta_0$. Note that $ \mathbf{y} \cdot \nabla \theta = \div (\theta \mathbf{y})$.\\
	
	\noindent The Boussinesq system  models heat transfer in a viscous incompressible heat conducting fluid. It consists of  the Navier-Stokes equations (in the vector velocity $\mathbf{y}$) coupled with the convection-diffusion equation (for the scalar temperature $\theta$). The external body force $\mathbf{f}(x)$ and the heat source density $g(x) $ may render the overall system unstable in the technical sense described below by \eqref{1.21}. The goal of the paper is to  exploit  the  localized  controls $\mathbf{u}$ and $v$, sought to be finite dimensional and in  feedback form, in order to stabilize the overall system. As an additional benefit of our investigation, the feedback fluid component of $ \mathbf{u} $ will be of reduced dimension $ ( d-1) $ rather than $ d $, while the feedback heat component of $ \mathbf{u} $ will be of 1-dimensional. This is a consequence of the Unique Continuation Property expressed by Theorem \ref{T1.4} for the adjoint static problem. As far as practical applications are concerned, one may consider the situation of controlling the temperature and humidity in a bounded  restricted environment - see \cite{BH:2013}, \cite{BHH:2016} for an eloquent description  of  the physical phenomenon. Due to the physical  significance of  the Boussinesq system, the problem of its stabilization has been considered in the literature - with both localized and boundary controls - following of course prior developments concerning the Navier-Stokes model alone. See subsection \ref{subsec_1.4} for a view of the literature.\\
	
	\noindent \textbf{Motivation: why studying uniform stabilization of the Boussinesq problem (\ref{1.1}) in the Besov functional setting (in the fluid component) of the present paper?}\\
	
	\noindent In short: Stimulated by recent research achievements \cite{LPT.1}, \cite{LPT.2} on the uniform stabilization of the Navier-Stokes equations - to be elaborated below - the present paper sets the stage as a preliminary, needed step toward the authors' final goal of solving the Boussinesq uniform stabilization problem with \uline{finite dimensional feedback controls} acting \underline{tangentially on a part of the boundary} in dimension $d=3$. We next illustrate the above assertion. While several results have been obtained in 2 and 3 dimensional cases, the main open question after 20 years or so - in both the Navier-Stokes as well as the Boussinesq systems - is whether it is possible at all to uniformly stabilize the $3-d$ equations with \uline{finite dimensional feedback controls} acting on the \uline{boundary}. Indeed, the main obstacle in this $3-d$ case is due to compatibility conditions between boundary terms and initial conditions, resulting from the necessity of using sufficiently regular [differentiable]  solutions for $d=3$. As it was originally the case with the Navier-Stokes equations, the Boussinesq case in the literature considers until now \cite{BT:2011}, \cite{Lef}, \cite{RRR:2019}, \cite{Wang} the problem within the  Hilbert framework. Due to the Stokes non-linearity in $3-d$, this necessitates in the Navier-Stokes equations - whether alone or as a component of the Boussinesq system - the need of working with  sufficiently high order Sobolev spaces which recognize boundary conditions, and thus impose compatibility conditions. This has so far prevented the goal of achieving in the $3-d$ case the desired result of \uline{finite-dimensionality} of the \uline{tangential boundary} feedback stabilizing control in full generality, beyond the special case of  \cite{BLT1:2006}, \cite{LT2:2015}, where the initial conditions are compactly supported. See also \cite[below (4.122)]{BT:2011}. For Navier-Stokes equations, this $3d$-recognized long-standing open problem on the \uline{finite dimensionality} of the \uline{localized tangential boundary} feedback stabilizing control in full generality, was recently settled in the affirmative in \cite{LPT.2}. Its treatment strongly benefited - technically and conceptually - from the preliminary successful analysis of the test case on uniform stabilization of the ($d=2$ and) $d=3$ - Navier-Stokes equations, by means of \uline{finite dimensional} feedback stabilizing controls localized on an arbitrarily subset $ \omega $ of the interior $\Omega$, precisely in the same Besov setting \cite{LPT.1}. To this end, it was necessary to introduce a radically different functional setting. It is based on Besov spaces of tight indices, which do not recognize boundary conditions, see Remark 1.1, ``close" to the well-known critical space $\mathbf{L}^3$ in space and $L^1$ in time for $ d=3$. As said, it has the virtue of not recognizing boundary conditions, see Remark 1.1, while still being adequate for handling the N-S nonlinearity in $3-d$. Circumventing the obstacle of compatibility conditions on the boundary, which prevents (for $d \geq  3$) finite dimensionality of stabilizing boundary controls in all prior Hilbert-based treatments, has been a major predicament in $3d$-boundary feedback stabilization of N-S flows. In contrast the Besov setting develops a theory within non-Hilbertian structures, where solutions under considerations have higher integrability but almost no differentiability. Thus, such Besov space setting, which has already been used for well-posedness of systems of incompressible flows, eg \cite{DM:2015}, \cite{MZ:2000}, is here a forced necessity introduced by the authors in order to solve the uniform stabilization problem of the Navier-Stokes equations with a \underline{localized tangential boundary} feedback control which moreover is \underline{finite dimensional} also for $d=3$. Though the jump from localized interior to localized boundary (tangential) stabilizing \underline{finite dimensional} ($d=3$) controls offers serious additional challenges - conceptual as well as technical - the solution in \cite{LPT.1} of the former case proved very useful in achieving the solution of the latter case in \cite{LPT.2}. Accordingly, the present paper is the generalization to the Boussinesq system of the local interior stabilization problem of the Navier-Stokes problem in \cite{LPT.1}. It is intended to be equally enlightening and useful toward the final goal, which is the uniform stabilization problem of the Boussinesq system by \underline{finite dimensional} localized tangential \underline{boundary} controls; that is, the direct counterpart of the Navier-Stokes result in \cite{LPT.2}.\\
	
	\noindent Accordingly, in seeking a $3d$ extension of such tangential boundary stabilization problem from the $3d$-Navier-Stokes equations to the $3d$-Boussinesq system, in the present paper we first consider the new Besov-setting in the context of interior localized controls. \\
	
	\noindent As in these references, the analysis of seeking stabilizing finite dimensional controls is \uline{spectral based} \cite{BLT33:2006}, \cite{LT1:2015}, \cite{LT2:2015}, unlike the Riccati-based Hilbert approaches originally used in the Navier-Stokes equations \cite{BT:2004}, \cite{BLT1:2006}, \cite{BLT22:2006} and also in the Boussinesq case \cite{BHH:2016}, following \cite{LT3:2015}.

	\subsection{Stationary Boussinesq equations}
	
	\noindent Our starting point is the following result.
	
	\begin{thm} \label{Thm-1.1}
		Consider the following steady-state Boussinesq system in $\Omega$
		\begin{subequations}\label{1.2}
			\begin{align}
			- \nu \Delta \mathbf{y}_e + (\mathbf{y}_e \cdot \nabla) \mathbf{y}_e - \gamma (\theta_e - \bar{\theta}) \mathbf{e}_d+ \nabla \pi_e &= \mathbf{f}(x)   &\text{in } \Omega \label{1.2a}\\
			-\kappa \Delta \theta_e + \mathbf{y}_e \cdot \nabla \theta_e &= g(x) &\text{in } \Omega \label{1.2b}\\
			\begin{picture}(0,0)
			\put(-165,5){ $\left\{\rule{0pt}{34pt}\right.$}\end{picture}
			\text{div } \mathbf{y}_e &= 0   &\text{in } \Omega \label{1.2c}\\
			\mathbf{y}_e = 0, \ \theta_e &= 0 &\text{on } \partial \Omega. \label{1.2d}
			\end{align}
		\end{subequations}
		Let $1 < q < \infty$. For any $\mathbf{f},g \in \mathbf{L}^q(\Omega), L^q(\Omega)$, there exists a solution (not necessarily unique) $(\mathbf{y}_e,\theta_e, \pi_e) \in (\mathbf{W}^{2,q}(\Omega) \cap \mathbf{W}^{1,q}_{0}(\Omega)) \times (W^{2,q}(\Omega) \cap W^{1,q}_{0}(\Omega)) \times (W^{1,q}(\Omega)/\mathbb{R})$.
	\end{thm}
	\noindent See \cite{Ac}, \cite{AAC.1}, \cite{AAC.2} for $ q \ne 2$. In the Hilbert space setting, see \cite{CF:1980}, \cite{FT:1984}, \cite{S-R-R}, \cite{V-R-R}, \cite{Kim}.
	
	\subsection{A first quantitative description of the main goal of the present paper}\label{subsec-3}
	
	The starting point of the present paper is the following: that under a given external force $\mathbf{f}(x)$ for the fluid equations, a given heat source $g(x)$ for the thermal equation, and given viscosity coefficient $\nu$ and thermal conductivity $\kappa$, the equilibrium solution $\{\mathbf{y}_e, \theta_e\}$ is unstable, in a quantitative sense to be made more precise in sub-section \ref{subsec-8}, specifically in \eqref{1.21}. This will mean that the free dynamics \underline{linear} operator $\mathbb{A}_q$ defined in \eqref{1.18} - which has compact resolvent, and is the  generator of a s.c. analytic semigroup in the appropriate functional setting -  has $N$ unstable eigenvalues.\\
	
	\noindent The main goal of the present paper is then - at first qualitatively - \uline{to feedback stabilize the non-linear Boussinesq model (\ref{1.1}) subject to rough (non-smooth) initial conditions $\{\mathbf{y}_0,\theta_0\}$, in the vicinity of an (unstable) equilibrium solution $\{\mathbf{y}_e,\theta_e\}$ in (\ref{1.2})}, by means of a \uline{finite dimensional} localized feedback control pair $\{m\mathbf{u},mv\}$. Thus this paper pertains to the general issue of ``turbulence suppression or attenuation" in fluids. The general topic of turbulence suppression (or attenuation) in fluids has been the object of many studies over the years, mostly in the engineering literature – through experimental studies and via numerical simulation - and under different geometrical and dynamical settings. The references cited in the present paper by necessity pertain mostly to the mathematical literature, and most specifically on the \uline{localized interior control setting of problem (\ref{1.1})} \cite{Lef}, \cite{Wang}. A more precise description thereof is as follows: \textit{establish interior localized exponential stabilization of problem (\ref{1.1}) near an unstable equilibrium solution $\{\mathbf{y}_e,\theta_e\}$ by means of a \underline {finite dimensional localized}, spectral-based feedback control pair $\{m\mathbf{u},mv\}$, in the important case of initial condition $\mathbf{y}_0$ of low regularity}, as technically expressed by $\mathbf{y}_0$ being in a suitable Besov space with tight indices, and $\theta_0$ being in a corresponding $L^q$-space $ q > d $. More precisely, the resulting state space for the pair $\{\mathbf{y},\theta\}$, where uniform stabilization will be achieved is the space
	\begin{equation}\label{N4-111}
	\quad \VqpO \equiv \Bto \times \lqo, \ 1 < p < \frac{2q}{2q - 1} ; \ q > d , \, d = 2,3,
	\end{equation}
	where $ \Bto $ is a suitable subspace, see below in (\ref{A-1.16bb}), of the Besov space.
	\begin{equation}\label{I-1.3}
	\big( \mathbf{L}^q(\Omega), \mathbf{W}^{2,q}(\Omega) \big)_{1-\rfrac{1}{p},p} = \Bso, \quad 1 < p < \frac{2q}{2q -1}; \ q > d,\ d = 2,3,
	\end{equation}
	as a real interpolation space between $ \mathbf{L}^q ( \Omega)$ and $ \mathbf{W}^{2,q}( \Omega)  $.
	This setting will be further elaborated after introducing the Helmhotz decomposition below. In particular, local exponential stability for the velocity field $\mathbf{y}$ near an equilibrium solution $\mathbf{y}_e$ will be achieved in the topology of the Besov subspace $ \Bto $ in (\ref{A-1.16bb}). Note the tight index: $\ds 1 < p < \rfrac{6}{5} $ for $q > d = 3$, and $ 1 < p <  \rfrac{4}{3}$ for $ d=2 $. It will be documented below in Remark 1.1 that in such setting, the compatibility conditions on the boundary of the initial conditions \uline{are not recognized}. This feature is \uline{precisely our key objective within the stabilization problem}. The fundamental reason is that such feature will play a critical role in our aforementioned research project which consists in showing local uniform stabilization of the Boussinesq system, by means of finite dimensional, localized tangential boundary controls also in dimension dim$\Omega = d = 3$; the perfect counterpart of \cite{LPT.2}. In the case of the Navier-Stokes equations, uniform stabilization in the vicinity of a steady state solution by means of a localized tangential \uline{boundary} feedback control which in addition is also \uline{finite dimensional} for $d=3$ was a recognized open problem. It was recently settled in the affirmative in \cite{LPT.2} precisely in the Besov subspace $ \Bto $, $\ds 1 < p < \rfrac{6}{5} $ for $q > d = 3$. For $d=3$, such space is `close' to $\mathbf{L}^3$. Thus this departure from the Hilbert space theory - the latter present in all other publications related to the stabilization of Navier-Stokes or Boussinesq systems - allows one to trade differentiability with integrability. This then avoids the obstacle of compatibility conditions while preserving just enough regularity in order to carry-out the nonlinear analysis in $3-d$.\\
	
	\noindent \textbf{Criticality of the space $\mathbf{L}^3$ for $d=3$:} We now expand on the issue regarding the `criticality' of the space $\mathbf{L}^3(\Omega)$. In the case of the \uline{uncontrolled} N-S equations defined on the \uline{full space} $\BR^3$, extensive research efforts have recently lead to the discovery that the space $\mathbf{L}^3(\BR^3)$ is a `critical' space for the issue of well-posedness. Assuming that some divergence free initial data in $\mathbf{L}^3(\BR^3)$ would produce finite time singularity, then there exists a so-called minimal blow-up initial data in $\mathbf{L}^3(\BR^3)$ \cite{GKP:2010}, \cite{JS:2013}, \cite{ESS:1991}. More details in the context of the controlled 3-d Navier-Stokes equations are given in \cite{LPT.1}, \cite{LPT.2}.\\
	
	\noindent Thus, these latter two references manage to solve the uniform stabilization problem for the controlled N-S equations in a correspondingly related low-regularity function space setting. Our present paper is then an extension of \cite{LPT.1} to the Boussinesq system \eqref{1.1}. A further justification of our low-regularity level of the Besov space in (\ref{I-1.3}) is provided by  the final goal of our line of research. Based critically on said low-regularity level of the Besov space \eqref{A-1.16bb}, \uline{which does not recognize compatibility conditions on the boundary} of the initial conditions, Remark 1.1- we shall seek to solve in the affirmative a presently open problem by showing that uniform stabilization of the Boussinesq system is possible by localized tangential boundary feedback controls which moreover are finite dimensional also in dim $\Omega = d = 3$. This result will be an extension of \cite{LPT.2} in the case of the $3-d$ N-S equations.
	
	\subsection{Helmholtz decomposition}
	\label{subsec_1.5}
	A first difficulty one faces in extending the local exponential stabilization result for the interior localized problem (\ref{1.1}) from the Hilbert-space setting in \cite{Lef}, \cite{Wang} to the $\mathbf{L}^q$-based setting is the question of the existence of a Helmholtz (Leray) projection for the domain $\Omega$ in $\mathbb{R}^d$. More precisely: Given an open set $\Omega \subset \mathbb{R}^d$, the Helmholtz decomposition answers the question as to whether $\mathbf{L}^q(\Omega)$ can be decomposed into a direct sum of the solenoidal vector space $\lso$ and the space $\mathbf{G}^q(\Omega)$ of gradient fields. Here,	
	\begin{equation}\label{1.3}
	\begin{aligned}
	\lso &= \overline{\{\mathbf{y} \in \mathbf{C}_c^{\infty}(\Omega): \div \ \mathbf{y} = 0 \text{ in } \Omega \}}^{\norm{\cdot}_q}\\
	&= \{\mathbf{g} \in \mathbf{L}^q(\Omega): \div \ \mathbf{g} = 0; \  \mathbf{g}\cdot \nu = 0 \text{ on } \partial \Omega \},\\
	& \hspace{3cm} \text{ for any locally Lipschitz domain } \Omega \subset \mathbb{R}^d, d \geq 2 \\
	\mathbf{G}^q(\Omega) &= \{\mathbf{y} \in \mathbf{L}^q(\Omega):\mathbf{y} = \nabla p, \ p \in W_{loc}^{1,q}(\Omega) \} \ \text{where } 1 \leq q < \infty.
	\end{aligned}
	\end{equation}
	
	\noindent Both of these are closed subspaces of $\mathbf{L}^q$.
	
	\begin{definition}\label{Def-1.1}
		Let $1 < q < \infty$ and $\Omega \subset \mathbb{R}^n$ be an open set. We say that the Helmholtz decomposition for $\mathbf{L}^q(\Omega)$ exists whenever $\mathbf{L}^q(\Omega)$ can be decomposed into the direct sum (non-orthogonal)
		\begin{equation}
		\mathbf{L}^q(\Omega) = \lso \oplus \mathbf{G}^q(\Omega).\label{1.4}
		\end{equation}
		The unique linear, bounded and idempotent (i.e. $P_q^2 = P_q$) projection operator $P_q:\mathbf{L}^q(\Omega) \longrightarrow \lso$ having $\lso$ as its range and $\mathbf{G}^q(\Omega)$ as its null space is called the Helmholtz projection. Additional information of use in the present paper is given in Appendix \ref{app-B} below, say Proposition \ref{I-Prop-1.2}.
	\end{definition}
	
	\noindent This is an important property in order to handle the incompressibility condition $\div \ \mathbf{y} \equiv 0$. For instance, if such decomposition exists, the Stokes equation (say the linear version of (\ref{1.1a}) with control $\mathbf{u} \equiv 0$ and no coupling, thus $\theta \equiv 0$) can be formulated as an equation in the $\mathbf{L}^q$ setting. Here below we collect a subset of known results about Helmholtz decomposition. We refer to \cite[Section 2.2]{HS:2016}, in particular to the comprehensive Theorem 2.2.5 in this reference, which collects domains for which the Helmholtz decomposition is known to exist. These include the following cases:
	
	\begin{enumerate}[(i)]
		\item any open set $\Omega \subset \mathbb{R}^d$ for $q = 2$, i.e. with respect to the space $\mathbf{L}^2(\Omega)$; more precisely, for $q = 2$, we obtain the well-known orthogonal decomposition (in the standard notation, where $\nu$ = unit outward normal vector on $\Gamma$) \cite[Prop 1.9, p 8]{CF:1980}
		\begin{subequations}\label{I-1.6}
			\begin{align}
			\mathbf{L}^2(\Omega) &= \mathbf{H} \oplus \mathbf{H}^{\perp}\\
			\mathbf{H} &= \{ \boldsymbol{\phi} \in \mathbf{L}^2(\Omega): \div \ \boldsymbol{\phi} \equiv 0 \text{ in } \Omega; \ \boldsymbol{\phi} \cdot \nu \equiv 0 \text{ on } \Gamma \}\\
			\mathbf{H}^{\perp} &= \{ \boldsymbol{\psi} \in \mathbf{L}^2(\Omega): \boldsymbol{\psi} = \nabla h, \ h \in H^1(\Omega) \};
			\end{align}
		\end{subequations}
		\item a bounded $C^1$-domain in $\mathbb{R}^d$ \cite{FMM:1998}, $1 < q < \infty $ \cite[Theorem 1.1 p 107, Theorem 1.2 p 114]{Ga:2011} for $C^2$-boundary
		\item a bounded Lipschitz domain $\Omega \subset \mathbb{R}^d \ (d = 3)$ and for $\frac{3}{2} - \epsilon < q < 3 + \epsilon$ sharp range \cite{FMM:1998};
		\item a bounded convex domain $\Omega \subset \mathbb{R}^d, d \geq 2, 1 < q < \infty$ \cite{FMM:1998}.
	\end{enumerate}
	
	\noindent On the other hand, on the negative side, it is known that there exist domains $\Omega \subset \mathbb{R}^d$ such that the Helmholtz decomposition does not hold for some $q \neq 2$ \cite{MS:1986}.\\
	
	\noindent \textbf{Assumption (H-D)}. Henceforth in this paper, we assume that the bounded domain $\Omega \subset \mathbb{R}^d$ under consideration admits a Helmholtz decomposition for the values of $q, \ 1 < q < \infty$, here considered at first, for the linearized problem (\ref{1.19}) or \eqref{N1.23} below. The final results Theorem \ref{N-Thm-2.2} through \ref{N-Thm-2.5} for the non-linear problem (\ref{1.1}) will require $q > d$, see (\ref{3.24}), in the case of interest $d = 2, 3$. \\
	
	\noindent We can now provide further critical information on the Besov space  $ \Bto$    which is the fluid component of the state space (\ref{N4-111}) where uniform stabilization takes place.\\
	
	\noindent \textbf{Definition of Besov spaces $\mathbf{B}^s_{q,p}$ on domains of class $C^1$ as real interpolation of Sobolev spaces:} Let $m$ be a positive integer, $m \in \mathbb{N}, 0 < s < m, 1 \leq q < \infty,1 \leq p \leq \infty,$ then we define \cite[p 1398]{GGH:2012} the Besov space

	\begin{equation} \label{1.8}
	\mathbf{B}^{s}_{q,p}(\Omega) = (\mathbf{L}^q(\Omega),\mathbf{W}^{m,q}(\Omega))_{\frac{s}{m},p} 
	\end{equation}
	as a real interpolation space between $\mathbf{L}^q(\Omega)$ and $\mathbf{W}^{m,q}(\Omega)$.This definition does not depend on $\ds m \in \mathbb{N}$ \cite[p xx]{W:1985}. This clearly gives
	\begin{equation} \label{1.9}
	\mathbf{W}^{m,q}(\Omega) \subset \mathbf{B}_{q,p}^s(\Omega) \subset \mathbf{L}^q(\Omega) \quad \text{ and } \quad \norm{\mathbf{y}}_{\mathbf{L}^q(\Omega)} \leq C \norm{\mathbf{y}}_{\mathbf{B}_{q,p}^s(\Omega)}.
	\end{equation}

	\noindent We shall be particularly interested in the following special real interpolation space of the $\mathbf{L}^q$ and $\mathbf{W}^{2,q}$ spaces $\Big( m = 2, s = 2 - \frac{2}{p} \Big)$:
	\begin{equation}\label{A-1.151}
	\mathbf{B}^{2-\frac{2}{p}}_{q,p}(\Omega) = \big(\mathbf{L}^q(\Omega),\mathbf{W}^{2,q}(\Omega) \big)_{1-\frac{1}{p},p}.
	\end{equation}
	\noindent Our interest in (\ref{A-1.151}) is due to the following characterization \cite[Thm 3.4]{HA:2000}, \cite[p 1399]{GGH:2012}: if $A_q$ denotes the Stokes operator to be introduced in (\ref{N1.9}) below, then
	\begin{subequations}\label{A-1.16}
		\begin{align}
		\Big( \lso,\mathcal{D}(A_q) \Big)_{1-\frac{1}{p},p} &= \Big\{ \mathbf{g} \in \Bso : \text{ div } \mathbf{g} = 0, \ \mathbf{g}|_{\Gamma} = 0 \Big\} \quad \text{if } \frac{1}{q} < 2 - \frac{2}{p} < 2 \label{A-1.16aa}\\
		\Big( \lso,\mathcal{D}(A_q) \Big)_{1-\frac{1}{p},p} &= \Big\{ \mathbf{g} \in \Bso : \text{ div } \mathbf{g} = 0, \ \mathbf{g}\cdot \nu|_{\Gamma} = 0 \Big\} \equiv \Bto \label{A-1.16bb}\\
		&\hspace{3.5cm} \text{ if } 0 < 2 - \frac{2}{p} < \frac{1}{q}; \text{ or } 1 < p < \frac{2q}{2q - 1}.\nonumber
		\end{align}	
	\end{subequations}
	
	\begin{rmk}  \label{remark1-11}
		Notice that, in (\ref{A-1.16bb}), the condition $\ds \mathbf{g} \cdot \nu |_{\Gamma} = 0$ is an intrinsic condition of the space $\ds \lso$ in (\ref{1.3}), not an extra boundary condition as $\ds \mathbf{g}|_{\Gamma} = 0$ in (\ref{A-1.16aa}).
	\end{rmk}
	
	\noindent \textbf{Orientation:} As already noted,ultimately, we shall seek to obtain uniform feedback stabilization of the fluid component $\mathbf{y}$ in the Besov subspace $\Bt(\Omega)$, dim $\ds \Omega = d < q < \infty$, $ \ds 1 < p < \rfrac{2q}{2q-1}$, defined by real interpolation in \eqref{I-1.3}, (\ref{A-1.16bb}); The reason being that such a space \underline{ does not recognize boundary conditions}, as noted above in Remark \ref{remark1-11}. Analyticity and maximal regularity of the Stokes problem will require $p>1$.\\
	
	\noindent By way of orientation, we state at the outset two main points. For the linearized $\mathbf{w}$-problem (\ref{1.19}) below in the feedback form (\ref{2.9}) or \eqref{N4-7}, the corresponding well-posedness and global feedback uniform stabilization result, Theorem \ref{NThm-2.1} or Theorem \ref{Thm-2.1}, holds in general for $1 < q < \infty$. Instead, the final, main well-posedness and feedback uniform, local stabilization results, Theorems \ref{N-Thm-2.2} through \ref{N-Thm-2.5}, the latter for the nonlinear feedback problem (\ref{I-2.27}) or (\ref{I-2.28}) corresponding to the original problem (\ref{1.1}), will require $q > 3$ to obtain the embedding $\ds \mathbf{W}^{1,q}(\Omega) \hookrightarrow \mathbf{L}^{\infty}(\Omega)$ in our case of interest $d = 3$, see (\ref{3.24}), hence $\ds 1 < p < \rfrac{6}{5}$; and $q > 2$, in the $d = 2$-case, hence $\ds 1 < p < \rfrac{4}{3}$.

	\subsection{Translated Nonlinear Boussinesq Problem and its Abstract Model}
	
	\noindent \textbf{PDE Model:} We return to Theorem \ref{Thm-1.1} which provides an equilibrium triplet $\{\mathbf{y}_e, \theta_e, \pi_e\}$. Then, we translate by $\{\mathbf{y}_e, \theta_e, \pi_e\}$ the original Boussinesq problem (\ref{1.1}). Thus we introduce new variables
	\begin{align}\label{1.5a}
	\mathbf{z} = \mathbf{y} - \mathbf{y}_e\ \  (\mbox{a $d$-vector}), \quad h = \theta - \theta_e \ \ (\mbox{a scalar}), \quad \chi = \pi - \pi_e \ \  (\mbox{a scalar})
	\end{align}
	and obtain the translated problem
	\begin{subequations}    \label{1.5}
		\begin{align}
		\mathbf{z}_t - \nu \Delta \mathbf{z} + (\mathbf{y}_e \cdot \nabla)\mathbf{z} + (\mathbf{z} \cdot \nabla)\mathbf{y}_e + (\mathbf{z} \cdot \nabla) \mathbf{z} - \gamma h \mathbf{e}_d + \nabla \chi &= m\mathbf{u} \quad \text{ in } Q \label{1.5b}\\
		h_t - \kappa \Delta h + \mathbf{y}_e \cdot \nabla h + \mathbf{z} \cdot \nabla h + \mathbf{z} \cdot \nabla \theta_e &= mv \quad \text{ in } Q \label{1.5c}\\
		\begin{picture}(0,0)
		\put(-240,-1){ $\left\{\rule{0pt}{40pt}\right.$}\end{picture}
		\div \ \mathbf{z} &= 0  \quad \text{   in } Q \label{1.5d}\\
		\mathbf{z} = 0, \ h & = 0 \quad \text{ on } \Sigma \label{1.5e}\\
		\mathbf{z}(0,x) = \mathbf{z}_0 = \mathbf{y}_0 - \mathbf{y}_e, \quad h(0,x) = h_0 & = \theta - \theta_e \ \text{ on } \Omega \label{1.5f} \\[1mm]
		L_e(\mathbf{z}) = (\mathbf{y}_e \cdot \nabla)\mathbf{z} + (\mathbf{z} \cdot \nabla)\mathbf{y}_e \quad (\mbox{Oseen perturbation}). & \label{1.5g}
		\end{align}
	\end{subequations}
	
	\noindent \textbf{Abstract Nonlinear Translated Model.}	First, for $1 < q < \infty$ fixed, the Stokes operator $A_q$ in $\lso$ with Dirichlet boundary conditions  is defined by
	\begin{equation}    \label{N1.9}
	A_q \mathbf{z} = -P_q \Delta \mathbf{z}, \quad
	\mathcal{D}(A_q) = \mathbf{W}^{2,q}(\Omega) \cap \mathbf{W}^{1,q}_0(\Omega) \cap \lso.
	\end{equation}
	The operator $A_q$ has a compact inverse $A_q^{-1}$ on $\lso$, hence $A_q$ has a compact resolvent on $\lso$.
	
	\noindent Next, we define the heat operator $B_q$ in $L^q(\Omega)$ with homogeneous Dirichlet boundary conditions
	\begin{equation}\label{E1.15}
	B_q h = -\Delta h, \quad \mathcal{D}(B_q) = W^{2,q}(\Omega) \cap W^{1,q}_0(\Omega).
	\end{equation}
	\noindent Next, we introduce the first order operator $A_{o,q}$,
	\begin{equation}        \label{N1-11}
	A_{o,q} \mathbf{z} = P_q[(\mathbf{y}_e \cdot \nabla )\mathbf{z} + (\mathbf{z} \cdot \nabla )\mathbf{y}_e], \quad \mathcal{D}(A_{o,q}) = \mathcal{D}(A_q^{\rfrac{1}{2}}) \subset \lso,
	\end{equation}
	where the $\ds \calD(A^{\rfrac{1}{2}}_q)$ is defined explicitly by complex interpolation
	\begin{subequations}\label{A-1.211}	
		\begin{equation}\label{A-1.211a}
		[ \calD(A_q), \lso ]_{1-\alpha} = \calD(A_q^{\alpha}), \ 0 < \alpha < 1, \  1 < q < \infty;
		\end{equation}
		\begin{equation}\label{A-1.21b}
		[ \calD(A_q), \lso ]_{\frac{1}{2}} = \calD(A_q^{\rfrac{1}{2}}) \equiv \mathbf{W}_0^{1,q}(\Omega) \cap \lso.
		\end{equation}
	\end{subequations}
	Thus, $A_{o,q}A_q^{-\rfrac{1}{2}}$ is a bounded operator on $\lso$, and thus $A_{o,q}$ is bounded on $\calD(A_q^{\rfrac{1}{2}})$
	\begin{equation*}
	\norm{A_{o,q}f} = \norm{A_{o,q} A_q^{-\rfrac{1}{2}} A_q^{-\rfrac{1}{2}} A_q f} \leq C_q \norm{A_q^{\rfrac{1}{2}} f}, \quad f \in \calD(A_q^{\rfrac{1}{2}}).
	\end{equation*}
	This leads to the definition of the Oseen operator for the fluid
	\begin{equation}\label{1.9.new}
	\calA_q  = - (\nu A_q + A_{o,q}), \quad \calD(\calA_q) = \calD(A_q) \subset \lso.
	\end{equation}
	\noindent Next, we introduce the first order operator $B_{o,q}$, corresponding to $B_q$:
	\begin{equation}
	B_{o,q} h = \mathbf{y}_e \cdot \nabla h, \quad \mathcal{D}(B_{o,q}) = \mathcal{D}(B_q^{\rfrac{1}{2}}) \subset \lqo.
	\end{equation}
	Thus, $B_{o,q}B_q^{-\rfrac{1}{2}}$ is a bounded operator on $ L^q (\Omega )$, and thus $B_{o,q}$ is bounded on $\calD(B_q^{\rfrac{1}{2}})$
	\begin{equation*}
	\norm{B_{o,q}f} = \norm{B_{o,q} B_q^{-\rfrac{1}{2}} B_q^{-\rfrac{1}{2}} B_q f} \leq C_q \norm{B_q^{\rfrac{1}{2}} f}, \quad f \in \calD(B_q^{\rfrac{1}{2}}).
	\end{equation*}
	This leads to the definition of the following operator for the heat
	\begin{equation}\label{E1.20}
	\calB_q  = - (\kappa B_q + B_{o,q}), \quad \calD(\calB_q) = \calD(B_q) \subset \lqo.
	\end{equation}
	Then, we define the projection of the nonlinear portion of the fluid operator in (\ref{1.5b})
	\begin{equation}
	\calN_q(\mathbf{z}) = P_q [(\mathbf{z} \cdot \nabla) \mathbf{z}], \quad \calD(\calN_q) = \mathbf{W}^{1,q}(\Omega) \cap \mathbf{L}^{\infty} (\Omega) \cap \lso.	
	\end{equation}
	(Recall that $W^{1,q}(\Omega) \hookrightarrow L^{\infty} (\Omega)$ for $q > d = \mbox{dim } \Omega$ \cite[p. 74]{SK:1989}).\\
	
	\noindent Then, we define the nonlinear coupled term of the heat equation as
	\begin{equation}\label{1.13}
	\calM_q[\mathbf{z}](h) = \mathbf{z} \cdot \nabla h, \quad \calD(\calM_q[\mathbf{z}]) = W^{1,q}(\Omega) \cap L^{\infty} (\Omega). 	
	\end{equation}
	Finally, we define the coupling linear terms as bounded operators on $\lqo, \lso$ respectively, $q > d$:
	\begin{align}
	\text{[from the NS equation]} \quad \calC_{\gamma} h &= -\gamma P_q (h \mathbf{e}_d), \quad \calC_{\gamma} \in \calL (L^q(\Omega),\lso), \label{N1-17}   \\
	\text{[from the heat equation]} \quad \calC_{\theta_e} \mathbf{z} &=  \mathbf{z} \cdot \nabla \theta_e, \quad \calC_{\theta_e} \in \calL(\lso,\lqo). \label{N1-18}
	\end{align}
	\noindent Next we apply the Helmholtz projector $P_q$ on the coupled N-S equation (\ref{1.5b}), invoke the operators introduced above and obtain the following abstract version of the controlled Boussinesq system
	\begin{subequations}\label{1.16}
		\begin{align}
		\frac{d\mathbf{z}}{dt} - \calA_q \mathbf{z} + \calN_q \mathbf{z} + \calC_{\gamma} h &= P_q(m\mathbf{u}) &\text{ in } \lso\\
		\quad \frac{dh}{dt} - \calB_q h + \calM_q[\mathbf{z}] h + \calC_{\theta_e} \mathbf{z} &= mv &\text{ in } {L}^q(\Omega)\\
		\begin{picture}(0,0)
		\put(-120,15){ $\left\{\rule{0pt}{45pt}\right.$}\end{picture}
		\mathbf{z}(x,0) &= \mathbf{z}_0(x) &\text{ in } \lso\\
		h(x,0) &= h_0(x) &\text{ in } \lqo.
		\end{align}
	\end{subequations}
	
	\noindent or in matrix form in $\lso \times \lqo \equiv \mathbf{W}^q_{\sigma}(\Omega)$
	\begin{subequations}\label{1.17}
		\begin{align}
		\frac{d}{dt}
		\bbm \mathbf{z} \\ h \ebm &= \bbm \calA_q & -\calC_{\gamma} \\ -\calC_{\theta_e} & \calB_q \ebm \bbm \mathbf{z} \\ h \ebm - \bbm  \calN_q  & 0  \\ 0 & \calM_q[\mathbf{z}] \ebm \bbm \mathbf{z} \\ h \ebm + \bbm P_q (m\mathbf{u}) \\ mv \ebm \\[2mm]
		\begin{picture}(0,0)
		\put(-20,14){ $\left\{\rule{0pt}{40pt}\right.$}\end{picture}
		\bbm \mathbf{z}(0) \\ h(0) \ebm &= \bbm \mathbf{z}_0 \\ h_0 \ebm \in \mathbf{W}^q_{\sigma}(\Omega)
		\end{align}
	\end{subequations}	
	
	\begin{multline}\label{1.18}
	\BA_q = \bbm \calA_q & -\calC_{\gamma} \\ -\calC_{\theta_e} & \calB_q \ebm : \mathbf{W}^q_{\sigma}(\Omega) = \lso \times \lqo \supset \calD(\BA_q) = \calD(\calA_q) \times \calD(\calB_q) \\ = (\mathbf{W}^{2,q}(\Omega) \cap \mathbf{W}^{1,q}_{0}(\Omega) \cap \lso) \times (W^{2,q}(\Omega) \cap W^{1,q}_{0}(\Omega)) \longrightarrow \mathbf{W}^q_{\sigma}(\Omega).
	\end{multline}	
	
	\noindent Properties of the operator $\ds \BA_q$ will be given in the next Section \ref{subsec-8}.
	
	\subsection{The linearized $\mathbf{w}$-problem of the translated $\mathbf{z}$-model}
	\noindent Next, still for $1 < q < \infty$, we introduce the linearized controlled system of the translated model (\ref{1.5}), or (\ref{1.16}), (\ref{1.17}) in the variable $\mathbf{w} = \{\mathbf{w}_f, w_h\} \in \lso \times \lqo \equiv \mathbf{W}^q_{\sigma}(\Omega)$:
	\begin{equation}\label{1.19}
	\frac{d}{dt} \bbm\mathbf{w}_f \\ w_h \ebm = \BA_q \bbm \mathbf{w}_f \\ w_h \ebm + \bbm  P_q(m\mathbf{u}) \\ mv \ebm \text{ in } \mathbf{W}^q_{\sigma}(\Omega),
	\end{equation}
	
	\noindent with I.C. $\{\mathbf{w}_f(0), w_h(0)\} \in \mathbf{W}^q_{\sigma}(\Omega) = \lso \times \lqo$. Its PDE version is, recalling \eqref{1.5}
	\begin{subequations}\label{N1.23}
		\begin{align}
		\frac{d}{dt}\mathbf{w}_f - \nu \Delta \mathbf{w}_f + L_e(\mathbf{w}_f) - \gamma w_h \mathbf{e}_d + \nabla \chi &= m \mathbf{u}   &\text{ in } Q \label{N1.23a}\\
		\frac{d}{dt}w_h - \kappa \Delta w_h + \mathbf{y}_e \cdot \nabla w_h + \mathbf{w}_f \cdot \nabla w_h + \mathbf{w}_f \cdot \nabla \theta_e &= mv &\text{ in } Q \label{N1.23b}\\
		\begin{picture}(0,0)
		\put(-200,10){ $\left\{\rule{0pt}{50pt}\right.$}\end{picture}
		\text{div } \mathbf{w}_f &= 0   &\text{ in } Q \label{N1.23c}\\
		\mathbf{w}_f \equiv 0, \ w_h &\equiv 0 &\text{ on } \Sigma \label{N1.23d}\\
		\mathbf{w}_f(0,\cdot) = \mathbf{w}_{f,0}; \quad w_h(0,\cdot) & = w_{h,0} &\text{ on } \Omega. \label{N1.23e}
		\end{align}
	\end{subequations}	
	
	\subsection{Properties of the operator $\BA_q$ in \eqref{1.18}}
	\label{subsec-8}
	We shall use throughout the following notation (recall \eqref{1.3}, (\ref{1.18}) and \eqref{N4-111})
	\begin{equation}\label{N1-63}
	\mathbf{W}^q_{\sigma}(\Omega) \equiv \lso \times \lqo \quad
	\mathbf{V}^{q,p}(\Omega) \equiv  \Bto \times \lqo.
	\end{equation}
	
	\begin{rmk}\label{remark1-1}
		By using the maximal regularity of the heat equation, instead of the state space $ \VqpO $ in (\ref{N4-111}), (\ref{N1-63}) we could take the state space  $\mathbf{V}_{b}^{q,p}$ to be the product of two Besov spaces, i.e. $\mathbf{V}_{b}^{q,p}(\Omega) \equiv  \Bt \times B^{ 2-2/p}_{q,p}( \Omega) $. Here, the second Besov component is the real interpolation between $\lqo$ and $\calD(\calB_q)$, see \cite{PSch2001}. This remark applies to all results involving $V^{q,p}$, but it will not necessarily be noted explicitly case by case, in order not to overload the notation.
	\end{rmk}
	
	\noindent Accordingly, we shall look at the operator $\BA_q$ in \eqref{1.18} as defined on either space
	\begin{equation}        \label{N1-64}
	\BA_q: \mathbf{W}^q_{\sigma}(\Omega) \supset \calD(\BA_q) \to \mathbf{W}^q_{\sigma}(\Omega) \qquad
	\mbox{or} \qquad
	\BA_q: \mathbf{V}^{q,p}(\Omega) \supset \calD(\BA_q) \to \mathbf{V}^{q,p}(\Omega).
	\end{equation}
	In Section \ref{new_sec_3}, we shall omit specification of $\Omega$ and simply write $\mathbf{W}^q_{\sigma}$. The following result collects basic properties of the operator $\BA_q$. It is essentially a corollary of Theorems \ref{A-Thm-1.4} and \ref{A-Thm-1.5} in Appendix \ref{app-A} for the Oseen operator $\calA_q$, as similar results hold for the operator $\calB_q$, while the operator $\calC_{\gamma}$ and $\calC_{\theta_e}$ in the definition \eqref{1.18} of $\BA_q$ are bounded operators, see \eqref{N1-17}, \eqref{N1-18}.
	
	\begin{thm}\label{I-Thm-1.7}
		With reference to the Operator $\BA_q$ in \eqref{1.18}, \eqref{N1-64}, the following properties hold true:
		\begin{enumerate}[(i)]
			\item $\ds \BA_q$ is the generator of strongly continuous analytic semigroup on either $\mathbf{W}^q_{\sigma}(\Omega)$ or $\mathbf{V}^{q,p}(\Omega)$ for $t > 0$;
			\item $\BA_q$ possesses the $L^p$-maximal regularity property on either $\mathbf{W}^q_{\sigma}(\Omega)$ or $\mathbf{V}^{q,p}(\Omega)$ over a finite interval:
			\begin{equation}
			\BA_q \in MReg_{p} (L^p(0,T;*)), \ 0 < T < \infty, \quad
			(*) = \mathbf{W}^q_{\sigma}(\Omega) \mbox{ or } \mathbf{V}^{q,p}(\Omega).
			\end{equation}
			\item $\ds \BA_q$ has compact resolvent on either $\mathbf{W}^q_{\sigma}(\Omega)$ or $\mathbf{V}^{q,p}(\Omega)$.
		\end{enumerate}
	\end{thm}
	
	\noindent Analyticity of $\ds e^{\calA_qt}$ (resp. $e^{\calB_qt}$) in $\lso$ (resp. $L^q(\Omega)$) implies analyticity of $\ds e^{\calA_qt}$ (resp. $\ds e^{\calB_qt}$) on $\ds \calD(\calA_q) = \calD(A_q)$ (resp. $\ds \calD(\calB_q) = \calD(B_q)$), hence analyticity of $\ds e^{\calA_qt}$ (resp. $\ds e^{\calB_qt}$) on the interpolation space $\ds \Bto$  in (\ref{A-1.16bb}). (or in $\ds \Bso$) in (\ref{A-1.151}) (resp. \eqref{I-1.3} in the scalar case).\\
	
	\noindent For the notation of, and the results on, maximal regularity, see \cite{HA:1995}, \cite{Dore:2000}, \cite{GGH:2012}, \cite{HS:2016}, \cite{KW:2001}, \cite{KW:2004}, \cite{PS:2016}, \cite{We:2001}, \cite{W:1985}, etc. In particular, we recall that on a Banach space, maximal regularity implies analyticity of the semigroup, \cite{DeS} but not conversely.We refer to Appendix \ref{app-A}.\\
	
	\noindent \textbf{Basic assumption:} By Theorem \ref{I-Thm-1.7}, the operator $ \BA_q $ in (\ref{1.18}) has the eigenvalues (spectrum) located in a triangular sector of well-known type. Then our basic assumption - which justifies the present paper - is that such operator $\ds \BA_q$ is unstable: that is $\BA_q$ has a finite number, say $N$, of eigenvalues $\lambda_1, \lambda_2 ,\lambda_3 ,\dots,\lambda_N$ on the complex half plane $\{ \lambda \in \mathbb{C} : Re~\lambda \geq 0 \}$ which we then order according to their real parts, so that	
	\begin{equation}\label{1.21}
	\ldots \leq Re~\lambda_{N+1} < 0 \leq Re~\lambda_N \leq \ldots \leq Re~\lambda_1,
	\end{equation}
	
	\noindent each $\lambda_i, \ i=1,\dots,N$, being an unstable eigenvalue repeated according to its geometric multiplicity $\ell_i$. Let $M$ denote the number of distinct unstable eigenvalues $\lambda_i$ of $\BA_q$, so that $\ell_i$ is equal to the dimension of the eigenspace corresponding to $\lambda_i$. Instead, $\ds N = \sum_{i = 1}^{M} N_i$ is the sum of the corresponding algebraic multiplicity $N_i$ of $\lambda_i$, where $N_i$ is the dimension of the corresponding generalized eigenspace.\\
	
	\begin{rmk} \label{newremark}
		This remark is inserted upon request of a referee. Condition \eqref{1.21} is intrinsic to the notion of `stabilization', whereby then one seeks to construct a feedback control that transforms an original unstable problem (with no control) into a stable one. However, as is well-known, the same entire procedure can be employed to \underline{enhance at will the stability} of an originally stable system ($ Re~\lambda_1 < 0 $) by feedback control. Regarding the issue of multiplicity of the eigenvalues, we recall the well-known classical result \cite{mi}, \cite{u1} that the eigenvalues of the Laplacian are simple generically with respect to the domain.  	
	\end{rmk}
	
	\noindent The ability to stabilize   the Boussinesq system \eqref{1.1} in the vicinity of the (unstable) equilibrium solution $\{\mathbf{y}_e, \theta_e\}$, by means of a finite dimensional feedback control pair $\{\mathbf{u},v\}$ localized on the arbitrarily small sub-domain $\omega$, depends (non-surprisingly, according to past experience \cite{BLT1:2006}, \cite{LT1:2015}, \cite{LT2:2015}, \cite{RT:2009}--\cite{RT2:2009}) on an appropriate   Unique Continuation Property. To this end, we quote two results. They are part of several UCP proved in \cite{TW.1}.
	\ifdefined\xxxx
	Approach I seems more natural, but it requires the new Unique Continuation Property \eqref{2.1} $\to$ \eqref{2.2} which appears to be unknown at present, except for special cases.
	
	Approach II manages to circumvent the aforementioned Unique Continuation Property and simply requires the known Unique Continuation Property \eqref{N1-69} $\to$ \eqref{N1-70}. The corresponding approach seems to be less direct and less elegant. I will be provided by a companion paper by the authors.
	\bigskip
	
	\noindent {\bf APPROACH I: Solution of the uniform stabilization of the Boussinesq system \eqref{1.1} subject to the validity of a Unique Continuation Property \eqref{2.1} $\to$ \eqref{2.2} below.}
	\bigskip
	
	\fi
	
	\begin{thm} \label{T1.3}
		Let $\ds \Phi = \bbm \boldsymbol{\varphi} \\ \psi \ebm \in [\mathbf{W}^{2,q}(\Omega) \cap \lso] \times W^{2,q}(\Omega), \ \pi \in W^{1,q}(\Omega)$, $ \boldsymbol{\varphi} = [\varphi_1, \ldots, \varphi_d], \psi = \mbox{scalar, satisfy the following overdetermined problem}$
		\begin{subequations}\label{2.111}
			\begin{align}
			- \nu \Delta \boldsymbol{\varphi} + L_e(\boldsymbol{\varphi}) + \nabla \pi - \gamma \psi \mathbf{e}_d&= \lambda \boldsymbol{\varphi}   &\text{in } \Omega      \label{2.11a} \\
			- \kappa \Delta \psi + \mathbf{y}_e \cdot \nabla \psi + \boldsymbol{\varphi} \cdot \nabla \theta_e &= \lambda \psi &\text{in } \Omega \label{2.11b}\\
			\begin{picture}(0,0)
			\put(-130,6){ $\left\{\rule{0pt}{35pt}\right.$}\end{picture}
			\div \ \boldsymbol{\varphi} &= 0   &\text{in } \Omega  \label{2.11c} \\
			\boldsymbol{\varphi} = 0, \ \psi &= 0 &\text{on } \omega.    \label{2.11e}
			\end{align}		
		\end{subequations}
		Then
		\begin{equation}\label{2.2}
		\boldsymbol{\varphi} = 0, \psi = 0, \pi = const \text{ in } \Omega.
		\end{equation}
		
	\end{thm}
	\noindent We point out that the B.C
	\begin{equation}\label{1.36}
	\boldsymbol{\varphi}|_{\Gamma} = 0, \ \psi|_{\Gamma} = 0 \text{ on } \Gamma
	\end{equation}
	are not needed. [If included, the resulting problem would be an eigenproblem for the Boussinesq operator with over/determination as in (\ref{2.11e})]. However, the UCP required in Section \ref{sec-4} to establish the controllability of the finite dimensional projected problem (\ref{1.27a}) via verification of the Kalman rank condition \eqref{2.8} involves the following adjoint problem.
	\begin{thm} \label{T1.4}
		\noindent {\bf (The required UCP).} In the same notation $ \{\boldsymbol{\varphi}, \psi, \pi \} $ of Theorem \ref{T1.3} with the same assumed regularity,let
		\begin{subequations}\label{2.1}
			\begin{align}
			- \nu \Delta 	\boldsymbol{\varphi} + L^{\ast}_e(\boldsymbol{\varphi}) + \nabla \pi + \psi \nabla \theta_e &= \lambda \boldsymbol{\varphi}   &\text{in } \Omega      \label{2.1a} \\
			- \kappa \Delta \psi + \mathbf{y}_e \cdot \nabla \psi - \gamma \boldsymbol{\varphi} \cdot  \mathbf{e}_d &= \lambda \psi &\text{in } \Omega \label{2.1b}\\
			\begin{picture}(0,0)
			\put(-125,4){ $\left\{\rule{0pt}{35pt}\right.$}\end{picture}
			\text{div } \boldsymbol{\varphi}&= 0   &\text{in } \Omega  \label{2.1c} \\
			\psi = 0 ,\,	\{\varphi_1,...,\varphi_{d-1} \} &= 0,  &\text{on } \omega.    \label{2.1e}
			\end{align}		
		\end{subequations}
		
		\noindent Then
		\begin{equation}\label{1.38}
		\boldsymbol{\varphi} = 0, \psi = 0, \pi = const \text{ in } \Omega.
		\end{equation}
	\end{thm}
	\begin{rmk} \label{remark1-21}
		Again, for problem (\ref{2.1}a-b-c), the B.Cs such as (\ref{1.36}) are \underline{not} needed. We note that in  (\ref{2.1a}) the definition of the adjoint of $ 	L_e(\mathbf{z}) = (\mathbf{y}_e \cdot \nabla)\mathbf{z} + (\mathbf{z} \cdot \nabla)\mathbf{y}_e$ in (\ref{1.5g}) is $ L^{\ast}_e(\boldsymbol{\varphi}) \equiv (\mathbf{y}_e \cdot \nabla) \boldsymbol{\varphi} + \nabla^{\perp} \mathbf{y}_e \cdot \boldsymbol{\varphi}$. We also note that in taking the adjoint of the operator $ \BA_q$ in (\ref{1.18}) the coupling operators  $\calC_{\gamma}$ and $\calC_{\theta_e}$ switch places into their adjoints. See $ \BA_q^*$ in (\ref{B11.18}) of Appendix \ref{app-B}. This fact is critical in obtaining in the adjoint UCP of Theorem \ref{T1.4} that only the first $ d-1 $ components $ \{\varphi_1,...,\varphi_{d-1} \}$ of $ \boldsymbol{\varphi}$ need to be assumed as vanishing in $ \omega $ in (\ref{2.1e}). See Appendix \ref{app-B} for more details. Notice, in fact, that if condition $\psi = 0 \text{ in } \omega $ is used in (\ref{2.1b}), it then follows that the last component $ \varphi_{d} = 0 $ in $ \omega$ because of the form of the vector $ \mathbf{e}=(0,0...,1) $. This along with (\ref {2.1e}) yield the same overdetermined conditions on $ \omega $ as in (\ref{2.11e}).  In this sense, Theorem \ref{T1.4}  is essentially a corollary of the proof of Theorem \ref{T1.3} and this justifies the insertion of Theorem \ref{T1.3} here. See \cite{TW.1} for proofs via Carleman type estimates. Here, moreover it is shown that for the case where the subset $ \omega $ is in particular subtended by an arbitrary portion of the boundary $ \Gamma $ and satisfies some additional geometrical conditions, then only $(d-2)$ components of the vector $ \boldsymbol{\varphi}$ need to be assumed to vanish on $ \omega $, thus further relaxing condition (\ref{2.1e}). In these cases, however we also need a boundary condition weaker than $ \boldsymbol{\varphi}|_{\Gamma} = 0 $. Details are given in \cite[Theorems 1.6 and 1.7]{TW.1}. These results for the adjoint Boussinesq static problem are in line with the open-loop controllability results in \cite{cg2}, \cite{Gur:2006}, \cite{C-G-I-P1} \cite{cl1}.
	\end{rmk}
	
	\begin{rmk} \label{remark1-2}
		References \cite{fabre}, \cite{fl1} provide a UCP for the Stokes problem with implications on approximate controllability. The UCP property  of Theorem \ref{T1.3}  in the case $\psi  \equiv 0 $ [i.e. Oseen's operator in $\boldsymbol{\varphi}$] has been shown via Carleman's estimates in \cite{BT:2004} by first transforming $\Omega$ in a ``bent" half-space with a parabolic boundary, next selecting the Melrose-Sjostrand form for the Laplacian, and finally applying the Carleman estimates in integral form from \cite{H}. A different proof, directly on $\Omega$, and this time with no use of the condition $\boldsymbol{\varphi}\equiv 0$ on $\Gamma$, was later given in \cite{RT2:2009}, also via use of (different) Carleman-type estimates for the Laplacian. The direct (on $ \Omega $) Carlemann-type estimate proof in \cite{RT2:2009}, with additional technical modification, can be made to work and establish both Theorem \ref{T1.3} and Theorem \ref{T1.4} (\cite{TW.1}). In effect, such reference \cite{TW.1} provides two additional UCPs, say for the original problem (\ref{2.111} a-b-c) with over determination this time on a small subdomain $ \omega $ subtended by an arbitrary portion of the boundary. A proof yielding, say Theorem \ref{T1.3} with $ d = 2 $ and with limited regularity of the solution was given in \cite{RRR:2019}. Theorem \ref{T1.4} is in line with an ``observability inequality" for the time dependent problem (\ref{1.1}) needed in the study of local controllability to the origin or to a trajectory given in \cite{fgip}. It improves on the prior observability inequality in \cite{Gur:2006}. We shall invoke Theorem \ref{T1.4} in Section \ref{sec-4}, to verify the Kalman rank condition \eqref{2.8}. The proof of the implication: Theorem \ref{T1.4} $ \implies $\eqref{2.8} is given Appendix \ref{app-B}. We can now state the main results of the present paper.\\
		
	\end{rmk}
	
	
	\ifdefined\xxxxx
	
	\begin{rmk}
		(i) If the thermal equilibrium solution $\theta_e \equiv \overline{\theta} = 0$, then the answer is: Yes. This is so since, with $\theta_e \equiv 0$, the thermal equation $ -\kappa \Delta \psi + y_e \cdot \nabla \psi = \lambda \psi$ in (\ref{2.1b}) along with $\psi = 0$ in $\omega$ in (\ref{2.1e}) implies $\psi \equiv 0$ in $\Omega$ by a standard result in elliptic theory. Next, the $\varphi$-equation \eqref{2.1a} in $\Omega$ is now uncoupled and overdetermined because of $\varphi \equiv 0$ in $\omega$. It is known that [\quad RT] it then implies $\varphi \equiv 0$ in $\Omega$, $\pi \equiv$ const in $\Omega$, according to the following Unique Continuation Property of the eigenvalue Oseen equation:
		\begin{thm}     \label{RT_thm}
			Let $\{\varphi, \pi\} \in W^{2,q}(\Omega) \times W^{1,q}(\Omega)$ satisfy
			\begin{subequations}\label{N1-69}
				\begin{align}
				- \nu \Delta \varphi + L_e(\varphi) + \nabla \pi &= \lambda \varphi   &\text{in } \Omega      \label{N1-69a} \\
				\begin{picture}(0,0)
				\put(-120,-5){ $\left\{\rule{0pt}{30pt}\right.$}\end{picture}
				\text{div } \varphi &= 0   &\text{in } \Omega\\
				\varphi  &= 0 &\text{on } \omega
				\end{align}		
			\end{subequations}
			Then
			\begin{equation}
			\label{N1-70}
			\varphi = 0,\qquad
			\pi = const \text{ in } \Omega
			\end{equation}
		\end{thm}	
		(ii) The case of point (i) arises for instance when the data of the steady state solution \eqref{1.2} satisfy
		\begin{equation}
		\label{N1-71}
		g(x) \equiv 0,\quad
		f(x) = \text{conservative} = \nabla F(x) \mbox{ for some } F(x), \quad
		\overline{\theta} \equiv 0.
		\end{equation}
		Then, under \eqref{N1-71}, the pair $\{y_e \equiv 0, \theta_e \equiv \overline{\theta} \equiv 0\}$ in $\Omega$ is a solution of system \eqref{1.2} along with $\pi = $ constant in $\Omega$.
		\hfill \qedsymbol
	\end{rmk}	\fi
	\noindent Recall that the vector $\{ \boldsymbol{\varphi}, \psi \} = \{ \varphi_1,\varphi_2,...,\varphi_d,\psi \}$ in $\lso $ has $( d+ 1 )$ coordinates, the first $ d$  co-ordinates correspond to the fluid space, while the last coordinate corresponds to the heat space. Motivated by the UCP of Theorem \ref{T1.4}, we shall introduce the space $  \widehat{ \mathbf{L}}^q_{\sigma} (\Omega) $
	\begin{multline}\label{1.39}
	\widehat{ \mathbf{L}}^q_{\sigma} (\Omega) \equiv \text{ the space obtained from } \mathbf{L}^q_{\sigma} (\Omega) \text{ after omitting the } \\ d\text{-coordinate from the vectors of }\mathbf{L}^q_{\sigma} (\Omega).
	\end{multline}

	\section{Main results}\label{new_sec_2}
	
	As in our past work \cite{BT:2004}, \cite{BLT1:2006}, \cite{LT2:2015}, \cite{LT3:2015}, \cite{LPT.1}, \cite{LPT.2}, we shall henceforth let $\lso$ denote the complexified space $\lso + i\lso$, and similarly for $\lqo$, whereby then we consider the extension of the linearized problem (\ref{1.19}) to such complexified space $\lso \times \lqo$. Thus, henceforth, $ \mathbf{w}$ will mean $ \mathbf{w} + i \tilde{\mathbf{w}}$, $\mathbf{u}$ will mean $\mathbf{u} + i \tilde{\mathbf{u}}$, $\mathbf{w}_0$ will mean $\mathbf{w}_0 + i \tilde{\mathbf{w}}_0$. Our results would be given in this complexified setting. How to return to real-valued formulation of the results was done in these past reference (see e.g. \cite{BT:2004}, \cite{BLT1:2006}, \cite{LT2:2015}, \cite[Section 2.7]{LPT.1}). Because of space constraints, such real-valued statements will not be explicitly listed on the present paper. We refer to the above references. A main additional feature of the results below is that the feedback control $ \mathbf{u}^1_k$ corresponding to the fluid equation is of \underline{reduced} dimension: that is, of dimension $(d-1)$ rather than of dimension $d$. This is due to the UCP of Theorem \ref{T1.4}, as noted in Remark \ref{remark1-21}.
	
	\subsection{Global well-posedness and uniform exponential stabilization of the linearized $\mathbf{w}$-problem \eqref{1.19} on either the space $ \mathbf{W}^q_{\sigma}(\Omega) \equiv \lso \times \lqo$ or the space $ \mathbf{V}^{q,p}(\Omega) \equiv  \Bt \times \lqo$, $1 < q < \infty, 1<p<\frac{2q}{2q-1}$.}
	\label{subsec-2.1}
	See also Remark \ref{remark1-1}.
	\begin{thm}\label{NThm-2.1}
		Let the operator $\BA_q$ in \eqref{1.18} have $N$ possibly repeated unstable eigenvalues $\ds \{ \lambda_j \}_{j = 1}^N$  as in (\ref{1.21}), of which $M$ are distinct. Let $\ell_i$ denote the geometric multiplicity of $\lambda_i$. Set $K = \sup \{ \ell_i; i = 1, \cdots, M \}$. Let $\ds (\mathbf{W}^q_{\sigma})^u_N$ be the $N$-dimensional subspace of $\mathbf{W}^q_{\sigma}(\Omega)$ defined in \eqref{1.23} below. Recall the space $ \widehat{ \mathbf{L}}^q_{\sigma} ( \Omega )$ from (\ref{1.39}) and let likewise $\ds (\widehat{\mathbf{W}}^q_{\sigma})^u_N$ be the space obtained from $ (\mathbf{W}^q_{\sigma})^u_N$ by omitting the $ d $ coordinate from the vectors of $\ds  (\mathbf{W}^q_{\sigma})^u_N $. Then, one may construct a finite dimensional feedback operator $F$, such that, with $\ds \mathbf{w}=\{ \mathbf{w}_f,w_h\}, \mathbf{w}_N = P_N\mathbf{w}$ and
		
		\begin{equation}\label{N-2.1}
		\bbm \mathbf{u} \\[2mm] v \ebm
		=
		\bbm F^1(\mathbf{w}) \\[2mm] F^2(\mathbf{w}) \ebm, \
		F(\mathbf{w}) = 	\bbm P_q(mF^1(\mathbf{w})) \\[2mm] mF^2(\mathbf{w}) \ebm
		= \bbm P_q \Big( m \big( \sum_{k = 1}^{K} (\mathbf{w}_N, \mathbf{p}_k)_{\omega} \mathbf{u}^1_k \big) \Big) \\[2mm] m \big( \sum_{k = 1}^{K} (\mathbf{w}_N, \mathbf{p}_k)_{\omega} u^2_k \big)\ebm,
		\end{equation}
		\noindent with vectors $[\mathbf{u}^1_k, u^2_k] \in (\widehat{\mathbf{W}}^q_{\sigma})^u_N \subset \widehat{\mathbf{L}}^{q}_{\sigma}(\Omega) \times L^q(\Omega)$, and $\mathbf{p}_k \in ((\mathbf{W}^q_{\sigma})^u_N)^* \subset \lo{q'} \times L^{q'}(\Omega)$ (Remark \ref{rmkB.1}), then the $\mathbf{w}$-problem (\ref{1.19}) can be rewritten in feedback form on $\mathbf{W}^q_{\sigma}( \Omega )  $ as follows	
		\begin{multline}\label{2.9}
		\frac{d \mathbf{w}}{dt} = \frac{d}{dt} \bbm \mathbf{w}_f \\[2mm] w_h \ebm = \BA_q \mathbf{w} + F(\mathbf{w}) = \BA_q \mathbf{w} + \bbm P_q \Big( m \big( \sum_{k = 1}^{K} (\mathbf{w}_N, \mathbf{p}_k)_{\omega} \mathbf{u}^1_k \big) \Big) \\[2mm] m \big( \sum_{k = 1}^{K} (\mathbf{w}_N, \mathbf{p}_k)_{\omega} u^2_k \big)\ebm = \BA_{F,q} \mathbf{w}; \\ \mathbf{w}(0) = \mathbf{w}_0,
		\end{multline}
		
		with $\ds \calD(\BA_{F,q}) = \calD(\BA_q)$ where the operator $\ds \BA_{F,q}$ in (\ref{2.9}) has the following properties:
		
		\begin{enumerate}[(i)]
			\item It is the generator of a s.c. analytic semigroup $\ds e^{\BA_{F,q}t}$ in the space $ \mathbf{W}^q_{\sigma}(\Omega) \equiv \lso \times \lqo$ as well as in the space $\mathbf{V}^{q,p}(\Omega) \equiv  \Bt \times \lqo$, see also Remark \ref{remark1-1}.
			\item It is uniformly (exponentially) stable in  either of these spaces
			\begin{equation}\label{2.10}
			\norm{e^{\BA_{F,q}t} \mathbf{w}_0}_{(\cdot)} \leq C_{\gamma_0} e^{-\gamma_0 t}\norm{\mathbf{w}_0}_{(\cdot)}
			\end{equation}
			\noindent where $(\cdot)$ denotes either $\ds \lso \times \lqo \equiv \SqsO$ or else $\ds \Bto \times \lqo \equiv \VqpO$. In \eqref{2.10}, $\gamma_0$ is any positive number such that $Re~\lambda_{N+1} < -\gamma_0 < 0$.
			
			\item Finally, $\BA_{F,q}$ has maximal $L^p$-regularity up to $T = \infty$ on either of these spaces:
			\begin{equation}\label{2.11}
			\BA_{F,q} \in MReg_{p} (L^p(0,\infty; \ \cdot \ )), \text{ where } (\cdot) \text{ denotes}
			\begin{picture}(0,0)
			\put(-250,-10){ $\left\{\rule{0pt}{15pt}\right.$}\end{picture}
			\end{equation}
			\hspace{1cm} either $\ds \lso \times \lqo \equiv \SqsO$ or else $\ds \Bto \times \lqo \equiv \VqpO$.
			\hfill \qedsymbol
		\end{enumerate}
	\end{thm}
	\begin{rmk}
		
		We note explicitly that the vector $ \mathbf{u}^1_k $ acting on the fluid $ d-$dimensional component $ \mathbf{u}$ is of reduced dimension $(d-1)$, rather than $ d $.\end{rmk}
	
	\noindent The proof of Theorem \ref{NThm-2.1} begins in Section \ref{new_sec_3} and proceed through Section \ref{sec-4}.
	
	\subsection{Local well-posedness and uniform (exponential) null stabilization of the translated nonlinear $\{\mathbf{z}, h\}$-problem (\ref{1.5}) or (\ref{1.16}) by means of a finite dimensional explicit, spectral based feedback control pair localized on $\omega$} \label{new-subsec-2.2}
	
	Starting with the present subsection, the nonlinearity of problem (\ref{1.1}) will impose for $d = 3$ the requirement $q > 3$, see (\ref{3.24}) below. As our deliberate goal is to obtain the stabilization result for the fluid component $\mathbf{y}$ in the space $\ds \Bto$ \underline{which does not recognize boundary conditions}, Remark \ref{remark1-11}, then the limitation $\ds p < \rfrac{2q}{2q - 1}$ of this space applies, see \eqref{A-1.16bb}. In conclusion, our well-posedness and stabilization results will hold under the restriction $\ds q > 3, 1 < p < \rfrac{6}{5}$ for $d = 3$, and $\ds q > 2, 1 < p < \rfrac{4}{3}$ for $d = 2$.
	
	\begin{thm}\label{N-Thm-2.2}
		Let $\ds d = 2,3, \ q > d, \ 1 < p < \frac{2q}{2q-1}$. Consider the nonlinear $\{\mathbf{z},h\}$-problem \eqref{1.16} in the following feedback form
		
		\begin{equation}\label{N2.5}
		\frac{d}{dt} \bbm \mathbf{z} \\[2mm] h \ebm = \BA_{F,q} \bbm \mathbf{z} \\[2mm] h \ebm - \bbm \calN_q & 0 \\[2mm] 0 & \calM_q[\mathbf{z}] \ebm \bbm \mathbf{z} \\[2mm] h \ebm; \quad \BA_{F,q} = \bbm \calA_q & -\calC_{\gamma} \\[2mm] -\calC_{\theta_e} & \calB_q \ebm + F= \BA_q + F
		\end{equation}
		
		\noindent see \eqref{2.9}. Specifically,
		
		\begin{subequations}\label{N2.6}
			\begin{align}
			\frac{d\mathbf{z}}{dt} - \calA_q \mathbf{z} + \calC_{\gamma} h + \calN_q \mathbf{z} &= P_q \Bigg( m \bigg( \sum_{k = 1}^{K} \bigg(P_N \bbm \mathbf{z} \\ h \ebm, \mathbf{p}_k \bigg)_{\omega} \mathbf{u}^1_k \bigg) \Bigg) \label{N2.6a}\\
			\begin{picture}(0,0)
			\put(-20,20){ $\left\{\rule{0pt}{40pt}\right.$}\end{picture}
			\frac{dh}{dt} - \calB_q h + \calC_{\theta_e} \mathbf{z} + \calM_q[\mathbf{z}]h &= m \bigg( \sum_{k = 1}^{K} \bigg(P_N \bbm \mathbf{z} \\ h \ebm, \mathbf{p}_k \bigg)_{\omega} u^2_k \bigg) \label{N2.6b}
			\end{align}
		\end{subequations}
		i.e. subject to a feedback control of the same structure as in the linear $\mathbf{w}=\{\mathbf{w}_f, w_h\}$-dynamics \eqref{2.9}. Here, $\mathbf{p}_k, \mathbf{u}_k = [\mathbf{u}_k^1, u_k^2]$ are the same vectors as those constructed in Theorem \ref{NThm-2.1} and appearing in \eqref{2.9}.In particular $ \mathbf{u}^1_k $ has $(d-1)$ components. There exists a positive constant $r_1 > 0$ (identified in the proof below in (\ref{3.31})) such that if
		\begin{equation}\label{N3.11}
		\norm{\{ \mathbf{z}_0, h_0 \}}_{\Bto \times \lqo} < r_1,
		\end{equation}
		then problem \eqref{N2.6} defines a unique (fixed point) non-linear semigroup solution in the space
		\begin{equation}\label{N3.12}
		\{\mathbf{z},h\} \in \xipqs \times \xipq \equiv L^p(0,\infty; \calD(\BA_{F,q}) = \calD(\calA_q) \times \calD(\calB_q)) \cap W^{1,p}(0,\infty;\mathbf{W}^q_{\sigma}(\Omega)),
		\end{equation}
	\end{thm}
	where we have set, via \eqref{1.9.new}, \eqref{N1.9} for $ \calA_q $ and \eqref{E1.20}, \eqref{E1.15} for $\calB_q $:
	\begin{align} \label{E2.9}
	\mathbf{z} \in \xipqs & \equiv L^p(0,\infty; \calD(\calA_{q}) \cap W^{1,p}(0,\infty;\mathbf{L}^q_{\sigma}(\Omega)) \\[3mm] \nonumber
	& \equiv L^p(0,\infty; \mathbf{W}^{2,q}(\Omega) \cap \mathbf{W}^{1,q}_0 (\Omega )\cap \mathbf{L}^q_{\sigma}(\Omega))\cap W^{1,p}(0,\infty;\mathbf{L}^q_{\sigma}(\Omega))
	\end{align}
	\begin{align} \label{E2.10}
	h  \in \xipq & \equiv L^p(0,\infty; \calD(\calB_q) \cap W^{1,p}(0,\infty;L^q(\Omega)) \\[3mm] \nonumber
	& \equiv L^p(0,\infty; W^{2,q}(\Omega) \cap W^{1,q}_0 (\Omega ))\cap W^{1,p}(0,\infty;L^q(\Omega)).
	\end{align}
	Moreover, we have
	\begin{align}
	\xipqs & \hookrightarrow C([0, \infty ];\Bso) \label{E2.111}\\[3mm]
	\xipq & \hookrightarrow  C([0, \infty ];L^q(\Omega))\label{E2.11}
	\end{align}
	so that
	\begin{equation}\label{E2.13}
	\xipqs \times \xipq \hookrightarrow C([0, \infty ];\mathbf{V}^{q,p}(\Omega)),
	\end{equation}
	recalling the  embedding,  called trace theorem \cite[Theorem 4.10.2, p 180, BUC for $T=\infty$]{HA:2000}, \cite{PS:2016}, as in \cite[Eq. (1.30)]{LPT.1}. See also \eqref{A-1.30} in Appendix \ref{app-A}.\\
	
	\noindent The space $ \xipqs \times  \xipq $  defined above is the space of $L^p$-maximal regularity for the generator $\BA_{F,q}$. $\xipqs$ is the space of $L^p$-maximal regularity of the Stokes operator, $ A_q$ in \eqref{N1.9}, see \eqref{A-1.28} and also \eqref{A-1.30} in Appendix \ref{app-A}. $\xipq$ is the space of $L^p$-maximal regularity of the generator $\calB_q$ in \eqref{E1.15}. The proof of Theorem \ref{N-Thm-2.2} is given in Section \ref{new-sec-5}.
	
	\begin{thm}
		\label{NThm-2.3}
		In the situation of Theorem \ref{N-Thm-2.2}, we have that such solution is uniformly stable in the space
		$\mathbf{V}^{q,p}(\Omega) \equiv  \Bto \times \lqo$, see also Remark \ref{remark1-1}: there exists $\widetilde{\gamma}>0, M_{\widetilde{\gamma}} > 0$ such that said solution satisfies
		\begin{equation}
		\label{N2-9}
		\norm{\bbm \mathbf{z} \\ h \ebm(t)}_{\mathbf{V}^{q,p}(\Omega)} \leq M_{\widetilde{\gamma}} e^{-\widetilde{\gamma} t}         \norm{\bbm \mathbf{z}_0 \\ h_0 \ebm}_{\mathbf{V}^{q,p}(\Omega)}, \quad
		t \geq 0.
		\end{equation}
		
	\end{thm}
	\noindent The proof of Theorem \ref{NThm-2.3} is given in Section \ref{new-sec-6}.
	\subsection{Local well-posedness and uniform (exponential) stabilization of the original nonlinear $\{\mathbf{y},\theta\}$-problem (\ref{1.1}) in a neighborhood of an unstable equilibrium solution $\{\mathbf{y}_e, \theta_e\}$, by means of a finite dimensional explicit, spectral based feedback control pair localized on $\omega$}
	\label{new-subsec-2.3}
	
	\noindent The result of this subsection is an immediate corollary of sub-section \ref{new-subsec-2.2}.
	
	\begin{thm}\label{N-Thm-2.5}
		Let $1 < p < \rfrac{6}{5}, q > 3, d = 3$; and $1 < p < \rfrac{4}{3}, q > 2, d = 2$. Consider the original Boussinesq problem (\ref{1.1}). Let $\{\mathbf{y}_e,\theta_e\}$ be a given unstable equilibrium solution pair as guaranteed by Theorem \ref{Thm-1.1} for the steady state problem (\ref{1.2}). For a constant $\rho > 0$, let the initial condition $\{\mathbf{y}_0,\theta_0\}$ in (\ref{1.1e}) be in $\mathbf{V}^{q,p}(\Omega) \equiv  \Bto \times \lqo$, see also Remark \ref{remark1-1}, and satisfy
		\begin{equation}\label{N2.10}
		\boldsymbol{\calV}_{\rho} \equiv \Big\{ \{\mathbf{y}_0,\theta_0\} \in \mathbf{V}^{q,p}(\Omega): \norm{\mathbf{y}_0 - \mathbf{y}_e}_{\Bto} + \norm{\theta_0 - \theta_e}_{L^q(\Omega)} \leq \rho \Big\}, \quad \rho > 0.
		\end{equation}
		\noindent If $\rho > 0$ is sufficiently small, then
		\begin{enumerate}[(i)]
			\item for each $\{\mathbf{y}_0,\theta_0\} \in \boldsymbol{\calV}_{\rho}$, there exists an interior finite dimensional feedback control pair \label{I-Thm-2.5.i}
			\begin{equation}\label{I-2.26}
			\bbm \mathbf{u} \\[6mm] v \ebm
			=
			\bbm F^1\left(\bbm \mathbf{y}-\mathbf{y}_e \\ \theta-\theta_e \ebm\right) \\[3mm] F^2\left(\bbm \mathbf{y}-\mathbf{y}_e \\ \theta-\theta_e \ebm\right) \ebm
			=
			F\left(\bbm \mathbf{y}-\mathbf{y}_e \\ \theta-\theta_e \ebm\right)
			=
			\sum_{k = 1}^{K} \left( P_N \bbm \mathbf{y}-\mathbf{y}_e \\ \theta-\theta_e \ebm, \mathbf{p}_k\right)_{\omega} \mathbf{u}_k
			\end{equation}
			\noindent that is, of the same structure as in the translated $\{\mathbf{z},h\}$-problem (\ref{N2.6}), with the same vectors $\ds p_k, \mathbf{u}_k$ in (\ref{2.9}), such that the closed loop problem corresponding to (\ref{1.1})
			\begin{subequations}\label{I-2.27}
				\begin{align}
				\mathbf{y}_t - \nu \Delta \mathbf{y} +  (\mathbf{y} \cdot \nabla)\mathbf{y} - \gamma (\theta - \bar{\theta}) \mathbf{e}_d + \nabla \pi & = m\left(F^1\left(\bbm \mathbf{y}-\mathbf{y}_e \\ \theta-\theta_e \ebm\right)\right) + \mathbf{f}(x)  \text{ in } Q \label{I-2.27a}\\
				\theta_t - \kappa \Delta \theta + \mathbf{y} \cdot \nabla \theta &= m\left(F^2\left(\bbm \mathbf{y}-\mathbf{y}_e \\
				\theta-\theta_e \ebm\right)\right) + g(x)  \text{ in } Q   \label{new-I-2.27b}\\
				\div \ \mathbf{y} &= 0  \text{ in } Q \label{I-2.27b}\\
				\begin{picture}(125,-30)\put(-20,40){$\left\{\rule{0pt}{60pt}\right.$}\end{picture}
				\mathbf{y} =0, \quad \theta &= 0 \text{ on } \Sigma \label{I-2.27c}\\
				\mathbf{y}|_{t = 0} = \mathbf{y}_0, \quad \theta|_{t = 0} &= \theta_0 \text{ in } \Omega, \label{I-2.27d}
				\end{align}
			\end{subequations}
			\noindent rewritten abstractly after application of the Helmholtz projection $P_q$ as
			\begin{subequations}\label{I-2.28}
				\begin{align}
				\mathbf{y}_t + \nu A_q \mathbf{y} +  \calN_q \mathbf{y} + C_{\gamma} (\theta-\bar{\theta}) &= P_q \left[ m\left(F^1\left(\bbm \mathbf{y}-\mathbf{y}_e \\ \theta-\theta_e \ebm\right)\right) + \mathbf{f}(x) \right] \label{I-2.28a}\\
				\begin{picture}(125,0)	\put(-20,-15){$\left\{\rule{0pt}{60pt}\right.$}\end{picture}
				&= P_q \left[ m \left( \sum_{k=1}^{K} \left(P_N\bbm \mathbf{y}-\mathbf{y}_e \\ \theta-\theta_e \ebm, \mathbf{p}_k \right)_{\omega} \mathbf{u}^1_k \right) + \mathbf{f}(x) \right]\label{I-2.28b}\\
				\theta_t - \kappa B_q \theta + \calM_q[\mathbf{y}]\theta &= m \left( \sum_{k=1}^{K} \left(P_N\bbm \mathbf{y}-\mathbf{y}_e \\ \theta-\theta_e \ebm, \mathbf{p}_k \right)_{\omega} u^2_k \right) + g(x)   \label{new-I-2.28b}\\
				\{\mathbf{y}(0) = \mathbf{y}_0,\ \theta(0) = \theta_0 \} &\in  \Bto \times \lqo  \equiv \mathbf{V}^{q,p}(\Omega) \label{I-2.28d}
				\end{align}
			\end{subequations}
			\noindent has a unique solution $\{\mathbf{y},\theta\} \in C \big([0,\infty); \mathbf{V}^{q,p}(\Omega) \equiv  \Bto \times \lqo\big).$ See also Remark \ref{remark1-1}.
			Thus, here as in the proceeding Theorem \ref{NThm-2.1} through Theorem \ref{NThm-2.3}, we have
			\begin{align} \label{E2.19}
			& K = \sup \{ \ell_i; i = 1, \cdots, M \}, \mathbf{u}_k=[\mathbf{u}^1_k, u^2_k] \\[3mm]
			\begin{picture}(120,0)	\put(90,15){$\left\{\rule{0pt}{25pt}\right.$}\end{picture}
			& \mathbf{u}^1_k=[ u^{(1)}_k, u^{(2)}_k,...,u^{(d-1)}_k]. \nonumber
			\end{align}
			
			\item \label{thm2.4-ii}
			Moreover, such solution exponentially stabilizes the equilibrium solution $\{\mathbf{y}_e, \theta_e\}$ in the space $\Bto \times \lqo  \equiv \mathbf{V}^{q,p}(\Omega)$: there exist constants $\widetilde{\gamma} > 0$ and $M_{\widetilde{\gamma}} \geq 1$ such that said solution satisfies
			\begin{multline}\label{I-2.29}
			\norm{\mathbf{y}(t) - \mathbf{y}_e}_{\Bto} + \norm{\theta(t) - \theta_e}_{\lqo}
			\leq
			M_{\widetilde{\gamma}} e^{- \widetilde{\gamma} t} \left(\norm{\mathbf{y}_0 - \mathbf{y}_e}_{\Bto} + \right.\\ \left. \norm{\theta_0 - \theta_e}_{\lqo}\right),
			\end{multline}
			$ t \geq 0, \{\mathbf{y}_0,\theta_0\} \in \boldsymbol{\calV}_{\rho}$.
			\noindent Once the neighborhood $\boldsymbol{\calV}_{\rho}$ is obtained to ensure the well-posedness, then the values of $M_{\widetilde{\gamma}}$ and  $\widetilde{\gamma}$ do not depend on $\calV_{\rho}$ and $\widetilde{\gamma}$ can be made arbitrarily large through a suitable selection of the feedback operator $F$.
		\end{enumerate}	
	\end{thm}
	\noindent See Remark \ref{I-Rmk-9.1} comparing $\wti{\gamma}$ in (\ref{I-2.29}) with $\gamma_0$ in (\ref{2.10}).
	
	\subsection{Comparison with the literature}	\label{subsec_1.4}
	
	\noindent $\mathbf{1}$. With reference to both the ``Motivation" of Subsection \ref{subsec_1.1} as well as Subsection \ref{subsec-3}, it was already emphasized that all prior literature on the problem of feedback stabilization of either the Navier-Stokes equations or, subsequently, the Boussinesq system is carried out in a Hilbert-Sobolev setting. As already noted, this treatment is inadequate to obtain \underline{finite dimensionality} in full generality of the localized \underline{tangential boundary } feedback control for the 3d-Navier-Stokes equations. This obstacle then motivated the introduction of the $\mathbf{L}^q$-Sobolev/Besov setting in \cite{LPT.2}  with tight indices, see (\ref{I-1.3}), (\ref{A-1.16bb}), that does not recognize boundary conditions, see Remark \ref{remark1-11} , in order to solve affirmatively such open problem on the \underline{finite dimensionality} in 3d-Navier-Stokes \underline{tangential boundary feedback} stabilization. Reference  \cite{LPT.1} on localized interior controls sets the preparatory stage for the more demanding boundary control case in \cite{LPT.2}.\\
	
	\noindent $\mathbf{2}$. In the present contribution, following \cite{LPT.1,LPT.2}, the analysis of the original non-linear problem is carried out in the context of the property of \uline{maximal regularity of the linearized feedback} $\mathbf{w}$-problem (\ref{2.9}) or (\ref{N4-7}); that is, of the operator  $\BA_{F,q}$ in (\ref{N4-7}).  This applies to both well-posedness (in Section \ref{new-sec-5}) as well as stabilization (in Section \ref{new-sec-6}) of the nonlinear feedback problem. This is in contrast with prior treatments, such as \cite{BT:2004}, \cite{Barbu}, which often rely on chopping the nonlinearity and carrying out a limit process. The \uline{maximal regularity approach} introduced for stabilization problems in \cite{LPT.1} is cleaner and more desirable both technically and conceptually. On the other hand, applicability of Maximal Regularity requires well balanced spaces. Recent developments in this particular area \cite{We:2001,SS:2007,S:2006,PS:2016,KW:2001,KW:2004}, etc. were critical for carrying out our analysis. \\
	
	\noindent $\mathbf{3}$. We also point out that, as in \cite{LT1:2015,LT2:2015} and \cite{LPT.1,LPT.2} for the Navier-Stokes equations, in the case of the Boussinesq system, the number of needed controls will be related to the more desirable \underline{geometric} multiplicity, not the larger \underline{algebraic} multiplicity as in prior treatments, \cite{BT:2004}, \cite[p 276]{Barbu} of the unstable eigenvalues, by using the classical test for controllability of a system in Jordan form \cite{CTC:1984}, \cite[p 464]{B-M1}. This allows one also to obtain constructively an explicit form of the finite dimensional feedback control; and, moreover, to show that the \uline{feedback control acting on the fluid may be taken of reduced dimension}: one unit less, i.e. $d-1$, than the fluid component of dimension $d$. This is due to the Unique Continuation Property of the adjoint problem, as stated in Theorem \ref{T1.4}, Remark \ref{remark1-21}. This an additional contribution of the present paper. This dimension reduction in the closed-loop feedback stabilizing fluid control is in line with the open-loop controllability results in \cite{cg2}, \cite{Gur:2006}, \cite{C-G-I-P1} \cite{cl1}. The authors thank a referee for suggesting to investigate this possibility, which turned out to be successful. Such investigation also led to paper \cite{TW.1}.\\
	
	\noindent $\mathbf{4}$. We now review the literature on the feedback stabilization of the Boussinesq system. The first contribution is due to \cite{Wang} via internal feedback controls in the Hilbert setting $\mathbf{H}\times L^2(\Omega$), where $\mathbf{H}$ is defined in (\ref{A-1.21b}). \uline{The feedback controllers are both infinite-dimensional} of the type: $\mathbf{u}= -k(\mathbf{y}-\mathbf{y}_e)$ and $ v= -k(\theta-\theta_e),$ for large $k$, under technical assumptions on the localization of the controls. Paper \cite{Lef} also studies the feedback stabilization problem with, eventually, localized controls, which again are \uline{infinite dimensional}; eg defined in terms of a sub-differential. These results are in stark contrast with the present work, where - as noted in point 3 above - first, we establish that the stabilizing control is  finite dimensional; and, second, we show that its fluid component is of reduced dimension, as expressed only by means of $(d-1)$-components. Finally paper \cite{RRR:2019} studies the feedback stabilization problem with, this time, mixed boundary and claims to provide a rigorous treatment of problems studied in \cite{BH:2013}. \\
	
	\noindent$\mathbf{5}$. We now comment on the additional technical and conceptual difficulties in extending the uniform stabilization problem of  \cite{LPT.1} for the Navier-Stokes equations to the present Boussinesq system. First, the general strategy to attack the problem of feedback stabilization for \underline{parabolic-like dynamics} with feedback controllers of any type was introduced in \cite{RT:1975}. It has since become a standard approach in the literature, which was followed in a variety of parabolic problems \cite{LT1:1983}, \cite{LT22.1983}. In its implementation, however, this general strategy encounters an array of technical challenges that are highly feedback-dependent: the most demanding cases are when the feedback controllers are based on the boundary and even with boundary actuator, subject perhaps to additional requirements (arbitrarily small portion of the boundary, etc). This general strategy is used in \cite{LPT.1} as well as in the present paper, with the space $ \Bto $ in (\ref{A-1.16bb}) being the space of the fluid component where the stabilization is established. Such Besov space with tight indices has the critical property of being  'under the radar of compatibility conditions'. The presence of fluid-thermal coupling produces two level of problems: at the finite dimensional level, based an a UCP with special features; at the infinite dimensional level, based on establishing maximal regularity on a non-Hilbert setting, therefore stronger than analyticity of the underlying semigroup. As described below in point 6, the mathematical machinery used in the present paper could handle an array of other settings: first, a variety of different Boundary Conditions rather than just the non-slip Dirichlet boundary conditions \eqref{1.1d}; second, much more demanding 'seriously unbounded' coupling operators between the fluid and the thermal equations than the bounded coupling operators - $ \calC_{\gamma}$ and $\calC_{\theta_e}$ in (\ref{N1-17}), (\ref{N1-18}) - which are offered by the physically based Boussinesq model (\ref{1.1}).
	A key feedback - and problem - dependent source of difficulty encountered in the present Boussinesq model over the Navier-Stokes case of \cite{LPT.1} arises, as usual, at the level of testing and establishing controllability of the  finite dimensional projected $\mathbf{w}_N$- equation in (\ref{1.27a}); that is, in verifying the full rank condition (\ref{bbN4-5}). It took some time in the literature to realize that in the cases of "challenging" feedbacks such as the localized feedback of \cite{BT:2004} for the Navier-Stokes equation or else boundary feedback (with sensors and/actuators imposed as being boundary traces), a natural way of checking the resulting Kalman rank condition is to fall into an appropriate Unique Continuation Property for an adjoint problem. In doing so, as noted in Remark \ref{remark1-21}, we extract the benefit of obtaining the closed-loop feedback control acting on the fluid equation of reduced dimension: $(d-1)$ rather than $d$. As noted in point 3 above, this conclusion is in line with the open-loop controllability results in \cite{cg2}, \cite{Gur:2006}, \cite{C-G-I-P1}, \cite{cl1}. It was a referee who suggested that we carry out such a study. \\
	
	\noindent $\mathbf{6}$. In response to a referee's suggestion, we address the following \underline{Question}:
	Can our mathematical treatment cover other boundary conditions and more general coupled terms between the fluid and the thermal equations? \\
	
	\underline{The answer is in the affirmative}. We may e.g. cover all the B.C.s listed in the 2016 Birkh\"{a}user treatise \cite{PS:2016}, page 338. They refer to a very general Stokes problem (Laplacian replaced by normally strongly elliptic differential operator on a bounded domain,  in fact even with a domain with compact boundary, say of class $ \mathcal{C}^3 $). They are referred to as (mathematical definitions are given in \cite{PS:2016} page 338): (i) no-slip (the ones in our present treatment); (ii) pure slip; (iii) outflow; (iv) free. We now justify our assertion. Whatever the boundary conditions and the coupling terms, one ends up with the $ 2 \times 2 $ operator matrix of the type of $\BA_{F,q}$ in (\ref{N2.5}): the two operators $\calA_q$ and $\calB_q$ on the main diagonal are the free dynamic operators of each equation, while the two operators $\calC_{\gamma}$ and $ \calC_{\theta_e}$ on the secondary diagonal are the coupling terms. What does our mathematical treatment need? What it needs is that the resulting $ 2 \times 2 $ operator matrix possess the $ L^p$ maximal regularity on the desired functional setting [(fluid space) $\times$ (heat space)]; eventually, $\VqpO \equiv \Bto \times \lqo $ in (\ref{N4-111}). Hence, a-fortiori, be the generator of a s.c. analytic semigroup on such setting. Thus, as a first point, we need the statement of Theorem \ref{Thm-2.2}. In the present setting of problem (\ref{1.1}), the conclusion of Theorem \ref{Thm-2.2} is obtained in the basis of the following reasons:\\
	
	1a) The operators $\calA_q $ and $ \calB_q$ on the main diagonal have $L^p$ maximal regularity on each appropriate functional setting, on the basis of, ultimately, Solonnikov's old result for the Stokes problem for the fluid and also the corresponding known property for the heat operator. Thus the corresponding \underline{diagonal} $ 2 \times 2$ operator matrix has such $L^p$ maximal regularity on the desired cross functional setting. \\
	
	1b) Moreover, in the present physical model (\ref{1.1}), the two coupling operators  $ \calC_{\gamma} $ and $ \calC_{\theta_e}$ on the secondary diagonal are \underline{bounded} operators, see   (\ref{N1-17}), (\ref{N1-18}), But \underline{boundedness} of these two coupling operators is \underline{much more} than it is needed mathematically. \\
	
	\noindent \textbf{Conclusion}: We can take other homogeneous boundary conditions for the fluid and the heat equations \uline{as long as they produce the required $ L^p$ maximal regularity properties}, on the space $  \VqpO \equiv \Bto \times \lqo $ in (\ref{N4-111}). Regarding the fluid equation, each of the above B.C.s (i), (ii), (iii), (iv) does yield $ L^p$ maximal regularity of the Stokes problem and even for much more general elliptic operators than the Laplacian. See \cite[Theorem 7.3.1, p339]{PS:2016}.The Stokes operator is enough as the  $ L^p$ maximal regularity can then be transferred to the Oseen operator. As to the thermal equation, other B.C.s such as Neumann or Robin, beside the Dirichlet B.C., will likewise yield $ L^p$ maximal regularity for the thermal Laplacian. Moreover, the coupling operators $ \calC_{\gamma} $ and $ \calC_{\theta_e}$ may be much worse than \underline{bounded} perturbations as in (\ref{N1-17}), (\ref{N1-18}).It is well known \cite{Dore:2000},  \cite{KW:2001} that if an operator enjoys $ L^p$ maximal regularity, then adding a perturbation up to the same level minus epsilon preserve $ L^p$ maximal regularity (same as for analyticity). This will allow to greatly extend the model, purely in mathematical grounds, not necessarily based on physical grounds. Our model \ref{1.1} is the one which we found in many PDE-papers, beyond the control/stabilization areas.\\
	
	\noindent $\mathbf{7.}$  \underline{Question:} With boundary conditions leading to global controllability, is it possible - a referee asked - to get global stabilization with controls as in this paper? It is not easy to answer this question in full generality. First, there are several works on local exact controllability to the origin or to trajectories, see \cite{C-G-I-P1}, \cite{C-1}, \cite{C-2}, \cite{cg1}, \cite{cg2}, \cite{fgip}, \cite{Fernandezcara}, \cite{F}, \cite{Gur:2006} for an incomplete sample of recent works, using radically different techniques from those of the present paper (save for the underlying UCP of Theorem \ref{T1.4}, such as Carleman type inequalities for both Navier-Stokes or Boussinesq equations. On the other hand, the issue of deducing stabilization from controllability in the case of finite dimensional non linear ODE was much investigated in the 70-80's. As is well known, there is  a revealing, simple counter example of Brockett: a system of dimension $ n = 3 $ with $ m = 2 $ controls that is locally controllable, (though its linearization is not controllable), symmetric, yet it violates the condition of asymptotic stability \cite{Br}, \cite[p 88]{Za}. We also point out that there is extensive work at the PDE-level (or infinite-dimensional level) which shows that controllability implies stabilization under an additional condition (which is often true in many PDE-problems but it is challenging to check) \cite{L-T.4}, \cite{L-T.5}, \cite{L-T.6}, \cite{L-T.7}. However these works refer to hyperbolic or hyperbolic-like dynamics (such as wave equations, dynamic system of elasticity, Schr\"{o}dinger equations or plate equations).	
	
	\section{ Beginning with the proof of Theorem \ref{NThm-2.1}, Spectral decomposition of the linearized $\mathbf{w}$-problem \eqref{1.19} or \eqref{N1.23}}\label{new_sec_3}
	
	We return to the assumed starting point of the present paper, which is that the free dynamics operator $\BA_q$ in the open-loop controlled linear system \eqref{1.19} is unstable. Its properties are collected in Theorem \ref{I-Thm-1.7}. Accordingly, its eigenvalues satisfy the statement which includes their location in \eqref{1.21}. Denote by $P_N$ and $P_N^*$ (which actually depend on $q$) the projections given explicitly by \cite[p 178]{TK:1966}, \cite{BT:2004}, \cite{BLT1:2006}
	\begin{subequations}\label{1.22}
		\begin{align}
		\label{1.22a} P_N &= -\frac{1}{2 \pi i}\int_{\boldsymbol{\Gamma}}\left( \lambda I - \BA_q \right)^{-1}d \lambda : \mathbf{W}^q_{\sigma} \text{ onto } (\mathbf{W}^q_{\sigma})^u_N \subset \mathbf{L}^q_{\sigma}(\Omega) \times L^q(\Omega) \\
		\label{1.22b} P_N^* &= -\frac{1}{2 \pi i}\int_{\bar{\boldsymbol{\Gamma}}}\left( \lambda I - \BA_q^* \right)^{-1}d \lambda : (\mathbf{W}^q_{\sigma})^* \text{ onto } [(\mathbf{W}^q_{\sigma})^u_N]^* \subset \mathbf{L}^{q'}_{\sigma}(\Omega) \times L^{q'}(\Omega),
		\end{align}
	\end{subequations}
	
	\noindent where, here ${\boldsymbol{\Gamma}}$ (respectively, its conjugate counterpart $\bar{{\boldsymbol{\Gamma}}}$) is a smooth closed curve that separates the unstable spectrum from the stable spectrum of $\BA_q$ (respectively, $\BA_q^*$). As in \cite[Sect 3.4, p 37]{BLT1:2006}, following \cite{RT:1975}, \cite{RT:1980}, we decompose the space $\mathbf{W}^q_{\sigma}=\mathbf{W}^q_{\sigma}(\Omega) \equiv \lso \times \lqo$ into the sum of two complementary subspaces (not necessarily orthogonal):	
	\begin{multline}\label{1.23}
	\mathbf{W}^q_{\sigma} = (\mathbf{W}^q_{\sigma})^u_N \oplus (\mathbf{W}^q_{\sigma})^s_N; \ (\mathbf{W}^q_{\sigma})^u_N \equiv P_N \mathbf{W}^q_{\sigma} ;\ (\mathbf{W}^q_{\sigma})^s_N \equiv (I - P_N) \mathbf{W}^q_{\sigma}; \\ \text{ dim } (\mathbf{W}^q_{\sigma})^u_N = N
	\end{multline}
	
	\noindent where each of the spaces $(\mathbf{W}^q_{\sigma})^u_N$ and $(\mathbf{W}^q_{\sigma})^s_N$ (which depend on $q$, but we suppress such dependence) is invariant under $\BA_q$, and let	
	\begin{equation}\label{1.24}
	\BA^u_{q,N} = P_N \BA_q = \BA_q |_{(\mathbf{W}^q_{\sigma})^u_N} ; \quad \BA^s_{q,N} = (I - P_N) \BA_q = \BA_q |_{(\mathbf{W}^q_{\sigma})^s_N}
	\end{equation}
	
	\noindent be the restrictions of $\BA_q$ to $(\mathbf{W}^q_{\sigma})^u_N$ and $(\mathbf{W}^q_{\sigma})^s_N$ respectively. The original point spectrum (eigenvalues) $\{ \lambda_j \}_{j=1}^{\infty} $ of $\BA_q$ is then split into two sets	
	\begin{equation}\label{1.25}
	\sigma (\BA^u_{q,N}) = \{ \lambda_j \}_{j=1}^{N}; \quad  \sigma (\BA^s_{q,N}) = \{ \lambda_j \}_{j=N+1}^{\infty},
	\end{equation}
	
	\noindent and $(\mathbf{W}^q_{\sigma})^u_N$ is the generalized eigenspace of $\BA^u_{q,N}$ in (\ref{1.24}). The system (\ref{1.19}) on $\mathbf{W}^q_{\sigma} \equiv \lso \times \lqo$ can accordingly be decomposed as	
	\begin{equation}\label{1.26}
	\mathbf{w} = \mathbf{w}_N + \boldsymbol{\zeta}_N, \quad \mathbf{w}_N = P_N \mathbf{w}, \quad \boldsymbol{\zeta}_N = (I-P_N)\mathbf{w}.
	\end{equation}
	
	\noindent After applying $P_N$ and $(I-P_N)$ (which commute with $\BA_q$) on (\ref{1.19}), we obtain via (\ref{1.24})	
	\begin{subequations}\label{1.27} 
		\begin{equation}\label{1.27a}
		\text{on } (\mathbf{W}^q_{\sigma})^u_N: \mathbf{w}'_N - \BA^u_{q,N} \mathbf{w}_N = P_N \bbm P_{q} (m\mathbf{u}) \\ mv \ebm; \quad  \mathbf{w}_N(0) = P_N \bbm \mathbf{w}_f (0) \\ w_h(0) \ebm
		\end{equation}
		\begin{equation}\label{1.27b}
		\text{on } (\mathbf{W}^q_{\sigma})^s_N: \boldsymbol{\zeta}'_N - \BA^s_{q,N} \boldsymbol{\zeta}_N = (I-P_N) \bbm P_{q} (m\mathbf{u}) \\  mv \ebm; \quad  \boldsymbol{\zeta}_N(0) = (I - P_N) \bbm \mathbf{w}_f (0) \\ w_h(0) \ebm
		\end{equation}
	\end{subequations}
	\noindent respectively.
	
	\section{Proof of Theorem \ref{NThm-2.1}. Global well-posedness and uniform exponential stabilization of the linearized $\mathbf{w}$-problem (\ref{1.19}) on the space $\mathbf{W}^q_{\sigma}(\Omega) \equiv \lso \times \lqo$ or the space $\mathbf{V}^{q,p}(\Omega) \equiv  \Bto \times \lqo$}\label{sec-4}

	\noindent {\bf Orientation:} We  shall appeal to several  technical developments in \cite{LPT.1}, where the case of N-S equations has been studied. We will have to adopt and transfer some of its procedures to the  case of the Boussinesq system. Thus, properties such as maximal regularity and the entire  development for uniform stabilization of the linearized Boussinesq dynamics need to be established. It is worth noticing that while analyticity and maximal regularity are equivalent properties in the Hilbert setting \cite{DeS} this is not so in the Banach setting where maximal regularity is a more general and more delicate property. Moreover, it is only after we establish  uniform stabilization that we can claim maximal regularity up to infinity. This needs to be asserted by direct analysis of singular integrals. To proceed, we recall the state space $\ds \mathbf{W}^q_{\sigma}(\Omega) = \lso \times \lqo \equiv (\mathbf{W}^q_{\sigma})^u_N \oplus (\mathbf{W}^q_{\sigma})^s_N$. The unstable uncontrolled operator is $\ds \BA_q$ in (\ref{1.18}).\\
	
	\noindent The same strategy employed in \cite{LPT.1} leading to the linearized $\mathbf{w}$-dynamics of the translated non-linear $\mathbf{z}$-problem in terms of finite dimensional feedback controls can be followed  now, where $\mathbf{w}$ and $\mathbf{z}$ are augmented vectors consisting of a fluid vector component and a scalar thermal component. Thus we first seek to establish that the finite dimensional projection - the $\mathbf{w}_N$-equation presently in (\ref{1.27a}) is controllable on $(\mathbf{W}^q_{\sigma})^u_N$, hence exponentially stabilizable with an arbitrarily large decay rate \cite[p. 44]{Za}. Next, one then examines the corresponding $\boldsymbol{\zeta}_N$-equation, presently (\ref{1.27b}), where the arbitrarily large decay rate of the feedback $\mathbf{w}_N$-equation combined with the exponential stability on $(\mathbf{W}^q_{\sigma})^s_N$ of the s.c. analytic semigroup $\ds e^{\BA_{q,N}^s t}$ yields the desired result. All this works thanks to the present corresponding Unique Continuation Property, this time of the Boussinesq system, that is Theorem \ref{T1.4}. This implies the Kalman algebraic condition \eqref{bbN4-5}. A proof is given in Appendix \ref{app-B}. Such condition \eqref{bbN4-5} is equivalent to  linear independence of certain vectors occurring in (\ref{bbN4-5}) below, so that the finite dimensional projected $\mathbf{w}_N$-dynamics satisfies the controllability condition of Kalman or Hautus. See corresponding cases in \cite{RT:2009}, \cite{RT:2008}. More precisely, the counterpart of the analysis of \cite[Section 3 and 4]{LPT.1} has now some novel aspects which account for the present additional property that the UCP of Theorem \ref{T1.4} involves only the first $(d-1)$ components of the fluid vector. This conceptual advantage will be responsible for heavier notation as described below.\\
	
	\noindent For each $i = 1, \dots, M$, we denote by  $\ds \{\boldsymbol{\Phi}_{ij}\}_{j=1}^{\ell_i},  \{\boldsymbol{\Phi}^*_{ij}\}_{j=1}^{\ell_i}$ the normalized, linearly independent eigenfunctions of $\BA_q$, respectively $\BA_q^*$, say, on
	\begin{multline}
	\mathbf{W}^q_{\sigma}(\Omega) \equiv \lso \times \lqo
	\mbox{ and } \\
	(\mathbf{W}^q_{\sigma}(\Omega))^* \equiv (\lso)' \times (\lqo)' = \mathbf{L}^{q'}_{\sigma}(\Omega) \times L^{q'}(\Omega), \quad
	\frac{1}{q} + \frac{1}{q'} = 1,
	\end{multline}
	\noindent (where in the last equality we have invoked  Remark \ref{rmkB.1} in Appendix \ref{app-B}) corresponding to the $M$ distinct unstable eigenvalues $\lambda_1, \ldots, \lambda_M$ of $\BA_q$ and $\overline{\lambda}_1, \ldots, \overline{\lambda}_M$ of $\BA_q^*$ respectively, either on $\mathbf{W}^q_{\sigma}$ or on $\mathbf{V}^{q,p}$:
	\begin{align}
	\BA_q \boldsymbol{\Phi}_{ij} &= \lambda_i \boldsymbol{\Phi}_{ij} \in \calD(\BA_q) = [\mathbf{W}^{2,q}(\Omega) \cap \mathbf{W}^{1,q}_0(\Omega) \cap \lso] \times [W^{2,q}(\Omega) \cap W^{1,q}_0(\Omega)]  \label{2.3}\\
	\BA_q^* \boldsymbol{\Phi}_{ij}^* &= \bar{\lambda}_i \boldsymbol{\Phi}_{ij}^* \in \calD(\BA_q^*) = [\mathbf{W}^{2,q'}(\Omega) \cap \mathbf{W}^{1,q'}_0(\Omega) \cap \lo{q'}] \times [W^{2,q'}(\Omega) \cap W^{1,q'}_0(\Omega)]. \label{2.4}
	\end{align}
	We now express the eigenvectors $ \boldsymbol{\Phi^*_{ij}}$ in terms of their coordinates, as $(d+1)$ vectors:
	
	\begin{equation} \label{B4-19}
	\boldsymbol{\Phi^*_{ij}} = \{ \boldsymbol{\varphi^*_{ij}} , \psi^*_{ij} \} = \{ \varphi_{ij}^{*(1)} ,\varphi_{ij}^{*(2)},...,\varphi_{ij}^{*(d-1)}, \varphi_{ij}^{*(d)}, \psi^*_{ij} \} \,, \text{ a } (d+1)\text{-vector}.
	\end{equation}
	With reference to (\ref{B4-19}), we introduce the following corresponding d-vector
	\begin{equation} 
	\boldsymbol{\widehat{\Phi }^*_{ij} = \{ \boldsymbol{\widehat{\varphi}^*_{ij}}} , \psi^*_{ij} \} = \{ \varphi_{ij}^{*(1)} ,\varphi_{ij}^{*(2)},...\varphi_{ij}^{*(d-1)},  \psi^*_{ij} \} \,, \text{ a } d\text{-vector},
	\end{equation}
	obtained from $  \boldsymbol{\Phi^*_{ij}}$ by omitting the d-component $\varphi_{ij}^{*(d)}$ of the vector $\boldsymbol{\Phi^*_{ij}} $. Next,we construct the following matrix $U_i$ of size $\ell_i \times K, \ K =\sup \{ \ell_i: i = 1, \dots, M \}$
	\begin{equation}\label{Bb2.5}
	U_i =
	\begin{bmatrix}
	(\mathbf{u}_1,\widehat{\Phi}_{i1}^*)_{\omega} & \dots & (\mathbf{u}_K,\widehat{\Phi}_{i1}^*)_{\omega} \\[1mm]
	(\mathbf{u}_1,\widehat{\Phi}_{i2}^*)_{\omega} & \dots & (\mathbf{u}_K,\widehat{\Phi}_{i2}^*)_{\omega} \\
	\vdots & \ddots & \vdots \\
	(\mathbf{u}_1,\widehat{\Phi}_{i \ell_i}^*)_{\omega} & \dots & (\mathbf{u}_K,\widehat{\Phi}_{i \ell_i}^*)_{\omega} \\
	\end{bmatrix}
	:\ell_i \times K.
	\end{equation}
	Here we have set
	\begin{subequations}
		\begin{equation}
		\mathbf{u}_k = [ \mathbf{u}_k^1 , u_k^2 ] = [ (u_k^1)^{(1)}, (u_k^1)^{(2)}...(u_k^1)^{(d-1)},u_k^2] \in (\widehat{\mathbf{W}}^q_{\sigma})^u_N \subset \widehat{\mathbf{L}}^{q}_{\sigma}(\Omega) \times L^q(\Omega)
		\end{equation}
		\begin{multline}
		\widehat{ \mathbf{L}}^q_{\sigma} (\Omega) \equiv \text{ the space obtained from } \mathbf{L}^q_{\sigma} (\Omega) \text{ after omitting the } \\ d\text{-coordinate from the vectors of } \widehat{ \mathbf{L}}^q_{\sigma} (\Omega)
		\end{multline}
		\begin{multline}\label{E4.7c}
		(\widehat{\mathbf{W}}^q_{\sigma})^u_N \equiv \text{ the space obtained from } (\mathbf{W}^q_{\sigma})^u_N \text{ after omitting the }\\ d\text{-coordinate from the vectors of } (\mathbf{W}^q_{\sigma})^u_N.
		\end{multline}
		
	\end{subequations}
	In \eqref{Bb2.5}, we have defined the duality pairing over $ \omega $ as
	
	\begin{align}
	(\mathbf{u}_k, \widehat{\Phi}^*_{ij})_{\omega} &= \bpm \bbm \mathbf{u}^1_k \\ u^2_k \ebm, \bbm \widehat{\boldsymbol{\varphi}}^*_{ij} \\ \psi^*_{ij} \ebm \epm_{\omega} = \int_{\omega} [ \mathbf{u}^1_k \cdot \widehat{\boldsymbol{\varphi}}^*_{ij} + u^2_k \psi^*_{ij} ] d \omega \nonumber \\
	&= (\mathbf{u}^1_k , \widehat{\boldsymbol{\varphi}}^*_{i1})_{\widehat{\mathbf{L}}^q(\omega), \widehat{\mathbf{L}}^{q'}(\omega)} + (u^2_k, \psi^*_{i1})_{L^q(\omega), L^{q'}(\omega)} \label{4-8}
	\\
	&= \int_{\omega }
	\begin{bmatrix}
	( u_k^1)^{(1)}\\
	( u_k^1)^{(2)}\\
	\vdots \\
	( u_k^1)^{(d-1)}\\
	u_k^2
	\end{bmatrix} \cdot
	\begin{bmatrix}
	\varphi_{ij}^{*(1)}\\
	\varphi_{ij}^{*(2)}\\
	\vdots \\
	\varphi_{ij}^{*(d-1)}\\
	\psi^*_{ij}
	\end{bmatrix} \, d\omega.
	\end{align}
	The controllability Kalman/Hautus algebraic condition of the finite-dimensional projected $\mathbf{w}_N$-equation in (\ref{1.27a}) is given by
	\begin{equation}    \label{bbN4-5}
	\text{rank } U_i = \text{full} = \ell_i, \quad i = 1, \dots, M.
	\end{equation}
	\noindent It is proved in Appendix \ref{app-B} that the UCP of Theorem \ref{T1.4}$\implies $ rank condition \eqref{bbN4-5}. As a consequence, we obtain the following most critical result counterpart of \cite[Theorem 6.1]{LPT.1}.
	
	\begin{thm}\label{Thm-2.1}
		Let the operator $\BA_q$ have $N$ possibly repeated unstable eigenvalues $\{ \lambda_j \}_{j = 1}^N$ as in \eqref{1.21} of which $M$ are distinct. Let $\ell_i$ denote the geometric multiplicity of $\lambda_i$. Set $K = \sup \{ \ell_i; i = 1, \cdots, M \}$. Then, one may construct a feedback operator $F$
		
		\begin{equation}\label{2.8}
		F(\mathbf{w}) = \bbm P_q \Big( m \big( \sum_{k = 1}^{K} (\mathbf{w}_N, \mathbf{p}_k)_{\omega} \mathbf{u}^1_k \big) \Big) \\[2mm] m \big( \sum_{k = 1}^{K} (\mathbf{w}_N, \mathbf{p}_k)_{\omega} u^2_k \big)\ebm
		\end{equation}
		\noindent with vectors $[\mathbf{u}^1_k, u^2_k] \in (\widehat{\mathbf{W}}^q_{\sigma})^u_N \subset  \widehat{\mathbf{L}}^{q}_{\sigma}(\Omega) \times L^q(\Omega)$, (so that $ \mathbf{u}^1_k $ is $ (d-1)$ dimensional, see \eqref{E4.7c}) and $\mathbf{p}_k \in ((\mathbf{W}^q_{\sigma})^u_N)^* \subset \lo{q'} \times L^{q'}(\Omega)$ such that the $\mathbf{w}$-problem (\ref{1.19}) can be rewritten in feedback form on $\mathbf{W}^q_{\sigma}(\Omega)$ or $\mathbf{V}^{q,p}(\Omega)$ as follows ($\mathbf{w}_N = P_N\mathbf{w}$):
		
		\begin{multline}\label{N4-7}
		\frac{d \mathbf{w}}{dt} = \frac{d}{dt} \bbm \mathbf{w}_f \\[2mm] w_h \ebm = \BA_q \mathbf{w} + F(\mathbf{w}) = \BA_q \mathbf{w} + \bbm P_q \Big( m \big( \sum_{k = 1}^{K} (\mathbf{w}_N, \mathbf{p}_k)_{\omega} \mathbf{u}^1_k \big) \Big) \\[2mm] m \big( \sum_{k = 1}^{K} (\mathbf{w}_N, \mathbf{p}_k)_{\omega} u^2_k \big)\ebm = \BA_{F,q} \mathbf{w};\\ \mathbf{w}(0) = \mathbf{w}_0
		\end{multline}
		
		with $\ds \calD(\BA_{F,q}) = \calD(\BA_q)$ where the operator $\ds \BA_{F,q}$ in (\ref{2.9}) or (\ref{N4-7}) has the following properties:
		
		\begin{enumerate}[(i)]
			\item It is the generator of a s.c. analytic semigroup $\ds e^{\BA_{F,q}t}$ in the space $\mathbf{W}^q_{\sigma}(\Omega) = \lso \times \lqo$ as well as in the space $\mathbf{V}^{q,p}(\Omega) = \Bto \times \lqo$.
			\item It is uniformly (exponentially) stable in  either of these spaces
			\begin{equation}\label{N4-9}
			\norm{e^{\BA_{F,q}t} \mathbf{w}_0}_{(\cdot)} \leq C_{\gamma_0^{}} e^{-\gamma_0 t}\norm{\mathbf{w}_0}_{(\cdot)}
			\begin{picture}(0,0)
			\put(-230,-10){ $\left\{\rule{0pt}{22pt}\right.$}\end{picture}
			\end{equation}
			\hspace{.3cm} where $(\cdot)$ denotes either $\ds \lso \times \lqo \equiv \SqsO$ or else $\ds \Bto \times \lqo \equiv \VqpO$.\\
			In \eqref{N4-9}, $\gamma_0$ is any positive number such that $Re~\lambda_{N+1} < -\gamma_0 < 0$.
			\item Finally, $\BA_{F,q}$ has maximal $L^p$-regularity up to $T = \infty$ on either of these spaces:
			\begin{equation}\label{N4-10}
			\BA_{F,q} \in MReg (L^p(0,\infty; \ \cdot \ )) \text{ where } (\cdot) \text{ denotes}
			\begin{picture}(0,0)
			\put(-240,-10){ $\left\{\rule{0pt}{18pt}\right.$}\end{picture}
			\end{equation}
			\hspace{1cm} either $\ds \lso \times \lqo \equiv \SqsO$ or else $\ds \Bto \times \lqo \equiv \VqpO$ \qedsymbol
		\end{enumerate}
	\end{thm}
	The PDE version of the closed-loop abstract model \eqref{2.9} or \eqref{N4-7}, is (refer to the open-loop \eqref{N1.23})
	\begin{subequations}\label{N4.10}
		\begin{align}
		\frac{d}{dt}\mathbf{w}_f - \nu \Delta \mathbf{w}_f + L_e(\mathbf{w}_f) - \gamma w_h \mathbf{e}_d + \nabla \chi &=  m \bigg( \sum_{k = 1}^{K} (\mathbf{w}_N, \mathbf{p}_k)_{\omega} \mathbf{u}^1_k \bigg)   \text{ in } Q \label{N4.10a}\\
		\frac{d}{dt}w_h - \kappa \Delta w_h + \mathbf{y}_e \cdot \nabla w_h + \mathbf{w}_f \cdot \nabla w_h + \mathbf{w}_f \cdot \nabla \theta_e &= m \bigg( \sum_{k = 1}^{K} (\mathbf{w}_N, \mathbf{p}_k)_{\omega} u^2_k \bigg) \text{ in } Q \label{N4.10b}\\
		\begin{picture}(0,0)
		\put(-205,30){ $\left\{\rule{0pt}{70pt}\right.$}\end{picture}
		\text{div } \mathbf{w}_f &= 0   \text{ in } Q \label{N4.10c}\\
		\mathbf{w}_f \equiv 0, \ w_h &\equiv 0 \text{ on } \Sigma \label{N4.10d}\\
		\mathbf{w}_f(0,\cdot) = \mathbf{w}_{f,0}; \quad w_h(0,\cdot) & = w_{h,0} \text{ on } \Omega. \label{N4.10e}
		\end{align}
	\end{subequations}	
	\noindent What property (iii) in Theorem \ref{Thm-2.1} means explicitly is singled out in the next result. This result will be critically used in the subsequent non linear analysis of Section \ref{new-sec-5} and \ref{new-sec-6}.
	\begin{thm}\label{Thm-2.2}
		With respect to the operator $\ds \BA_{F,q}$ in (\ref{2.9}) or \eqref{N4-7} of Theorem \ref{Thm-2.1}, recall
		\begin{equation}\label{N4-11}
		\SqsO \equiv \lso \times \lqo, \quad \VqpO \equiv \Bto \times \lqo, \ 1 < p < \frac{2q}{2q - 1}.
		\end{equation}
		\noindent Then, the following properties hold true.		
		\begin{enumerate}[(i)]
			\item \begin{align}\label{2.13}
			\mathbf{F} &\longrightarrow \int_{0}^{t} e^{\BA_{F,q}(t - \tau)}\mathbf{F}(\tau) d \tau: \text{ continuous } \nonumber\\
			L^p(0,\infty; \SqsO) &\longrightarrow L^p(0,\infty; \calD(\BA_{F,q})=\calD(\calA_q) \times \calD(\calB_q))
			\end{align}
			\noindent whereby then automatically from \eqref{N4-7}
			\begin{equation}\label{2.14}
			\hspace{3cm}	L^p(0,\infty; \SqsO) \longrightarrow W^{1,p}(0, \infty; \SqsO)
			\end{equation}
			\noindent and ultimately, on the space of maximal regularity for $\mathbf{w}=\{\mathbf{w}_f,w_h\}$ w.r.t. the operator $ \BA_{F,q}$, see \eqref{N3.12}
			\begin{equation}\label{2.15}
			\hspace{1cm} L^p(0,\infty; \SqsO) \longrightarrow \xipqs \times \xipq \equiv L^p(0,\infty; \calD(\BA_{F,q})) \cap W^{1,p}(0, \infty; \SqsO).
			\end{equation}		
			\item The s.c. analytic uniformly stable semigroup $\ds e^{\BA_{F,q}t}$	on the space $\VqpO \equiv \Bto \times \lqo$ satisfies
			\begin{multline}\label{2.16}
			e^{\BA_{F,q}t}: \text{ continuous } \Bto \times \lqo \equiv \VqpO \longrightarrow \\ \xipqs \times \xipq  \equiv L^p(0, \infty; \calD (\BA_{F,q})) \cap W^{1,p}(0,\infty;\SqsO)
			\end{multline}	
			\noindent or equivalently, recalling $ \calD(\BA_{F,q}) = \calD(\BA_q) = \calD( \calA_q) \times \calD(\calB_q)$ in \eqref{1.18}, or \eqref{E2.9},\eqref{E2.10}:
			\begin{align}\label{2.17}
			&\longrightarrow \xipqs \times \xipq  \\[2mm] \nonumber
			&\equiv  L^p(0,\infty,[\mathbf{W}^{2,q}(\Omega) \cap \mathbf{W}^{1,q}_0(\Omega) \cap \lso] \times \\[2mm]
			& \qquad [W^{2,q}(\Omega) \cap W^{1,q}_0(\Omega)])
			\cap W^{1,p}(0,\infty;\SqsO) \hookrightarrow C([0,\infty] ; \mathbf{V}^{q,p}(\Omega)) \label{E4.22}.
			\end{align}
			see \eqref{E2.111}-\eqref{E2.13}.
			\noindent Equivalently, in summary, with reference to the $\mathbf{w}$-equation (\ref{2.9}), $ \mathbf{w}=\{ \mathbf{w}_f, w_h\}$ for the following estimates hold true for $\ds 1 < p < \frac{2q}{2q - 1}$:
			\begin{align}
			C_0 \norm{\mathbf{w}}_{C(0,\infty; \VqpO)} &\leq \norm{\mathbf{w}}_{\xipqs \times \xipq } + \norm{\pi}_{\yipq} \nonumber\\
			&\leq \norm{\mathbf{w}_t}_{L^p(0,\infty; \SqsO)} + \norm{\BA_{F,q} w}_{L^p(0,\infty; \SqsO)} + \norm{\pi}_{\yipq} \nonumber\\
			& \leq C_1 \big\{ \norm{\mathbf{F}}_{L^p(0,\infty; \SqsO)} + \norm{\mathbf{w}_0}_{\VqpO} \big\}. \label{2.18}
			\end{align}
		\end{enumerate}
		where $ \yipq $ is defined in \eqref{A-1.29} of Appendix \ref{app-A}.
	\end{thm}
	\noindent For the first inequality in \eqref{2.18} we have recalled the embedding, called trace theorem \cite[Theorem 4.10.2, p. 180, BUC for $T=\infty$]{HA:2000}, \cite{PS:2016}, as in \cite[Eq. (1.30)]{LPT.1}. See also \eqref{3.23} below.
	
	\newsavebox{\zhvector}
	\begin{lrbox}{\zhvector} $\bbm \mathbf{z} \\ h \ebm$ \end{lrbox}
	
	\section{Proof of Theorem \ref{N-Thm-2.2}. Well-posedness on $\xipqs \times \xipq$ of the non-linear \usebox{\zhvector}-dynamics in feedback form}\label{new-sec-5}
	
	\noindent In this section we return to the translated non-linear $\ds \bbm \mathbf{z} \\ h \ebm$-dynamics (\ref{1.16}) or (\ref{1.17}) and apply to it the feedback control pair $\{\mathbf{u},v\}$
	
	\begin{equation}\label{3.1}
	\bbm P_q m(\mathbf{u}) \\[5mm] m(v) \ebm = \bbm P_q \Bigg( m \bigg( \sum_{k = 1}^{K} \bigg(P_N \bbm \mathbf{z} \\ h \ebm, \mathbf{p}_k \bigg)_{\omega} \mathbf{u}^1_k \bigg) \Bigg) \\[2mm] m \bigg( \sum_{k = 1}^{K} \bigg(P_N \bbm \mathbf{z} \\ h \ebm, \mathbf{p}_k \bigg)_{\omega} u^2_k \bigg) \ebm
	=
	F\left(\bbm \mathbf{z} \\ h \ebm\right)
	\end{equation}
	
	\noindent that is, of the same structure as the feedback $F$ in (\ref{2.8}) identified on the RHS of the linearized $\ds \mathbf{w} = \bbm \mathbf{w}_f \\ w_h \ebm$-dynamics (\ref{N4-7}). This is the feedback operator which produced the s.c. analytic, uniformly stable semigroup $\ds e^{\BA_{F,q}t}$ on $\ds \SqsO = \lso \times \lqo$ or on $ \mathbf{V}^{q,p}(\Omega) \equiv  \Bt \times \lqo$, (Theorem \ref{Thm-2.1}) possessing $L^p$-maximal regularity on these spaces up to $T = \infty$; see (\ref{2.11}) or \eqref{N4-10}. Thus, returning to (\ref{1.16}) or (\ref{1.17}), in this section we consider the following feedback nonlinear problem, see \eqref{N2.5}
	
	\begin{equation}\label{3.2}
	\frac{d}{dt} \bbm \mathbf{z} \\[2mm] h \ebm = \BA_{F,q} \bbm \mathbf{z} \\[2mm] h \ebm - \bbm \calN_q & 0 \\[2mm] 0 & \calM_q[\mathbf{z}] \ebm \bbm \mathbf{z} \\[2mm] h \ebm; \quad \BA_{F,q} = \bbm \calA_q & -\calC_{\gamma} \\[2mm] -\calC_{\theta_e} & \calB_q \ebm +F = \BA_q + F
	\end{equation}
	
	\noindent specifically
	
	\begin{subequations}\label{3.3}
		\begin{align}
		\frac{d\mathbf{z}}{dt} - \calA_q \mathbf{z} + \calC_{\gamma} h + \calN_q \mathbf{z} &= P_q \Bigg( m \bigg( \sum_{k = 1}^{K} \bigg(P_N \bbm \mathbf{z} \\ h \ebm, \mathbf{p}_k \bigg)_{\omega} \mathbf{u}^1_k \bigg) \Bigg) \label{3.3a}\\
		\begin{picture}(0,0)
		\put(-20,20){ $\left\{\rule{0pt}{40pt}\right.$}\end{picture}
		\frac{dh}{dt} - \calB_q h + \calC_{\theta_e} \mathbf{z} + \calM_q[\mathbf{z}]h &= m \bigg( \sum_{k = 1}^{K} \bigg(P_N \bbm \mathbf{z} \\ h \ebm, \mathbf{p}_k \bigg)_{\omega} u^2_k \bigg). \label{3.3b}
		\end{align}
	\end{subequations}
	
	\noindent The variation of parameter formula for Eq (\ref{3.2}) is
	
	\begin{equation}\label{3.4}
	\bbm \mathbf{z} \\ h \ebm (t) = e^{\BA_{F,q}t} \bbm \mathbf{z}_0 \\ h_0 \ebm - \int_{0}^{t} e^{\BA_{F,q}(t - \tau)} \bbm \calN_q \mathbf{z}(\tau) \\ \calM_q[\mathbf{z}]h(\tau) \ebm d \tau.
	\end{equation}
	
	\noindent We already know from (\ref{2.10}) or
	\eqref{N4-9} that for $\ds \{\mathbf{z}_0, h_0\} \in \Bto \times \lqo \equiv \VqpO, \ 1 < p < \frac{2q}{2q-1}$, we have: there is $M_{\gamma_0}$ such that
	
	\begin{equation}\label{3.5}
	\norm{e^{\BA_{F,q}t} \bbm \mathbf{z}_0 \\ h_0 \ebm}_{\VqpO} \leq M_{\gamma_0} e^{-\gamma_0t} \norm{\bbm \mathbf{z}_0 \\ h_0 \ebm}_{\VqpO}, \ t \geq 0
	\end{equation}
	\noindent with $M_{\gamma_0}$ possibly depending on $p,q$. Maximal $L^p$-regularity properties corresponding to the solution operator formula (\ref{3.4}) were established in Theorem \ref{Thm-2.2}. Accordingly, for
	
	\begin{align}
	\mathbf{b}_0 &\equiv \{\mathbf{z}_0, h_0\} \in \Bto \times \lqo \equiv \VqpO \label{3.6}\\
	\mathbf{f} &\equiv \{\mathbf{f}_1, f_2\} \in \xipqs \times \xipq \equiv L^p(0, \infty, \calD (\BA_{F,q})) \cap W^{1,p}(0,\infty; \SqsO) \label{3.7}\\
	\calD(\BA_{F,q}) &= \calD(\BA_q) = \calD(\calA_q) \times \calD(\calB_q) = \nonumber\\
	& \hspace{2cm} [\mathbf{W}^{2,q}(\Omega) \cap \mathbf{W}^{1,q}_0(\Omega) \cap \lso] \times [W^{2,q}(\Omega) \cap W^{1,q}_0(\Omega)]\\
	\SqsO &= \lso \times \lqo \nonumber\\
	\xipqs &\equiv L^p(0,\infty; \mathbf{W}^{2,q}(\Omega) \cap \mathbf{W}^{1,q}_0(\Omega) \cap \lso) \cap W^{1,p}(0,\infty; \lso);  \tag{5.9a} \label{5.9a} \\
	\xipq &\equiv L^p(0,\infty; W^{2,q}(\Omega) \cap W^{1,q}_0(\Omega)) \cap W^{1,q}(0,\infty; \lqo)	    \tag{5.9b}, \label{5.9b}
	\end{align}
	\addtocounter{equation}{1}
	\noindent as in \eqref{N3.12}-\eqref{E2.10} or \eqref{2.17}, we define the operator
	\begin{equation}\label{3.10}
	\calF (\mathbf{b}_0, \mathbf{f}) \equiv e^{\BA_{F,q}t} \mathbf{b}_0 - \int_{0}^{t} e^{\BA_{F,q}(t - \tau)} \bbm \calN_q \mathbf{f}_1(\tau) \\ \calM_q[\mathbf{f}_1]f_2(\tau) \ebm d \tau.
	\end{equation}
	
	\noindent The main result of this section is Theorem \ref{N-Thm-2.2} restated as
	
	\begin{thm}\label{Thm-3.1}
		Let $\ds d = 2,3, \ q > d, \ 1 < p < \frac{2q}{2q-1}$. There exists a positive constant $r_1 > 0$ (identified in the proof below in (\ref{3.31})) such that if
		\begin{equation}\label{3.11}
		\norm{\mathbf{b}_0}_{\VqpO} = \norm{\{ \mathbf{z}_0, h_0 \}}_{\Bto \times \lqo} < r_1,
		\end{equation}
		then the operator $\calF$ in (\ref{3.10}) has a unique fixed point non-linear semigroup solution in $\xipqs \times \xipq$, see \eqref{N3.12}-\eqref{E2.10} or  (\ref{2.15})--(\ref{2.17}) or \eqref{5.9a}--\eqref{5.9b}
		\begin{equation}\label{3.12}
		\calF \bpm \bbm \mathbf{z}_0 \\ h_0 \ebm, \bbm \mathbf{z} \\ h \ebm \epm = \bbm \mathbf{z} \\ h \ebm, \ \text{or} \ \bbm \mathbf{z} \\ h \ebm (t) =  e^{\BA_{F,q}t} \bbm \mathbf{z}_0 \\ h_0 \ebm - \int_{0}^{t} e^{\BA_{F,q}(t - \tau)} \bbm \calN_q \mathbf{z}(\tau) \\ \calM_q[\mathbf{z}]h(\tau) \ebm d \tau,
		\end{equation}
		\noindent which therefore is the unique solution of problem (\ref{3.2}) = (\ref{3.3}) in $\ds \xipqs \times \xipq$ defined in \eqref{N3.12}-\eqref{E2.10} or (\ref{2.15})--(\ref{2.17}) or \eqref{5.9a}--\eqref{5.9b}.
	\end{thm}
	
	\noindent The proof of Theorem \ref{N-Thm-2.2} = Theorem \ref{Thm-3.1} is accomplished in two steps.\\	
	
	\noindent \underline{Step 1:}\\
	
	\begin{thm}\label{Thm-3.2}
		Let $d = 2,3,  \ q > d$ and $\ds 1 < p < \frac{2q}{2q - 1}$. There exists a positive constant $r_1 > 0$ (identified in the proof below (\ref{3.31})) and a subsequent constant $r > 0$ (identified in the proof below in (\ref{3.29})) depending on $r_1 > 0$ and a constant $C$ in (\ref{3.28}), such that with $\ds \norm{\mathbf{b}_0}_{\VqpO} < r_1$ as in (\ref{3.11}), the operator $\ds \calF(\mathbf{b}_0,\mathbf{f})$ in (\ref{3.10}) maps a ball $B(0,r)$ in $\ds \xipqs \times \xipq$ into itself. \qedsymbol
	\end{thm}
	
	\noindent Theorem \ref{Thm-3.1} will follow then from Theorem \ref{Thm-3.2} after establishing that\\
	
	\noindent \underline{Step 2:}\\
	
	\begin{thm}\label{Thm-3.3}
		Let $d = 2,3,  \ q > d$ and $\ds 1 < p < \frac{2q}{2q - 1}$. There exists a positive constant $r_1 > 0$ such that if $\ds \norm{\mathbf{b}_0}_{\VqpO} < r_1$ as in (\ref{3.11}), there exists a constant $0 < \rho_0 < 1$ (identified in (\ref{3.56})), such that the operator $\calF(\mathbf{b}_0, \mathbf{f})$ in (\ref{3.10}) defines a contraction in the ball $B(0,\rho_0)$ of $\xipqs \times \xipq$. \qedsymbol
	\end{thm}	
	
	\noindent The Banach contraction principle then establishes Theorem \ref{Thm-3.1}, once we prove Theorems \ref{Thm-3.2} and  \ref{Thm-3.3}. These are proved below. \\
	
	\noindent \underline{Proof of Theorem \ref{Thm-3.2}}\\
	
	\noindent \textit{Step 1:} We start from the definition (\ref{3.10}) of $\ds \calF(\mathbf{b}_0, \mathbf{f})$ and invoke the maximal regularity properties \eqref{2.16} for $\ds e^{\BA_{F,q}t}$ and (\ref{2.15}) for the integral term in  (\ref{3.10}). We then obtain from (\ref{3.10})
	
	\begin{align}
	\norm{\calF (\mathbf{b}_0, \mathbf{f})}_{\xipqs \times \xipq} &\leq \norm{e^{\BA_{F,q}t}\mathbf{b}_0}_{\xipqs \times \xipq} + \nonumber\\
	& \hspace{1cm}\norm{\int_{0}^{t} e^{\BA_{F,q}(t - \tau)} \bbm \calN_q \mathbf{f}_1 (\tau) \\[2mm] \calM_q[\mathbf{f}_1] f_2 (\tau) \ebm d \tau}_{\xipqs \times \xipq} \label{3.13}\\
	\leq C \bigg[ \norm{\mathbf{b}_0}_{\VqpO} + &\norm{\calN_q \mathbf{f}_1}_{\lplqs} + \norm{\calM_q[\mathbf{f}_1]f_2}_{\lplq} \bigg]. \label{3.14}
	\end{align}
	
	\noindent \textit{Step 2:} Regarding the term $\ds \calN_q \mathbf{f}_1$ we can invoke \cite[Eq (8.19)]{LPT.1} to obtain
	
	\begin{equation}
	\norm{\calN_q \mathbf{f}_1}_{\lplqs} \leq C \norm{\mathbf{f}_1}^2_{\xipqs}, \ \mathbf{f}_1 \in \xipqs \label{3.15}.
	\end{equation}
	
	\noindent Regarding the term $\ds \calM_q[\mathbf{f}_1]f_2$, we can trace the proof in \cite[from (8.10) $\longrightarrow $ (8.18)]{LPT.1} (which yielded estimate (\ref{3.14})). For the sake of clarity, we shall reproduce the computations in the present case with $\ds \calM_q [\mathbf{f}_1] f_2 = \mathbf{f}_1 \cdot \nabla f_2$, see (\ref{1.13}), \textit{mutatis mutandis}. We shall obtain
	
	\begin{equation}\label{3.16}
	\norm{\calM_q [\mathbf{f}_1] f_2}_{\lplq} \leq C \norm{\mathbf{f}_1}_{\xipqs} \norm{f_2}_{\xipq}, \ \mathbf{f}_1 \in \xipqs, f_2 \in \xipq.
	\end{equation}
	
	\noindent In fact, let us compute
	
	\begin{align}
	\norm{\calM_q [\mathbf{f}_1] f_2}_{\lplq}^p &\leq \int_{0}^{\infty} \norm{\mathbf{f}_1 \cdot \nabla f_2}^p_{L^q(\Omega)} dt \label{3.17}\\
	&\leq \int_{0}^{\infty} \bigg\{ \int_{\Omega} \abs{\mathbf{f}_1(t,x)}^q \abs{\nabla f_2(t,x)}^q d\Omega \bigg\}^{\rfrac{p}{q}}dt\\
	&\leq \int_{0}^{\infty} \bigg\{  \bigg[ \sup_{\Omega} \abs{\nabla f_2(t,x)}^q  \bigg]^{\rfrac{1}{q}} \bigg[\int_{\Omega} \abs{\mathbf{f}_1(t,x)}^q d\Omega  \bigg]^{\rfrac{1}{q}} \bigg\}^p dt\\
	&\leq \int_{0}^{\infty} \norm{\nabla f_2 (t, \cdot)}^p_{\mathbf{L}^{\infty}(\Omega)} \norm{\mathbf{f}_1(t, \cdot)}^p_{\lso} dt\\
	&\leq \sup_{0 \leq t \leq \infty} \norm{\mathbf{f}_1(t, \cdot)}^p_{\mathbf{L}_{\sigma}^q(\Omega)} \int_{0}^{\infty} \norm{\nabla f_2 (t, \cdot)}^p_{\mathbf{L}^{\infty}(\Omega)} dt\\
	&= \norm{\mathbf{f}_1}^p_{L^{\infty}(0,\infty; \lso)}\norm{\nabla f_2}^p_{L^p(0,\infty;\mathbf{L}^{\infty}(\Omega))}. \label{3.22}
	\end{align}
	\noindent See also \cite[Eq (8.14)]{LPT.1}.\\
	
	\noindent \textit{Step 3}: The following embeddings hold true (see the stronger Eq \eqref{E2.111}):
	\begin{enumerate}[(i)]
		\item \cite[Proposition 4.3, p 1406 with $\mu = 0, s = \infty, r = q$]{GGH:2012}  so that the required formula reduces to $1 \geq \rfrac{1}{p}$, as desired
		\begin{subequations}\label{3.23}
			\begin{align}
			\mathbf{f}_1 \in \xipqs \hookrightarrow \mathbf{f}_1 &\in L^{\infty}(0,\infty; \lso) \label{3.23a}\\
			\text{ so that, } \norm{\mathbf{f}_1}_{L^{\infty}(0,\infty; \lso)} &\leq C\norm{\mathbf{f}_1}_{\xipqs}; \label{3.23b}
			\end{align}
		\end{subequations}
		\item \cite[Theorem 2.4.4, p 74 requiring $C^1$-boundary]{SK:1989}
		\begin{equation}\label{3.24}
		W^{1,q}(\Omega) \subset L^{\infty}(\Omega) \text{ for q}>\text{dim }\Omega = d, \ d = 2,3,
		\end{equation}
	\end{enumerate}
	
	\noindent so that, with $p>1, q>d$:
	\begin{align}
	\norm{\nabla f_2}^p_{L^p(0,\infty; \mathbf{L}^{\infty}(\Omega))} &\leq C \norm{ \nabla f_2}^p_{L^p(0,\infty; \mathbf{W}^{1,q}(\Omega))} \leq C \norm{f_2}^p_{L^p(0,\infty; W^{2,q}(\Omega))} \label{3.25}\\
	&\leq C \norm{f_2}^p_{\xipq}. \label{3.26}
	\end{align}
	
	\noindent In going from (\ref{3.25}) to (\ref{3.26}) we have recalled the definition of $f_2 \in \xipq$ in \eqref{5.9b} or \eqref{E2.10}, as $f_2$ was taken at the outset in $\ds \calD(\calB_q) \subset L^q(\Omega)$. Then the sought-after final estimate (\ref{3.16}) of the nonlinear term $\ds \calM_q[\mathbf{f}_1]f_2$ is obtained from substituting (\ref{3.23b}) and (\ref{3.26}) into the RHS of (\ref{3.22}).\\
	
	\noindent \textit{Step 4:} Substituting estimates (\ref{3.15}) and (\ref{3.16}) on the RHS of (\ref{3.14}), we finally obtain
	
	\begin{equation}\label{3.27}
	\norm{\calF(\mathbf{b}_0, \mathbf{f}}_{\xipqs \times \xipq} \leq C \Big\{ \norm{\mathbf{b}_0}_{\VqpO} + \norm{\mathbf{f}_1}_{\xipqs} \big( \norm{\mathbf{f}_1}_{\xipqs} + \norm{f_2}_{\xipq} \big) \Big\}.
	\end{equation}
	
	\noindent See \cite[Eqt (8.20)]{LPT.1}.\\
	
	\noindent \textit{Step 5:} We now impose restrictions on the data on the RHS of (\ref{3.27}): $\mathbf{b}_0$ is in a ball of radius $r_1 > 0$ in $\ds \VqpO = \Bto \times \lqo$ and $\mathbf{f} = \{ \mathbf{f}_1, f_2\}$ lies in a ball of radius $r > 0$ in $\ds \xipqs \times \xipq$. We further demand that the final result $\ds \calF(\mathbf{b}_0, \mathbf{f})$ shall lie in a ball of radius $r > 0$ in $\ds \xipqs \times \xipq$. Thus, we obtain from \eqref{3.27}
	
	\begin{align}
	\norm{\calF(\mathbf{b}_0, \mathbf{f})}_{\xipqs \times \xipq} &\leq C \Big\{ \norm{\mathbf{b}_0}_{\VqpO} + \norm{\mathbf{f}_1}_{\xipqs} \big( \norm{\mathbf{f}_1}_{\xipqs} + \norm{f_2}_{\xipq} \big) \Big\} \nonumber\\
	&\leq C(r_1 + r \cdot r) \leq r. \label{3.28}
	\end{align}
	
	\noindent This implies
	\begin{equation}\label{3.29}
	Cr^2 - r + Cr_1 \leq 0 \quad \text{or} \quad \frac{1 - \sqrt{1-4C^2r_1}}{2C} \leq r \leq \frac{1 + \sqrt{1-4C^2r_1}}{2C}
	\end{equation}
	\noindent whereby
	\begin{equation}\label{3.30}
	\begin{Bmatrix}
	\text{ range of values of r }
	\end{Bmatrix}
	\longrightarrow \text{ interval } \Big[ 0, \frac{1}{C} \Big], \text{ as } r_1 \searrow 0,
	\end{equation}
	\noindent a constraint which is guaranteed by taking
	\begin{equation}\label{3.31}
	r_1 \leq \frac{1}{4C^2},\ C \text{ being the constant in } (\ref{3.28}) ( \text{w.l.o.g } C > 1/4 ).
	\end{equation}
	\noindent We have thus established that by taking $r_1$ as in (\ref{3.31})  and subsequently $r$ as in (\ref{3.29}), then the map
	\begin{multline}\label{3.32}
	\calF(\mathbf{b}_0, \mathbf{f}) \text{ takes: }
	\begin{Bmatrix}
	\text{ ball in } \VqpO \\
	\text{of radius } r_1
	\end{Bmatrix}
	\times
	\begin{Bmatrix}
	\text{ ball in } \xipqs \times \xipq \\
	\text{of radius } r
	\end{Bmatrix}
	\text{ into }\\
	\begin{Bmatrix}
	\text{ ball in } \xipqs \times \xipq \\
	\text{of radius } r
	\end{Bmatrix}, \ d < q, \ 1 < p < \frac{2q}{2q-1}.
	\end{multline}
	\noindent This establishes Theorem \ref{Thm-3.2}. \qedsymbol \\
	
	\noindent \underline{Proof of Theorem \ref{Thm-3.3}}\\
	
	\noindent \textit{Step 1:} For $\mathbf{f} = \{ \mathbf{f}_1, f_2 \}, \mathbf{g} = \{ \mathbf{g}_1, g_2 \}$ both in the ball of $\ds \xipqs \times \xipq$ of radius $r$ obtained in the proof of Theorem \ref{Thm-3.2}, we estimate from (\ref{3.10}):
	
	\begin{align}
	&\norm{\calF (\mathbf{b}_0, \mathbf{f}) - \calF(\mathbf{b}_0,\mathbf{g})}_{\xipqs \times \xipq} = \nonumber\\ & \hspace{2.5cm} \norm{\int_{0}^t e^{\BA_{F,q}(t - \tau)} \bbm \calN_q\mathbf{f}_1(\tau) - \calN_q \mathbf{g}_1(\tau) \\[2mm] \calM_q[\mathbf{f}_1]f_2(\tau) - \calM_q[\mathbf{g}_1]g_2(\tau)\ebm d \tau}_{\xipqs \times \xipq} \nonumber\\
	&\leq \wti{m} \Big[ \norm{\calN_q\mathbf{f}_1 - \calN_q\mathbf{g}_1}_{\lplqs} + \norm{\calM_q[\mathbf{f}_1]f_2 - \calM_q[\mathbf{g}_1]g_2}_{L^p(0,\infty;L^q(\Omega))} \Big] \label{3.33}
	\end{align}
	\noindent after invoking the maximal regularity property (\ref{2.15})--(\ref{2.17}).\\
	
	\noindent \textit{Step 2:} As to the first term or the RHS of (\ref{3.33}), we can invoke \cite[Eq (8.41)]{LPT.1} and obtain
	\begin{equation}\label{3.34}
	\norm{\calN_q\mathbf{f}_1 - \calN_q\mathbf{g}_1}_{\lplqs} \leq 2^{\rfrac{1}{p}}C^{\rfrac{1}{p}} \norm{\mathbf{f}_1 - \mathbf{g}_1}_{\xipqs} \big( \norm{\mathbf{f}_1}_{\xipqs} + \norm{\mathbf{g}_1}_{\xipqs} \big).
	\end{equation}
	\noindent Regarding the second term on the RHS of (\ref{3.33}) involving $\calM_q$, we can track the proof of \cite[from (8.28) to (8.41)]{LPT.1} (which yielded estimate (\ref{3.34})) \textit{mutatis mutandis}. For the sake of clarity, we shall reproduce the computations in the present case, recalling from (\ref{1.13}) that $\ds \calM_q[\mathbf{f}_1]f_2 = \mathbf{f}_1 \cdot \nabla f_2, \calM_q[\mathbf{g}_1]g_2 = \mathbf{g}_1 \cdot \nabla g_2$. We shall obtain
	\begin{multline}\label{3.35}
	\norm{\calM_q[\mathbf{f}_1]f_2 - \calM_q[\mathbf{g}_1]g_2}^p_{\lplq} \leq C \Big\{ \norm{\mathbf{f}_1 - \mathbf{g}_1}^p_{\xipqs} \norm{f_2}^p_{\xipq} \\+ \norm{\mathbf{g}_1}^p_{\xipqs} \norm{f_2 - g_2}^p_{\xipq} \Big\}.
	\end{multline}
	\noindent In fact, adding and subtracting
	\begin{align}
	\calM_q[\mathbf{f}_1]f_2 - \calM_q[\mathbf{g}_1]g_2 &= \mathbf{f}_1 \cdot \nabla f_2 - \mathbf{g}_1 \cdot \nabla g_2 \nonumber\\
	&= \mathbf{f}_1 \cdot \nabla f_2 - \mathbf{g}_1 \cdot \nabla f_2 + \mathbf{g}_1 \cdot \nabla f_2 - \mathbf{g}_1 \cdot \nabla g_2 \nonumber\\
	&= (\mathbf{f}_1 - \mathbf{g}_1)\cdot \nabla f_2 + \mathbf{g}_1 \cdot \nabla (f_2 - g_2) = A + B. \label{3.36}
	\end{align}
	\noindent Thus, using $(*): \ \abs{A + B}^p \leq 2^p \big[ \abs{A}^p + \abs{B}^p  \big]$ \cite[p. 12]{TL:1980}, we estimate
	\begin{align}
	\norm{\calM_q[\mathbf{f}_1]f_2 - \calM_q[\mathbf{g}_1]g_2}_{\lplq}^p  & = \int_{0}^{\infty} \bigg\{ \bigg[ \int_{\Omega} \abs{\mathbf{f}_1 \cdot \nabla f_2 - \mathbf{g}_1 \cdot \nabla g_2}^q d \Omega \bigg]^{\rfrac{1}{q}}\bigg\}^p dt\\
	\text{(by (\ref{3.36}))} \qquad &= \int_{0}^{\infty} \bigg[ \int_{\Omega} \abs{A+B}^q d \Omega \bigg]^{\rfrac{p}{q}}dt\\
	&\leq 2^p \int_{0}^{\infty} \bigg\{ \int_{\Omega} \big[\abs{A}^q + \abs{B}^q \big] d \Omega \bigg\}^{\rfrac{p}{q}}dt\\
	&= 2^p \int_{0}^{\infty} \bigg\{ \Big[ \int_{\Omega} \abs{A}^q d \Omega + \int_{\Omega} \abs{B}^q d \Omega   \Big]^{\rfrac{1}{q}} \bigg\}^p dt\\
	&= 2^p \int_{0}^{\infty} \bigg\{ \Big[ \norm{A}^q_{L^q(\Omega)} + \norm{B}^q_{L^q(\Omega)} \Big]^{\rfrac{1}{q}} \bigg\}^p dt\\
	(\text{by } (*) \mbox{ with } p \to \frac{1}{q}) \qquad
	&\leq 2^p \cdot 2^{\rfrac{p}{q}}\int_{0}^{\infty} \Big\{ \norm{A}_{L^q(\Omega)} + \norm{B}_{L^q(\Omega)} \Big\}^p dt\\
	(\text{by } (*) ) \qquad
	&\leq 2^{p+p+\rfrac{p}{q}}\int_{0}^{\infty} \Big[ \norm{A}^p_{L^q(\Omega)} + \norm{B}^p_{L^q(\Omega)} \Big] dt \label{3.43}\\
	&= 2^{p+p+\rfrac{p}{q}}\int_{0}^{\infty} \Big[ \norm{(\mathbf{f}_1 - \mathbf{g}_1)\cdot \nabla f_2}^p_{L^q(\Omega)}\nonumber \\ &\hspace{2cm}+\norm{\mathbf{g}_1 \cdot \nabla(f_2 - g_2)}^p_{L^q(\Omega)} \Big] dt \label{3.44}\\
	&\leq  2^{p+p+\rfrac{p}{q}}\int_{0}^{\infty} \Big\{ \norm{\mathbf{f}_1 - \mathbf{g}_1}^p_{\lso} \norm{\nabla f_2}^p_{L^q(\Omega)}\nonumber \\ &\hspace{2cm}+\norm{\mathbf{g}_1}^p_{\lso} \norm{\nabla(f_2 - g_2)}^p_{L^q(\Omega)} \Big\} dt, \label{3.45}
	\end{align}
	\noindent recalling the definitions of $A$ and $B$ in (\ref{3.36}) in passing from (\ref{3.43}) to (\ref{3.44}). \cite[This is the counterpart Eq (8.36)]{LPT.1}. Next we proceed by majorizing the $L^q(\Omega)$-norm for $\nabla f_2$ and $\nabla (f_2 - g_2)$ by the $L^{\infty}(\Omega)$-norm. We obtain
	\begin{align}
	\norm{\calM_q[\mathbf{f}_1]f_2 - \calM_q[\mathbf{g}_1]g_2}_{\lplq}^p &\leq \text{RHS of (\ref{3.45})} \nonumber\\
	&\hspace{-3.5cm} \leq  c 2^{p+p+\rfrac{p}{q}} \bigg\{ \sup_{0 \leq t \leq \infty} \norm{\mathbf{f}_1 - \mathbf{g}_1}^p_{\lso} \int_{0}^{\infty} \norm{\nabla f_2}^p_{L^{\infty}(\Omega)}\nonumber \\ & \hspace{-1.5cm} + \sup_{0 \leq t \leq \infty} \norm{\mathbf{g}_1}^p_{\lso} \int_{0}^{\infty} \norm{\nabla(f_2 - g_2)}^p_{L^{\infty}(\Omega)} dt \bigg\} \label{3.46}\\
	&\hspace{-3.5cm} \leq c 2^{p+p+\rfrac{p}{q}} \Big\{\norm{\mathbf{f}_1 - \mathbf{g}_1}^p_{L^{\infty}(0,\infty; \lso)} \norm{\nabla f_2}^p_{L^p(0,\infty; L^{\infty}(\Omega))}\nonumber \\
	&\hspace{-1.85cm} + \norm{\mathbf{g}_1}^p_{L^{\infty}(0,\infty; \lso)} \norm{\nabla(f_2 - g_2)}^p_{L^p(0,\infty; L^{\infty}(\Omega))} \Big\}\label{3.47}\\
	& \hspace{-5.5cm} \text{ by (\ref{3.23b}) and (\ref{3.26})}  \nonumber \\
	& \hspace{-3.5cm}\leq C \Big\{ \norm{\mathbf{f}_1 - \mathbf{g}_1}_{\xipqs}^p \norm{f_2}_{\xipq}^p + \norm{\mathbf{g}_1}_{\xipqs}^p  \norm{f_2 - g_2}_{\xipq}^p \Big\},  \label{3.48}
	\end{align}
	\noindent counterpart of \cite[Eq (8.39)]{LPT.1}.\\
	
	\noindent \textit{Step 3:} We substitute estimate (\ref{3.34}) and estimate (\ref{3.48}) on the RHS of (\ref{3.33}) and obtain via $(*)$
	
	\begin{align}
	\norm{\calF(\mathbf{b}_0,\mathbf{f}) - \calF(\mathbf{b}_0,\mathbf{g})}^p_{\xipqs \times \xipq} &\leq 2^p \wti{m}^p  \Big\{ \norm{\calN_q \mathbf{f}_1 - \calN_q \mathbf{g}_1}^p_{\lplq} \nonumber
	\\&+ \norm{\calM_q [\mathbf{f}_1] f_2 - \calM_q[\mathbf{g}_1] g_2}^p_{L^p(0,\infty;L^q(\Omega))} \Big\} \label{3.49}\\
	&\hspace{-2.5cm} \leq C_p \wti{m}^p  \Big\{ \norm{\mathbf{f}_1 - \mathbf{g}_1}^p_{\xipqs} \Big( \norm{\mathbf{f}_1}_{\xipqs} + \norm{\mathbf{g}_1}_{\xipqs} \Big)^p \nonumber
	\\&\hspace{-2.2cm} + \norm{\mathbf{f}_1 - \mathbf{g}_1}^p_{\xipqs} \norm{f_2}^p_{\xipq} + \norm{\mathbf{g}_1}^p_{\xipqs} \norm{f_2 - g_2}^p_{\xipq} \Big\}. \label{3.50}
	\end{align}
	\noindent This is the counterpart of \cite[Eq (8.42)]{LPT.1}.\\
	
	\noindent \textit{Step 4:} Next pick the points $\mathbf{f} = \{ \mathbf{f}_1, f_2 \}$ and $\mathbf{g} = \{ \mathbf{g}_1, g_2 \}$ in a ball of $\ds \xipqs \times \xipq$ of radius $R$:
	\begin{align}
	\norm{\mathbf{f}}_{\xipqs \times \xipq} &= \norm{\mathbf{f}_1}_{\xipqs} + \norm{f_2}_{\xipq} < R \label{3.51}\\
	\norm{\mathbf{g}}_{\xipqs \times \xipq} &= \norm{\mathbf{g}_1}_{\xipqs} + \norm{g_2}_{\xipq} < R. \label{3.52}
	\end{align}
	\noindent Then (\ref{3.51}),\eqref{3.52} used (\ref{3.50}) implies
	\begin{align}
	\norm{\calF(\mathbf{b}_0,\mathbf{f}) - \calF(\mathbf{b}_0,\mathbf{g})}^p_{\xipqs \times \xipq} &\leq C_p \Big\{ \norm{\mathbf{f}_1 - \mathbf{g}_1}^p_{\xipqs} \big[ (2R)^p + R^p \big] \nonumber \\
	&\hspace{2cm}+ \norm{f_2 - g_2}^p_{\xipq} R^p \Big\}\\
	&\leq C_p \Big\{ \norm{\mathbf{f}_1 - \mathbf{g}_1}^p_{\xipqs} \big[ (2^p + 1)R^p \big] \nonumber \\
	&\hspace{2cm}+ \norm{f_2 - g_2}^p_{\xipq} R^p \Big\} \nonumber\\
	\big( K_p = C_p (2^p + 1) \big) \quad &\leq K_p R^p \Big\{ \norm{\mathbf{f}_1 - \mathbf{g}_1}^p_{\xipqs} + \norm{f_2 - g_2}^p_{\xipq} \Big\}\nonumber\\
	[a>0, b>0, \ a^p + b^p \leq (a + b)^p] \quad &\leq K_p R^p \Big\{ \norm{\mathbf{f}_1 - \mathbf{g}_1}_{\xipqs} + \norm{f_2 - g_2}_{\xipq} \Big\}^p \label{3.54}
	\end{align}
	using $K_p = C_p(2^p + 1)$ and $a^p + b^p \leq (a+b)^p$ for $a>0, b>0$.
	Finally,
	\begin{align}
	\norm{\calF(\mathbf{b}_0,\mathbf{f}) - \calF(\mathbf{b}_0,\mathbf{g})}_{\xipqs \times \xipq} &\leq  K_p^{\rfrac{1}{p}} R \Big\{ \norm{\mathbf{f}_1 - \mathbf{g}_1}_{\xipqs} + \norm{f_2 - g_2}_{\xipq} \Big\} \nonumber\\
	&= \rho_0 \norm{\mathbf{f} - \mathbf{g}}_{\xipqs \times \xipq} \label{3.55}
	\end{align}
	\noindent and $\ds \calF(\mathbf{b}_0,\mathbf{f})$ is a contraction on the space $\xipqs \times \xipq$ as soon as
	\begin{equation}\label{3.56}
	\rho_0 = K_p^{\rfrac{1}{p}}R < 1 \quad \text{or} \quad R < \frac{1}{K_p^{\rfrac{1}{p}}}.
	\end{equation}
	\noindent In this case, the map $\ds \calF(\mathbf{b}_0, \mathbf{f})$ defined in (\ref{3.10}) has a fixed point $\ds \bbm \mathbf{z} \\ h \ebm$ in $\xipqs \times \xipq$:
	\begin{equation}\label{3.57}
	\calF \bpm \mathbf{b}_0, \bbm \mathbf{z} \\ h \ebm \epm = \bbm \mathbf{z} \\ h \ebm, \ \text{or} \ \bbm \mathbf{z} \\ h \ebm (t) =  e^{\BA_{F,q}t} \bbm \mathbf{z}_0 \\ h_0 \ebm - \int_{0}^{t} e^{\BA_{F,q}(t - \tau)} \bbm \calN_q \mathbf{z}(\tau) \\[2mm] \calM_q[\mathbf{z}]h(\tau) \ebm d \tau,
	\end{equation}
	\noindent and such point $\ds \bbm \mathbf{z} \\ h \ebm$ is the unique solution of the translated non-linear system (\ref{3.2}), or (\ref{3.3}), with finite dimensional control
	\begin{equation*}
	\bbm \mathbf{u} \\[2mm] v \ebm = \bbm P_q \Bigg( m \bigg( \sum_{k = 1}^{K} \bigg(P_N \bbm \mathbf{z} \\ h \ebm, \mathbf{p}_k \bigg)_{\omega} \mathbf{u}^1_k \bigg) \Bigg) \\[2mm] m \bigg( \sum_{k = 1}^{K} \bigg(P_N \bbm \mathbf{z} \\ h \ebm, \mathbf{p}_k \bigg)_{\omega} u^2_k \bigg) \ebm
	\end{equation*}
	
	\noindent in feedback form, as described by Eq (\ref{3.2}). Theorem \ref{Thm-3.3} and hence Theorem \ref{Thm-3.1} are proved. \qedsymbol
	
	\begin{rmk} \label{R5.1}
		As $[\mathbf{u}^1_k, u^2_k] \in (\widehat{\mathbf{W}}^q_{\sigma})^u_N \subset  \widehat{\mathbf{L}}^{q}_{\sigma}(\Omega) \times L^q(\Omega)$, we have that the feedback control acting on the fluid variable $\mathbf{u}$ is of reduced dimension  $(d-1)$, see \eqref{E4.7c}.
	\end{rmk}
	
	\section{Proof of Theorem \ref{NThm-2.3}: local exponential decay of the non-linear $ [\mathbf{z},h]$ translated dynamics (\ref{3.2}) = (\ref{3.3}) with finite dimensional localized feedback controls and $(d-1)$ dimensional $\mathbf{u}^1_k $}
	\label{new-sec-6}
	
	\noindent In this section we return to the feedback problem (\ref{3.2}) = (\ref{3.3}), rewritten equivalently as in (\ref{3.4})
	\begin{equation}\label{4.1}
	\bbm \mathbf{z} \\ h \ebm (t) =  e^{\BA_{F,q}t} \bbm \mathbf{z}_0 \\ h_0 \ebm - \int_{0}^{t} e^{\BA_{F,q}(t - \tau)} \bbm \calN_q \mathbf{z}(\tau) \\[2mm] \calM_q[\mathbf{z}]h(\tau) \ebm d \tau,
	\end{equation}
	\noindent For $\mathbf{b}_0 = [\mathbf{z}_0, h_0]$ in a small ball of $\ds \VqpO = \Bto \times \lqo$, Theorem \ref{Thm-3.1}(= Theorem \ref{N-Thm-2.2}) provides a unique solution $\{\mathbf{z}, h\}$ in a small ball of $\ds \xipqs \times \xipq$. We recall from (\ref{N4-9})
	
	\begin{equation}\label{4.2}
	\norm{e^{\BA_{F,q}t} \bbm \mathbf{z}_0 \\ h_0 \ebm}_{\VqpO} \leq M_{\gamma_0} e^{-\gamma_0 t} \norm{\bbm \mathbf{z}_0 \\ h_0 \ebm}_{\VqpO}, t \geq 0
	\end{equation}
	\noindent $\ds \VqpO = \Bto \times \lqo$. Our goal is to show that for $\ds [\mathbf{z}_0, h_0]$ in a small ball of $\ds \VqpO$, problem (\ref{4.1}) satisfies the exponential decay
	\begin{equation}\label{4.3}
	\norm{\bbm \mathbf{z} \\ h \ebm(t)}_{\VqpO} \leq C_a e^{-a t} \norm{\bbm \mathbf{z}_0 \\ h_0 \ebm}_{\VqpO}, t \geq 0, \text{for some constants } a > 0, C = C_a \geq 1.
	\end{equation}
	
	\noindent \textit{Step 1:} Starting from (\ref{4.1}) and using the decay (\ref{4.2}), we estimate
	\begin{align}
	\norm{\bbm \mathbf{z}(t) \\ h(t) \ebm}_{\VqpO} &\leq  M_{\gamma_0}e^{-\gamma_0 t} \norm{\bbm \mathbf{z}_0 \\ h_0 \ebm}_{\VqpO} + \nonumber \\
	&\hspace{2cm} \sup_{0 \leq t \leq \infty} \norm{\int_{0}^{t} e^{\BA_{F,q}(t - \tau)} \bbm \calN_q \mathbf{z}(\tau) \\[2mm] \calM_q[\mathbf{z}]h(\tau) \ebm d \tau}_{\VqpO} \label{4.4}\\
	&\leq M_{\gamma_0}e^{-\gamma_0 t} \norm{\bbm \mathbf{z}_0 \\ h_0 \ebm}_{\VqpO} + \nonumber \\
	&\hspace{2cm} C_1 \norm{\int_{0}^{t} e^{\BA_{F,q}(t - \tau)} \bbm \calN_q \mathbf{z}(\tau) \\[2mm] \calM_q[\mathbf{z}]h(\tau) \ebm d \tau}_{\xipqs \times \xipq} \label{4.5}\\
	&\leq M_{\gamma_0}e^{-\gamma_0 t} \norm{\bbm \mathbf{z}_0 \\ h_0 \ebm}_{\VqpO} + C_2 \bigg[ \norm{\calN_q \mathbf{z}}_{\lplqs} + \nonumber \\
	&\hspace{5cm} \norm{\calM_q[\mathbf{z}]h}_{\lplq} \bigg] \label{4.6}\\
	&\leq M_{\gamma_0}e^{-\gamma_0 t} \norm{\bbm \mathbf{z}_0 \\ h_0 \ebm}_{\VqpO} + C_3 \bigg[ \norm{\mathbf{z}}^2_{\xipqs} + \norm{\mathbf{z}}_{\xipqs} \norm{h}_{\xipq} \bigg]. \label{4.7}
	\end{align}
	\noindent In going from (\ref{4.4}) to (\ref{4.5}) we have recalled the embedding $\ds \xipqs \times \xipq \hookrightarrow L^{\infty} \big(0, \infty; \Bto \times \lqo \big)$ from \eqref{E2.13} or \eqref{E4.22}. Next, in going from (\ref{4.5}) to (\ref{4.6}) we have used the maximal regularity property \eqref{2.13}--\eqref{2.16}. Finally, to go from (\ref{4.6}) to (\ref{4.7}), we have invoked (\ref{3.15}) for $\ds \calN_q \mathbf{z}$ and (\ref{3.16}) for $\ds \calM_q [\mathbf{z}]h$. Thus, the conclusion of \textit{Step 1} is
	\begin{align}
	\norm{\bbm \mathbf{z}(t) \\ h(t) \ebm}_{\VqpO} &\leq M_{\gamma_0}e^{-\gamma_0 t} \norm{\bbm \mathbf{z}_0 \\ h_0 \ebm}_{\VqpO} + C_3 \norm{\mathbf{z}}_{\xipqs} \bigg[ \norm{\mathbf{z}}_{\xipqs} + \norm{h}_{\xipq} \bigg] \label{4.8}\\
	&= M_{\gamma_0}e^{-\gamma_0 t} \norm{\bbm \mathbf{z}_0 \\ h_0 \ebm}_{\VqpO} + C_3 \norm{\mathbf{z}}_{\xipqs} \norm{[\mathbf{z}, h]}_{\xipqs \times \xipq}. \label{4.9}
	\end{align}
	
	\noindent \textit{Step 2:} We now return to (\ref{4.1}) and take the $\ds \xipqs \times \xipq$ norm across:
	\begin{align}
	\norm{\bbm \mathbf{z} \\ h \ebm(t)}_{\xipqs \times \xipq} &\leq  \norm{e^{\BA_{F,q}t} \bbm \mathbf{z}_0 \\ h_0 \ebm}_{\xipqs \times \xipq} + \nonumber \\
	&\hspace{1.5cm} \norm{\int_{0}^{t} e^{\BA_{F,q}(t - \tau)} \bbm \calN_q \mathbf{z}(\tau) \\[2mm] \calM_q[\mathbf{z}]h(\tau) \ebm d \tau}_{\xipqs \times \xipq} \label{4.10}\\
	(\text{by (\ref{2.16})}) \ &\leq M_1 \norm{\bbm \mathbf{z}_0 \\ h_0 \ebm}_{\VqpO} + C_3 \norm{\mathbf{z}}_{\xipqs} \norm{[\mathbf{z}, h]}_{\xipqs \times \xipq} \label{4.11}
	\end{align}
	\noindent by invoking (\ref{2.16}) on the first semigroup term in (\ref{4.10}) and the estimate from (\ref{4.5}) to (\ref{4.9}) on the second integral term in (\ref{4.10}). Thus (\ref{4.11}) is established and implies
	\begin{equation}\label{4.12}
	\norm{\bbm \mathbf{z} \\ h \ebm}_{\xipqs \times \xipq}  \Big[ 1 - C_3 \norm{\mathbf{z}}_{\xipqs} \Big]\leq M_1 \norm{\bbm \mathbf{z}_0 \\ h_0 \ebm}_{\VqpO}.
	\end{equation}
	\noindent [This is the counterpart of \cite[Eq (9.7)]{LPT.1}].\\
	
	\noindent \textit{Step 3:} The well-posedness Theorem \ref{Thm-3.1} says that
	\begin{equation}\label{4.13}
	\begin{Bmatrix}
	\text{ If } \norm{\bbm \mathbf{z}_0 \\ h_0 \ebm}_{\VqpO} \leq r_1 \\[4mm]
	\text{for } r_1 \text{ sufficiently small}
	\end{Bmatrix}
	\implies
	\begin{Bmatrix}
	\text{ The solution } [\mathbf{z},h] \text{ satisfies} \\[3mm]
	\norm{\bbm \mathbf{z} \\ h \ebm}_{\xipqs \times \xipq} \leq r
	\end{Bmatrix}
	\end{equation}
	
	\noindent where the constant $r$ satisfies the constraint (\ref{3.29}) in terms of the constant $r_1$ and some constant $C$ in (\ref{3.27}) or (\ref{3.28}) or (\ref{3.15}) or (\ref{3.16}). We seek to guarantee that we have
	\begin{equation}\label{4.14}
	\norm{\mathbf{z}}_{\xipqs} \leq \frac{1}{2C_3} < \frac{1}{2C},  \mbox{ hence }
	\frac{1}{2} \leq \Big[ 1 - C_3 \norm{\mathbf{z}}_{\xipqs} \Big],
	\end{equation}
	\noindent where we can always take the constant $C_3$ in \eqref{4.11} greater than the constant $C$ in \eqref{3.29}, \eqref{3.30}, \eqref{3.31}.Then \eqref{4.14} can be achieved by choosing $r_1 > 0$ sufficiently small. In fact, as $r_1 \searrow 0$. Eq. \eqref{3.30} shows that the interval $r_{\min} \leq r \leq r_{\max}$ of corresponding values of $r$ tends to the interval $[0,\frac{1}{C}]$. Thus, \eqref{4.14} can be achieved as $r_{\min} \searrow 0$: $0 < r_{\min} < r < \frac{1}{2C}$. Next (\ref{4.14}) implies that (\ref{4.12}) becomes
	\begin{equation}\label{4.15}
	\norm{\bbm \mathbf{z} \\ h \ebm}_{\xipqs \times \xipq} \leq 2M_1 \norm{\bbm \mathbf{z}_0 \\ h_0 \ebm}_{\VqpO} \leq 2M_1r_1.
	\end{equation}
	\noindent by \eqref{4.13}. Substituting (\ref{4.15}) in estimate (\ref{4.9}) then yields
	\begin{align}
	\norm{\bbm \mathbf{z}(t) \\ h(t) \ebm}_{\VqpO} &\leq M_{\gamma_0} e^{-\gamma_0 t}  \norm{\bbm \mathbf{z}_0 \\ h_0 \ebm}_{\VqpO} + 2M_1 C_3 \norm{\mathbf{z}}_{\xipqs} \norm{\bbm \mathbf{z}_0 \\ h_0 \ebm}_{\VqpO} \nonumber\\
	&\leq \what{M} \Big[ e^{-\gamma_0 t} + 4 \what{M}_1 C_3 r_1 \Big] \norm{\bbm \mathbf{z}_0 \\ h_0 \ebm}_{\VqpO} \label{4.16}
	\end{align}
	\noindent again by (\ref{4.15}) with $\ds \what{M} = \max \{ M_{\gamma_0}, M_1 \}$. This is the counterpart of \cite[Eq (9.16)]{LPT.1}.\\
	
	\noindent \textit{Step 4:} We now take $T$ sufficiently large and $r_1 > 0$ sufficiently small so that
	\begin{equation}\label{4.17}
	\beta = \what{M} e^{-\gamma_0 T} + 4 \what{M} C_3 r_1 < 1.
	\end{equation}
	\noindent Then (\ref{4.16}) implies
	\begin{equation}\label{4.18}
	\norm{\bbm \mathbf{z}(T) \\ h(T) \ebm}_{\VqpO} \leq \beta \norm{\bbm {\mathbf{z}}_0 \\ h_0 \ebm}_{\VqpO}
	\end{equation}
	\noindent $\ds \VqpO = \Bto \times \lqo$, and hence
	\begin{equation}\label{4.19}
	\norm{\bbm \mathbf{z}(nT) \\ h(nT) \ebm}_{\VqpO} \leq \beta \norm{\bbm \mathbf{z}((n-1)T) \\ h((n-1)T) \ebm}_{\VqpO} \leq \beta^n \norm{\bbm \mathbf{z}_0 \\ h_0 \ebm}_{\VqpO}.
	\end{equation}
	\noindent Since $\beta < 1$, the semigroup property of the evolution implies \cite{Bal:1981} that there are constants $M_{\wti{\gamma}} \geq 1, \wti{\gamma} >0$ such that
	\begin{equation}\label{4.20}
	\norm{\bbm \mathbf{z}(t) \\ h(t) \ebm}_{\VqpO} \leq M_{\wti{\gamma}} e^{-\wti{\gamma}t} \norm{\bbm \mathbf{z}_0 \\ h_0 \ebm}_{\VqpO}
	\end{equation}
	\noindent with $\ds \norm{[\mathbf{z}_0, h_0]}_{\VqpO}  \leq r_1 = $ small. This proves \eqref{N2-9}, i.e Theorem \ref{NThm-2.3} \qedsymbol
	
	\begin{rmk}\label{I-Rmk-9.1}
		The above computations - (\ref{4.17}) through (\ref{4.19}) - can be used to support qualitatively the intuitive expectation that ``the larger the decay rate $\gamma_0$ in (\ref{2.10}) =\eqref{N4-9}= (\ref{4.2}) of the linearized feedback $\mathbf{w}$-dynamics (\ref{2.9}), the larger the decay rate $\wti{\gamma}$ in (\ref{4.20}) of the nonlinear feedback $\{\mathbf{z},h\}$-dynamics (\ref{N2.5}) = (\ref{N2.6}) or \eqref{3.2},\eqref{3.3}; hence the larger the rate $\wti{\gamma}$ in (\ref{E2.19}) of the original $\{\mathbf{y},\theta\}$-dynamics in (\ref{I-2.27}) = (\ref{I-2.28})".\\
		
		\noindent The following considerations are somewhat qualitative. Let $S(t)$ denote the non-linear semigroup in the space $\mathbf{V}^{q,p}(\Omega) \equiv  \Bto \times \lqo$, with infinitesimal generator $\ds \big[ \BA_{_{F,q}} - \calN_q \big]$ describing the feedback $\{\mathbf{z},h\}$-dynamics (\ref{N2.5}) = (\ref{N2.6}), or \eqref{3.2},\eqref{3.3} as guaranteed by the well posedness Theorem \ref{N-Thm-2.2} = Theorem \ref{Thm-3.1}. Thus, $\bbm \mathbf{z}(t) \\ h(t) \ebm = S(t)\bbm \mathbf{z}_0 \\ h_0 \ebm$ on $\mathbf{V}^{q,p}(\Omega)$. By (\ref{4.17}), we can rewrite (\ref{4.18}) as:
		\begin{equation}\label{4.21}
		\norm{S(T)}_{\calL \big(\mathbf{V}^{q,p}(\Omega) \big)} \leq \beta < 1.
		\end{equation}
		\noindent It follows from \cite[p 178]{Bal:1981} via the semigroup property that
		\begin{equation}\label{4.22}
		- \wti{\gamma} \ \ \text{is just below} \ \ \frac{\ln \beta}{T} < 0.
		\end{equation}
		\noindent Pick $r_1 > 0$ in (\ref{4.17}) so small that $4\widehat{M}^2 C_3 r_1$ is negligible, so that $\beta$ is just above $\ds \widehat{M} e^{- \gamma_0 T}$, so $\ln \beta$ is just above $\ds \big[ \ln \widehat{M} - \gamma_0 T \big]$, hence
		\begin{equation}\label{4.23}
		\frac{\ln \beta}{T} \text{ is just above } \Big[ (-\gamma_0) + \frac{\ln \widehat{M}}{T}\Big].
		\end{equation}
		\noindent Hence, by (\ref{4.22}), (\ref{4.23}),
		\begin{equation}\label{4.24}
		\wti{\gamma} \sim \gamma_0 - \frac{\ln \widehat{M}}{T}
		\end{equation}
		\noindent and the larger $\gamma_0$, the larger is $\wti{\gamma}$, as desired.
	\end{rmk}

	\begin{appendices}
		
		\section{Some auxiliary results for the Stokes and Oseen operators: analytic semigroup generation, maximal regularity, domains of fractional powers}
		\label{app-A}
		
		\setcounter{equation}{0}
		
		\setcounter{theorem}{0}
		
		\renewcommand{\theequation}{{\rm A}.\arabic{equation}}
		\renewcommand{\thetheorem}{{\sc A}.\arabic{theorem}}
		
		In this section we collect some known results used in the paper. As a prerequisite of the present Appendix \ref{app-A}, we make reference to the paragraph  \textbf{Definition of Besov spaces $\mathbf{B}^s_{q,p}$ on domains of class $C^1$ as real interpolation of Sobolev spaces}, Eqts (\ref{1.8})-(\ref{A-1.16}) and Remark \ref{remark1-11}.

		\begin{enumerate}[(a)]
			\item \textbf{The Stokes and Oseen operators generate a strongly continuous analytic semigroup on $\lso$, $1 < q < \infty$}.\\
			
			\begin{thm}\label{A-Thm-1.2}
				Let $d \geq 2, 1 < q < \infty$ and let $\Omega$ be a bounded domain in $\mathbb{R}^d$ of class $C^3$. Then
				\begin{enumerate}[(i)]
					\item the Stokes operator $-A_q = P_q \Delta$ in (\ref{N1.9}), repeated here as
					\begin{equation}\label{A-1.17}
					-A_q \psi  = P_q \Delta \psi , \quad
					\psi \in \mathcal{D}(A_q) = \mathbf{W}^{2,q}(\Omega) \cap \mathbf{W}^{1,q}_0(\Omega) \cap \lso
					\end{equation}
					generates a s.c. analytic semigroup $e^{-A_qt}$ on $\lso$. See \cite{Gi:1981} and the review paper \cite[Theorem 2.8.5 p 17]{HS:2016}.			
					\item The Oseen operator $\calA_q$ in (\ref{1.9.new}) \label{A-Thm-1.2(iii)}
					\begin{equation}\label{A-1.18}
					\calA_q  = - (\nu A_q + A_{o,q}), \quad \calD(\calA_q) = \calD(A_q) \subset \lso
					\end{equation}
					generates a s.c. analytic semigroup $e^{\calA_qt}$ on $\lso$. This follows as $A_{o,q}$ is relatively bounded with respect to $A^{\rfrac{1}{2}}_q$,
					see (\ref{N1-11}): thus a standard theorem on perturbation of an analytic semigroup generator applies \cite[Corollary 2.4, p 81]{P:1983}.
					
					\item \begin{subequations}\label{A-1.19}
						\begin{align}
						0 \in \rho (A_q) &= \text{ the resolvent set of the Stokes operator } A_q\\
						\begin{picture}(0,0)
						\put(-40,10){ $\left\{\rule{0pt}{18pt}\right.$}\end{picture}
						A_q^{-1} &: \lso \longrightarrow \lso \text{ is compact}
						\end{align}
					\end{subequations}			
					\item The s.c. analytic Stokes semigroup $e^{-A_qt}$ is uniformly stable on $\lso$: there exist constants $M \geq 1, \delta > 0$ (possibly depending on $q$) such that
					\begin{equation}\label{A-1.20}
					\norm{e^{-A_qt}}_{\calL(\lso)} \leq M e^{-\delta t}, \ t > 0.
					\end{equation}
					
				\end{enumerate}
			\end{thm}
			\item \textbf{Domains of fractional powers, $\calD(A_q^{\alpha}), 0 < \alpha < 1$ of the Stokes operator $A_q$ on $\lso, 1 < q < \infty$}.
			We elaborate on (\ref{A-1.211}a-b)
			\begin{thm}\label{A-Thm-1.3}
				For the domains of fractional powers $\calD(A_q^{\alpha}), 0 < \alpha < 1$, of the Stokes operator $A_q$ in \eqref{A-1.17} = (\ref{N1.9}), the following complex interpolation relation holds true \cite{Gi:1985} and \cite[Theorem 2.8.5, p 18]{HS:2016}
				\begin{equation}\label{A-1.21}
				[ \calD(A_q), \lso ]_{1-\alpha} = \calD(A_q^{\alpha}), \ 0 < \alpha < 1, \  1 < q < \infty;
				\end{equation}
				in particular
				\begin{equation}\label{A-1.22}
				[ \calD(A_q), \lso ]_{\frac{1}{2}} = \calD(A_q^{\rfrac{1}{2}}) \equiv \mathbf{W}_0^{1,q}(\Omega) \cap \lso.
				\end{equation}
				Thus, on the space $\calD(A_q^{\rfrac{1}{2}})$, the norms
				\begin{equation}\label{A-1.23}
				\norm{\nabla \ \cdot \ }_{L^q(\Omega)} \text{ and } \norm{ \ }_{L^q(\Omega)}
				\end{equation}
				are related via Poincar\'{e} inequality.
			\end{thm}
			
			\item \textbf{The Stokes operator $-A_q$ and the Oseen operator $\calA_q, 1 < q < \infty$ generate s.c. analytic semigroups on the Besov space, from (\ref{A-1.16})} \label{A-Sec-1.10d}
			\begin{subequations}\label{A-1.24}
				\begin{align}
				\Big( \lso,\mathcal{D}(A_q) \Big)_{1-\frac{1}{p},p} &= \Big\{ \mathbf{g} \in \Bso : \text{ div } \mathbf{g} = 0, \ \mathbf{g}|_{\Gamma} = 0 \Big\} \quad \text{if } \frac{1}{q} < 2 - \frac{2}{p} < 2; \label{A-1.24a}\\
				\Big( \lso,\mathcal{D}(A_q) \Big)_{1-\frac{1}{p},p} &= \Big\{ \mathbf{g} \in \Bso : \text{ div } \mathbf{g} = 0, \ \mathbf{g}\cdot \nu|_{\Gamma} = 0 \Big\} \equiv \Bto \label{A-1.24b}\\
				&\hspace{6cm} \text{ if } 0 < 2 - \frac{2}{p} < \frac{1}{q}.\nonumber
				\end{align}	
			\end{subequations}
			Theorem \ref{A-Thm-1.2} states that the Stokes operator $-A_q$ generates a s.c. analytic semigroup on the space $\lso, \ 1 < q < \infty$, hence on the space $\calD(A_q)$ in (\ref{N1.9}) = (\ref{A-1.17}), with norm $\ds \norm{ \ \cdot \ }_{\calD(A_q)} = \norm{ A_q \ \cdot \ }_{\lso}$ as $0 \in \rho(A_q)$. Then, one obtains that the Stokes operator $-A_q$ generates a s.c. analytic semigroup on the real interpolation spaces in (\ref{A-1.24}). Next, the Oseen operator $\calA_q = -(\nu A_q + A_{o,q})$ in \eqref{A-1.18} = (\ref{1.9.new})  likewise generates a s.c. analytic semigroup $\ds e^{\calA_q t}$ on $\ds \lso$ since $A_{o,q}$ is relatively bounded w.r.t. $A_q^{\rfrac{1}{2}},$ as $A_{o,q}A_q^{-\rfrac{1}{2}}$ is bounded on $\lso$. Moreover $\calA_q$ generates a s.c. analytic semigroup on $\ds \calD(\calA_q) = \calD(A_q)$ (equivalent norms). Hence $\calA_q$ generates a s.c. analytic semigroup on the real interpolation space of (\ref{A-1.24}). Here below, however, we shall formally state the result only in the case $\ds 2-\rfrac{2}{p} < \rfrac{1}{q}$. i.e. $\ds  1 < p < \rfrac{2q}{2q-1}$, in the space $\ds \Bto$, as this does not contain B.C., Remark \ref{remark1-11}. The objective of the present paper is precisely to obtain stabilization results on spaces that do not recognize B.C.
			
			\begin{thm}\label{A-Thm-1.4}
				Let $1 < q < \infty, 1 < p < \rfrac{2q}{2q-1}$
				\begin{enumerate}[(i)]
					\item The Stokes operator $-A_q$ in \eqref{A-1.17} = (\ref{N1.9}) generates a s.c. analytic semigroup $e^{-A_qt}$ on the space $\Bt(\Omega)$ defined in (\ref{A-1.16bb}) = (\ref{A-1.24b}) which moreover is uniformly stable, as in (\ref{A-1.20}),
					\begin{equation}\label{A-1.25}
					\norm{e^{-A_qt}}_{\calL \big(\Bto\big)} \leq M e^{-\delta t}, \quad t > 0.
					\end{equation}
					\item The Oseen operator $\calA_q$ in \eqref{A-1.18}=(\ref{1.9.new}) generates a s.c. analytic semigroup $e^{\calA_qt}$ on the space $\Bto$ in (\ref{A-1.16bb}) = (\ref{A-1.24}).
				\end{enumerate}
			\end{thm}
			\item \textbf{Space of maximal $L^p$ regularity on $\lso$ of the Stokes operator $-A_q, \ 1 < p < \infty, \ 1 < q < \infty $ up to $T = \infty$.}
			We return to the dynamic Stokes problem in $\{\boldsymbol{\varphi}(t,x), \pi(t,x) \}$
			\begin{subequations}\label{A-1.26}
				\begin{align}
				\boldsymbol{\varphi}_t - \Delta \boldsymbol{\varphi} + \nabla \pi &= \mathbf{F} &\text{ in } (0, T] \times \Omega \equiv Q   \label{N!-13a}\\		
				div \ \boldsymbol{\varphi} &\equiv 0 &\text{ in } Q\\
				\begin{picture}(0,0)
				\put(-80,5){ $\left\{\rule{0pt}{35pt}\right.$}\end{picture}
				\left.\boldsymbol{\varphi} \right \rvert_{\Sigma} &\equiv 0 &\text{ in } (0, T] \times \Gamma \equiv \Sigma\\
				\left. \boldsymbol{\varphi} \right \rvert_{t = 0} &= \boldsymbol{\varphi}_0 &\text{ in } \Omega,
				\end{align}
			\end{subequations}
			
			rewritten in abstract form, after applying the Helmholtz projection $P_q$ to (\ref{N!-13a}) and recalling $A_q$ in \eqref{A-1.17} = (\ref{N1.9}) as
			\begin{equation}\label{A-1.27}
			\boldsymbol{\varphi}' + A_q \boldsymbol{\varphi} = \Fs \equiv P_q \mathbf{F}, \quad \boldsymbol{\varphi}_0 \in \lqaq.
			\end{equation}
			
			Next, we introduce the space of maximal regularity for $\{\boldsymbol{\varphi}, \boldsymbol{\varphi}'\}$ as \cite[p 2; Theorem 2.8.5.iii, p 17]{HS:2016}, \cite[p 1404-5]{GGH:2012}, with $T$ up to $\infty$:
			\begin{equation}\label{A-1.28}
			\mathbf{X}^T_{p,q, \sigma} = L^p(0,T;\calD(A_q)) \cap {W}^{1,p}(0,T;\lso)
			\end{equation}
			(recall \eqref{A-1.17}= (\ref{N1.9}) for $\calD(A_q)$) and the corresponding space for the pressure as
			\begin{equation}\label{A-1.29}
			Y^T_{p,q} = L^p(0,T;\widehat{W}^{1,q}(\Omega)), \quad \widehat{W}^{1,q}(\Omega) = W^{1,q}(\Omega) / \mathbb{R}.
			\end{equation}
			The following embedding, also called trace theorem, holds true \cite[Theorem 4.10.2, p 180, BUC for $T=\infty$]{HA:2000}, \cite{PS:2016}.
			\begin{equation}\label{A-1.30}
			\xtpqs \subset \xtpq \equiv L^p(0,T; \mathbf{W}^{2,q}(\Omega)) \cap W^{1,p}(0,T; \mathbf{L}^q(\Omega)) \hookrightarrow C \Big([0,T]; \Bso \Big).
			\end{equation}
			For a function $\mathbf{g}$ such that $div \ \mathbf{g} \equiv 0, \ \mathbf{g}|_{\Gamma} = 0$ we have $\mathbf{g} \in \xtpq \iff \mathbf{g} \in \xtpqs$.\\
			The solution of Eq (\ref{A-1.27}) is
			\begin{equation}\label{A-1.31}
			\boldsymbol{\varphi}(t) = e^{-A_qt} \boldsymbol{\varphi}_0 + \int_{0}^{t} e^{-A_q(t-s)} \Fs(\tau) d \tau.
			\end{equation}
			The following is the celebrated result on maximal regularity on $\lso$ of the Stokes problem due originally to Solonnikov \cite{VAS:1977} reported in \cite[Theorem 2.8.5.(iii) and Theorem 2.10.1 p24 for $\boldsymbol{\varphi}_0 = 0$]{HS:2016}, \cite{S:2006}, \cite[Proposition 4.1 , p 1405]{GGH:2012}.
			\begin{thm}\label{A-Thm-1.5}
				Let $1 < p,q < \infty, T \leq \infty$. With reference to problem (\ref{A-1.26}) = (\ref{A-1.27}), assume
				\begin{equation}\label{A-1.32}
				\Fs \in L^p(0,T;\lso), \ \boldsymbol{\varphi}_0 \in \Big( \lso, \calD(A_q)\Big)_{1-\frac{1}{p},p}.
				\end{equation}
				Then there exists a unique solution $\boldsymbol{\varphi} \in \xtpqs, \pi \in \ytpq$ to the dynamic Stokes problem (\ref{A-1.26}) or (\ref{A-1.27}), continuously on the data: there exist constants $C_0, C_1$ independent of $T, \Fs, \boldsymbol{\varphi}_0$ such that via (\ref{A-1.30})
				\begin{equation}\label{A-1.33}
				\begin{aligned}
				C_0 \norm{\boldsymbol{\varphi}}_{C \big([0,T]; \Bso \big)} &\leq \norm{\boldsymbol{\varphi}}_{\xtpqs} +  \norm{\pi}_{\ytpq}\\ &\equiv \norm{\boldsymbol{\varphi}'}_{L^p(0,T;\lso)} + \norm{A_q \boldsymbol{\varphi}}_{L^p(0,T;\lso)} +  \norm{\pi}_{\ytpq}\\
				&\leq C_1 \bigg \{ \norm{\Fs}_{L^p(0,T;\lso)}  + \norm{\boldsymbol{\varphi}_0}_{\big( \lso, \calD(A_q)\big)_{1-\frac{1}{p},p}} \bigg \}.
				\end{aligned}
				\end{equation}
				In particular,
				\begin{enumerate}[(i)]
					\item With reference to the variation of parameters formula (\ref{A-1.31}) of problem (\ref{A-1.27}) arising from the Stokes problem (\ref{A-1.26}), we have recalling (\ref{A-1.28}): the map
					\begin{align}
					\Fs &\longrightarrow \int_{0}^{t} e^{-A_q(t-\tau)}\Fs(\tau) d\tau \ : \text{continuous} \label{A-1.34}\\
					L^p(0,T;\lso) &\longrightarrow \xtpqs \equiv L^p(0,T; \calD(A_q)) \cap W^{1,p}(0,T; \lso) \label{A-1.35}.				
					\end{align}			
					\item The s.c. analytic semigroup $e^{-A_q t}$ generated by the Stokes operator $-A_q$ (see \eqref{A-1.17}= (\ref{N1.9})) on the space $\ds \Big( \lso, \calD(A_q)\Big)_{1-\frac{1}{p},p}$ (see statement below (\ref{A-1.24})) satisfies
					\begin{subequations}\label{A-1.36}
						\begin{multline}
						e^{-A_q t}: \ \text{continuous} \quad \Big( \lso, \calD(A_q)\Big)_{1-\frac{1}{p},p} \longrightarrow \xtpqs \equiv \\ L^p(0,T; \calD(A_q)) \cap W^{1,p}(0,T; \lso) \label{A-1.36a}.
						\end{multline}
						In particular via (\ref{A-1.24b}), for future use, for $1 < q < \infty, 1 < p < \frac{2q}{2q - 1}$, the s.c. analytic semigroup $\ds e^{-A_q t}$ on the space $\ds \Bto$, satisfies
						\begin{equation}
						e^{-A_q t}: \ \text{continuous} \quad \Bto \longrightarrow \xtpqs. \label{A-1.36b}
						\end{equation}
					\end{subequations}				
					\item Moreover, for future use, for $1 < q < \infty, 1 < p < \frac{2q}{2q - 1}$, then (\ref{A-1.33}) specializes to
					\begin{equation}\label{A-1.37}
					\norm{\boldsymbol{\varphi}}_{\xtpqs} + \norm{\pi}_{\ytpq} \leq C \bigg \{ \norm{\Fs}_{L^p(0,T;\lso)} + \norm{\boldsymbol{\varphi}_0}_{\Bto} \bigg \}.
					\end{equation}
				\end{enumerate}		
			\end{thm}
			
			\item \textbf{Maximal $L^p$ regularity on $\lso$ of the Oseen operator $\calA_q, \ 1 < p < \infty, \ 1 < q < \infty$, up to $T < \infty$.} We next transfer the maximal regularity of the Stokes operator $(-A_q)$ on $\lso$-asserted in Theorem \ref{A-Thm-1.5} into the maximal regularity of the Oseen operator $\calA_q = -\nu A_q - A_{o,q}$ in (\ref{A-1.18}) exactly on the same space $\xtpqs$ defined in (\ref{A-1.28}), however only up to $T < \infty$.
			
			\noindent Thus, consider the dynamic Oseen problem in $\{ \boldsymbol{\varphi}(t,x), \pi(t,x) \}$ with equilibrium solution $y_e$, see (\ref{1.2}):		
			\begin{subequations}\label{A-1.38}
				\begin{align}
				\boldsymbol{\varphi}_t - \Delta \boldsymbol{\varphi} + L_e(\boldsymbol{\varphi}) + \nabla \pi &= F &\text{ in } (0, T] \times \Omega \equiv Q \label{A-1.38a}\\		
				div \ \boldsymbol{\varphi} &\equiv 0 &\text{ in } Q\\
				\begin{picture}(0,0)
				\put(-110,7){$\left\{\rule{0pt}{35pt}\right.$}\end{picture}
				\left. \boldsymbol{\varphi} \right \rvert_{\Sigma} &\equiv 0 &\text{ in } (0, T] \times \Gamma \equiv \Sigma\\
				\left. \boldsymbol{\varphi}\right \rvert_{t = 0} &= \boldsymbol{\varphi}_0 &\text{ in } \Omega,
				\end{align}
			\end{subequations}
			\begin{equation}
			L_e(\boldsymbol{\varphi}) = (y_e . \nabla) \boldsymbol{\varphi} + (\boldsymbol{\varphi}. \nabla) y_e \hspace{3cm} \label{A-1.39}
			\end{equation}
			rewritten in abstract form, after applying the Helmholtz projector $P_q$ to (\ref{A-1.38a}) and recalling $\calA_q$ in (\ref{A-1.18}), as
			\begin{equation}\label{A-1.40}
			\boldsymbol{\varphi}_t = \calA_q \boldsymbol{\varphi} + P_q \mathbf{F} = - \nu A_q \boldsymbol{\varphi} - A_{o,q}\boldsymbol{\varphi}+ \Fs, \quad \boldsymbol{\varphi}_0 \in \big( \lso, \calD(A_q)\big)_{1-\frac{1}{p},p}
			\end{equation}
			whose solution is
			\begin{equation}\label{A-1.41}
			\boldsymbol{\varphi}(t) = e^{\calA_qt} \boldsymbol{\varphi}_0 + \int_{0}^{t} e^{\calA_q(t-\tau)} \Fs(\tau) d \tau,
			\end{equation}
			\begin{equation}\label{A-1.42}
			\boldsymbol{\varphi}(t) = e^{-\nu A_qt}\boldsymbol{\varphi}_0 + \int_{0}^{t} e^{-\nu A_q(t-\tau)} \Fs(\tau) d \tau - \int_{0}^{t} e^{- \nu A_q(t-\tau)} A_{o,q} \boldsymbol{\varphi}(\tau) d \tau.
			\end{equation}
			
			\begin{thm}\label{A-Thm-1.6}
				Let $1 < p,q < \infty, \ 0 < T < \infty$. Assume (as in (\ref{A-1.32}))
				\begin{equation}\label{A-1.43}
				\Fs \in L^p \big( 0, T; \mathbf{L}^q_{\sigma} (\Omega) \big), \quad \boldsymbol{\varphi}_0 \in \lqaq
				\end{equation}
				where $\calD(A_q) = \calD(\calA_q)$, see \eqref{A-1.18}=(\ref{1.9.new}). Then there exists a unique solution $\boldsymbol{\varphi} \in \xtpqs, \ \pi \in \ytpq$ of the dynamic Oseen problem (\ref{A-1.38}), continuously on the data: that is, there exist constants $C_0, C_1$ independent of $\Fs, \boldsymbol{\varphi}_0$ such that
				\begin{align}
				C_0 \norm{\boldsymbol{\varphi}}_{C \big([0,T]; \Bso \big)} &\leq \norm{\boldsymbol{\varphi}}_{\xtpqs} + \norm{\pi}_{\ytpq}\nonumber \\ &\equiv \norm{\boldsymbol{\varphi}'}_{L^p(0,T;\mathbf{L}^q(\Omega))} + \norm{A_q \boldsymbol{\varphi}}_{L^p(0,T;\mathbf{L}^q(\Omega))} + \norm{\pi}_{\ytpq}\\
				&\leq C_T \left \{ \norm{\Fs}_{L^p(0,T;\lso)}  + \norm{\boldsymbol{\varphi}_0}_{\lqaq} \right \}
				\end{align}
				where $T < \infty$. Equivalently, for $1 < p, q < \infty$
				\begin{enumerate}[i.]
					\item The map
					\begin{equation}
					\begin{aligned}
					\Fs \longrightarrow \int_{0}^{t} e^{\calA_q(t-\tau)}\Fs(\tau) d\tau \ : \text{continuous}&\\
					L^p(0,T;\lso) &\longrightarrow L^p \big(0,T;\calD(\calA_q) = \calD(A_q) \big)\label{A-1.46}
					\end{aligned}			
					\end{equation}
					where then automatically, see (\ref{A-1.40})
					\begin{equation}
					L^p(0,T;\lso) \longrightarrow W^{1,p}(0,T;\lso) \label{A-1.47}
					\end{equation}
					and ultimately via \eqref{A-1.28}
					\begin{equation}
					L^p(0,T;\lso) \longrightarrow \xtpqs \equiv L^p \big(0,T;\calD(A_q) \big) \cap W^{1,p}(0,T;\lso). \label{A-1.48}
					\end{equation}
					\item The s.c. analytic semigroup $e^{\calA_q t}$ generated by the Oseen operator $\calA_q$ (see \eqref{A-1.18}=(\ref{1.9.new})) on the space $\ds \lqaq $ satisfies for $1 < p, q < \infty$
					\begin{equation}
					e^{\calA_q t}: \ \text{continuous} \quad \lqaq \longrightarrow L^p \big(0,T;\calD(\calA_q) = \calD(A_q)  \big) \label{A-1.49}
					\end{equation}
					and hence automatically by (\ref{A-1.28})
					\begin{equation}
					e^{ \calA_q t}: \ \text{continuous} \quad \lqaq \longrightarrow \xtpqs. \label{A-1.50}
					\end{equation}
					In particular, for future use, for $1 < q < \infty, 1 < p < \frac{2q}{2q - 1}$, we have that the s.c. analytic semigroup $\ds e^{\calA_q t}$ on the space $\ds \Bto$, satisfies
					\begin{equation}
					e^{\calA_q t}: \ \text{continuous} \quad \Bto \longrightarrow L^p \big(0,T;\calD(\calA_q) = \calD(A_q)  \big), \ T < \infty. \label{A-1.51}
					\end{equation}
					and hence automatically
					\begin{equation}
					e^{ \calA_q t}: \ \text{continuous} \quad \Bto \longrightarrow \xtpqs , \ T < \infty. \label{A-1.52}
					\end{equation}
				\end{enumerate}
			\end{thm}
			\noindent A proof is given in \cite[Appendix B]{LPT.1}.
		\end{enumerate}
		\section{The UCP of Theorem \ref{T1.4} implies the controllability rank condition \eqref{bbN4-5}}\label{app-B}
		\setcounter{equation}{0}
		\setcounter{theorem}{0}
		\renewcommand{\theequation}{{\rm B}.\arabic{equation}}
		\renewcommand{\thetheorem}{{\c B}.\arabic{theorem}}
		
		We return to equations (\ref{2.3}) and (\ref{2.4}) giving eigenvalues/vectors of the operator $ 	\BA_q $ in (\ref{1.18}) and its adjoint $ \BA_q^* $:
		
		\begin{align}
		\BA_q \boldsymbol{\Phi}_{ij} &= \lambda_i \boldsymbol{\Phi}_{ij} \in \calD(\BA_q) = [\mathbf{W}^{2,q}(\Omega) \cap \mathbf{W}^{1,q}_0(\Omega) \cap \lso] \times [W^{2,q}(\Omega) \cap W^{1,q}_0(\Omega)]  \label{B2.3}\\
		\BA_q^* \boldsymbol{\Phi}_{ij}^* &= \bar{\lambda}_i \boldsymbol{\Phi}_{ij}^* \in \calD(\BA_q^*) = [\mathbf{W}^{2,q'}(\Omega) \cap \mathbf{W}^{1,q'}_0(\Omega) \cap \lo{q'}] \times [W^{2,q'}(\Omega) \cap W^{1,q'}_0(\Omega)]. \label{B2.4}
		\end{align}
		\noindent \underline{The operator $ \BA_q$ :} We recall from (\ref{1.18}), \eqref{N1-17}, \eqref{N1-18} that	
		\begin{multline}\label{B1.18}
		\BA_q = \bbm \calA_q & -\calC_{\gamma} \\ -\calC_{\theta_e} & \calB_q \ebm : \mathbf{W}^q_{\sigma}(\Omega) = \lso \times \lqo \supset \calD(\BA_q) = \calD(\calA_q) \times \calD(\calB_q) \\ = (\mathbf{W}^{2,q}(\Omega) \cap \mathbf{W}^{1,q}_{0}(\Omega) \cap \lso) \times (W^{2,q}(\Omega) \cap W^{1,q}_{0}(\Omega)) \longrightarrow \mathbf{W}^q_{\sigma}(\Omega).
		\end{multline}
		\begin{align}
		\quad \calC_{\gamma} h &= -\gamma P_q (h \mathbf{e}_d), \quad \calC_{\gamma} \in \calL (L^q(\Omega),\lso), \label{BN1-17}   \\
		\quad \calC_{\theta_e} \mathbf{z} &=  \mathbf{z} \cdot \nabla \theta_e, \quad \calC_{\theta_e} \in \calL(\lso,\lqo). \label{BN1-18}
		\end{align}
		With $ \boldsymbol{\Phi} = [\boldsymbol{\varphi}, \psi ],  $ the PDE-version of $ \BA_q  \boldsymbol{\Phi} = \lambda \boldsymbol{\Phi} $ is
		\begin{subequations}\label{B2.111}
			\begin{align}
			- \nu \Delta \boldsymbol{\varphi} + L_e(\boldsymbol{\varphi}) + \nabla \pi - \gamma \psi \mathbf{e}_d&= \lambda \boldsymbol{\varphi}   &\text{in } \Omega      \label{B2.11a} \\
			- \kappa \Delta \psi + \mathbf{y}_e \cdot \nabla \psi + \boldsymbol{\varphi} \cdot \nabla \theta_e &= \lambda \psi &\text{in } \Omega \label{B2.11b}\\
			\begin{picture}(0,0)
			\put(-140,3){ $\left\{\rule{0pt}{36pt}\right.$}\end{picture}
			\text{div } \boldsymbol{\varphi} &= 0   &\text{in } \Omega  \label{B2.11c} \\
			\boldsymbol{\varphi} = 0, \ \psi &= 0 &\text{on } \Gamma.    \label{B2.11e}
			\end{align}		
		\end{subequations}\\
		
		\noindent \underline{The operator $\BA_q^*$:} We return to the Helmholtz decomposition in (\ref{1.3}), (\ref{1.4}) and provide additional information. For $\mathbf{M} \subset\mathbf{L}^q(\Omega), \ 1 < q < \infty$, we denote the annihilator of $\mathbf{M}$ by
		
		\begin{equation}
		\mathbf{M}^{\perp} = \bigg \{ \mathbf{f} \in \mathbf{L}^{q'}(\Omega) : \int_{\Omega} \mathbf{f} \cdot \mathbf{g} \ d \Omega = 0, \text{ for all } \mathbf{g} \in \mathbf{M} \bigg \}
		\end{equation}
		\noindent where $q'$ is the dual exponent of $\ds q: \ \rfrac{1}{q} + \rfrac{1}{q'} = 1$.
		
		\begin{prop}\cite[Prop 2.2.2 p6]{HS:2016}, \cite[Ex. 16 p115]{Ga:2011}\label{I-Prop-1.2}, \cite{FMM:1998} \ \\
			Let $\Omega \subset \BR^d$ be an open set and let $1 < q < \infty$.
			\begin{enumerate}[a)]
				\item The Helmholtz decomposition exists for $\mathbf{L}^q(\Omega)$ if and only if it exists for $\mathbf{L}^{q'}(\Omega)$, and we have: (adjoint of $P_q$) = $P^*_q = P_{q'}$ (in particular $P_2$ is orthogonal), where $P_q$ is viewed as a bounded operator $\ds \mathbf{L}^q(\Omega) \longrightarrow \mathbf{L}^q(\Omega)$, and $\ds P^*_q = P_{q'}$ as a bounded operator $\ds \mathbf{L}^{q'}(\Omega) \longrightarrow \mathbf{L}^{q'}(\Omega), \ \rfrac{1}{q} + \rfrac{1}{q'} = 1$.
				\item Then, with reference to (\ref{1.4})
				\begin{subequations}
					\begin{equation}
					\Big[ \lso \Big]^{\perp} = \mathbf{G}^{q'}(\Omega) \text{ and } \Big[ \mathbf{G}^q(\Omega) \Big]^{\perp} = \mathbf{L}^{q'}_{\sigma}(\Omega) \label{I-A.2a}.
					\end{equation}
					\begin{rmk} \label{rmkB.1}
						\noindent Throughout the paper we shall use freely that
						\begin{equation}
						\big( \lso \big)' = \lo{q'}, \quad \frac{1}{q} + \frac{1}{q'} = 1 \label{I-A.2b}
						\end{equation}
						\noindent This can be established as follows. From (\ref{1.4}) write $\ds \lso$ as a factor space $\ds \lso = \mathbf{L}^q(\Omega) / \mathbf{G}^q(\Omega) \equiv \mathbf{X}/\mathbf{M}$ so that \cite[p 135]{TL:1980}.
						\begin{equation}
						\big( \lso \big)' = \big( \mathbf{L}^q(\Omega) / \mathbf{G}^q(\Omega) \big)' =  \big( \mathbf{X}/\mathbf{M} \big)' = \mathbf{M}^{\perp} = \Big[ \mathbf{G}^q(\Omega) \Big]^{\perp} = \lo{q'} \label{I-A.2c}.
						\end{equation}
						\noindent In the last step, we have invoked (\ref{I-A.2a}), which is also established in \cite[Lemma 2.1, p 116]{Ga:2011}. Similarly
						\begin{equation}
						\big(\mathbf {G}^q(\Omega) \big)' = \big(\mathbf{ L}^q(\Omega) / \lso \big)' = \Big[ \lso \Big]^{\perp} =\mathbf{ G}^{q'}(\Omega).
						\end{equation}
					\end{rmk}
				\end{subequations}			
			\end{enumerate}
		\end{prop}
		
		\noindent Next, let $ \boldsymbol{\varphi}^* \in \mathbf{L}^{q'}_{\sigma}(\Omega) $ and $ h \in \L^q( \Omega ) $. From (\ref{BN1-17}), we then compute in the noted duality pairings, via the above results, omitting at times the symbol $ \Omega $:
		\begin{align}
		\big( \calC_{\gamma} h , \boldsymbol{\varphi}^* \big)_{\lso,\mathbf{L}^{q'}_{\sigma}(\Omega)} & = \big( - \gamma P_q ( h  \mathbf{e}_d),  \boldsymbol{\varphi}^* \big)_{\lso,\mathbf{L}^{q'}_{\sigma}(\Omega)}=\big( ( h  \mathbf{e}_d),- \gamma(P_{q'} \boldsymbol{\varphi}^*) \big)_{\mathbf{L}^q, \mathbf{L}^{q'}}\\
		&=\big( h ,- \gamma(P_{q'} \boldsymbol{\varphi}^*) \cdot\mathbf{e}_d) \big)_{L^q, L^{q'}}  =  \big( h, \calC_{\gamma}^* \boldsymbol{\varphi}^* \big)_{L^q, L^{q'}}
		\end{align}
		Thus
		\begin{equation}\label{B-11}
		\calC_{\gamma}^* \boldsymbol{\varphi}^* =  - \gamma(P_{q'} \boldsymbol{\varphi}^*) \cdot\mathbf{e}_d ,\,\,   \calC_{\gamma}^* \in \mathcal{L}( \lo{q'} , L^{q'}( \Omega )).
		\end{equation}
		Next, let $\mathbf{z} \in \lso $, so that $ \div\,\mathbf{z} \equiv 0 , \mathbf{z} \cdot \nu = 0 $ on $ \Gamma $ by (\ref{1.3}). Let $ \psi^* \in L^{q'} ( \Omega ) $. For $ q > d $, recall $ \nabla \theta_e \in \mathbf{W}^{1,q}(\Omega) \hookrightarrow \mathbf{L}^{\infty}( \Omega ) $ from Theorem \ref{Thm-1.1}. Then, from (\ref{BN1-18}) we compute
		
		\begin{align}
		\big( \calC_{\theta_{e}} \mathbf{z} , \psi^* \big)_{L^q,L^{q'}} & = \big( \mathbf{z} \cdot \nabla \theta_e , \psi^* \big)_{L^q,L^{q'}} = \big( ( P_{q} \mathbf{z}) \cdot \nabla \theta_e , \psi^* \big)_{L^q,L^{q'}} \\
		&= \int_{\Omega} (P_{q} \mathbf{z}) \cdot ( \psi^* \nabla \theta_e)   \, d\Omega = \big( P_q \mathbf{z} , \psi^* \nabla \theta_e \big)_{L^q,L^{q'}} \\
		&=\big( \mathbf{z} , P_{q'} ( \psi^* \nabla \theta_e ) \big)_{\lso,\mathbf{L}^{q'}_{\sigma}(\Omega)} = \big( \mathbf{z}, \calC_{\theta_{e}}^* \psi^* \big)_{\lso,\mathbf{L}^{q'}_{\sigma}(\Omega)}.
		\end{align}
		Thus,
		\begin{equation}\label{B-15}
		\calC_{\theta_{e}}^* \psi^*= P_{q'} ( \psi^* \nabla \theta_e ) \in \mathbf{L}^{q'}_{\sigma}(\Omega).
		\end{equation}
		The adjoint $\BA_q^* $ of $ \BA_q $ in (\ref{B1.18}) is, by virtue of (\ref{B-11}) and (\ref{B-15})
		\begin{multline}\label{B11.18}
		\BA_q^* = \bbm \calA_q^* & -\calC_{\theta_e}^* \\ -\calC_{\gamma}^* & \calB_q^* \ebm : \mathbf{W}^{q'}_{\sigma}(\Omega) = \lsoo \times \lqoo \supset \calD(\BA_q^*) = \calD(\calA_q^*) \times \calD(\calB_q^*) \\ = (\mathbf{W}^{2,q'}(\Omega) \cap \mathbf{W}^{1,q'}_{0}(\Omega) \cap \lsoo) \times (W^{2,q'}(\Omega) \cap \mathbf{W}^{1,q'}_{0}(\Omega)) \longrightarrow \mathbf{W}^{q'}_{\sigma}(\Omega).
		\end{multline}
		
		With $ \boldsymbol{\Phi^*} = [\varphi^*, \psi^* ],  $ the PDE-version of  $ \BA_q^*  \boldsymbol{\Phi^*} = \lambda \boldsymbol{\Phi^*} $ is
		\begin{subequations}\label{BA2.111}
			\begin{align}
			-\nu \Delta \boldsymbol{\varphi^*} + L^*_e(\boldsymbol{\varphi^*})+ \psi^* \nabla \theta_e + \nabla \pi &= \lambda \boldsymbol{\varphi^*}   &\text{in } \Omega      \label{BA2.11a} \\
			-\kappa \Delta \psi^* + \mathbf{y}_e \cdot \nabla \psi^* - \gamma  \boldsymbol{\varphi^*}  \cdot\mathbf{e}_d &= \lambda \psi^* &\text{in } \Omega \label{BA2.11b}\\
			\begin{picture}(0,0)
			\put(-138,4){ $\left\{\rule{0pt}{35pt}\right.$}\end{picture}
			\text{div } \boldsymbol{\varphi^*} &= 0   &\text{in } \Omega  \label{BA2.11c} \\
			\boldsymbol{\varphi^*} = 0, \ \psi^* &= 0 &\text{on } \Gamma.    \label{BA2.11e} 
			\end{align}		
		\end{subequations}\\
		\uline{Proof that the UCP of Theorem \ref{T1.4} implies the controllability rank condition} \eqref{bbN4-5}). We return to (\ref{B2.3}) and (\ref{B2.4}) and express the eigenvectors in terms of their coordinates, as $(d+1)$ vectors.
		\begin{equation}
		\boldsymbol{\Phi_{ij}} = \{ \boldsymbol{\varphi}_{ij} , \psi_{ij} \} = \{ \varphi_{ij}^{(1)} ,\varphi_{ij}^{(2)},...,\varphi_{ij}^{(d-1)}, \varphi_{ij}^{(d)}, \psi_{ij} \} \,, \text{ a } (d+1)\text{-vector}
		\end{equation}
		\begin{equation} \label{B-19}
		\boldsymbol{\Phi^*_{ij}} = \{ \boldsymbol{\varphi^*_{ij}} , \psi^*_{ij} \} = \{ \varphi_{ij}^{*(1)} ,\varphi_{ij}^{*(2)},...,\varphi_{ij}^{*(d-1)}, \varphi_{ij}^{*(d)}, \psi^*_{ij} \} \,, \text{ a } (d+1)\text{-vector}.
		\end{equation}
		With reference to (\ref{B-19}), we introduce the following corresponding $d$-vector
		\begin{equation} \label{B-20}
		\boldsymbol{\widehat{\Phi }^*_{ij} = \{ \boldsymbol{\widehat{\varphi}^*_{ij}}} , \psi^*_{ij} \} = \{ \varphi_{ij}^{*(1)} ,\varphi_{ij}^{*(2)},...,\varphi_{ij}^{*(d-1)},  \psi^*_{ij} \} \,, \text{ a } d\text{-vector}
		\end{equation}
		obtained from $  \boldsymbol{\Phi^*_{ij}}$ by omitting the d-component $\varphi_{ij}^{*(d)}$ of the vector $\boldsymbol{\Phi^*_{ij}} $. Next, with reference to \eqref{Bb2.5},
		\noindent construct the following matrix $U_i$ of size $\ell_i \times K, \ K =\sup \{ \ell_i: i = 1, \dots, M \}$
		\begin{equation}\label{B2.5}
		U_i =
		\begin{bmatrix}
		(\mathbf{u}_1,\widehat{\Phi}_{i1}^*)_{\omega} & \dots & (\mathbf{u}_K,\widehat{\Phi}_{i1}^*)_{\omega} \\[1mm]
		(\mathbf{u}_1,\widehat{\Phi}_{i2}^*)_{\omega} & \dots & (\mathbf{u}_K,\widehat{\Phi}_{i2}^*)_{\omega} \\
		\vdots & \ddots & \vdots \\
		(\mathbf{u}_1,\widehat{\Phi}_{i \ell_i}^*)_{\omega} & \dots & (\mathbf{u}_K,\widehat{\Phi}_{i \ell_i}^*)_{\omega} \\
		\end{bmatrix}
		:\ell_i \times K.
		\end{equation}
		Here with  \begin{subequations}

			\begin{equation}
			\mathbf{u}_k = [ \mathbf{u}_k^1 , u_k^2 ] = [ ( u_k^1)^{(1)}, (u_k^1)^{(2)}...(u_k^1)^{(d-1)},u_k^2] \in \mathbf{\widehat{L}}^q_{\sigma}( \Omega) \times L^q(\Omega)
			\end{equation}
			\begin{equation}
			\widehat{ \mathbf{L}}^q_{\sigma} ( \Omega ) \equiv \text{ the space obtained from } \mathbf{L}^q_{\sigma} ( \Omega ) \text{ after omitting the } d\text{-coordinate},
			\end{equation}
		\end{subequations}
		we have defined the duality pairing over $ \omega $ as
		
		\begin{align}
		(\mathbf{u}_k, \widehat{\Phi}^*_{ij})_{\omega} &= \bpm \bbm \mathbf{u}^1_k \\ u^2_k \ebm, \bbm \widehat{\boldsymbol{\varphi}}^*_{ij} \\ \psi^*_{ij} \ebm \epm_{\omega} = \int_{\omega} [ \mathbf{u}^1_k \cdot \widehat{\boldsymbol{\varphi}}^*_{ij} + u^2_k \psi^*_{ij} ] d \omega \nonumber \\
		&= (\mathbf{u}^1_k , \widehat{\boldsymbol{\varphi}}^*_{i1})_{\widehat{\mathbf{L}}^q(\omega), \widehat{\mathbf{L}}^{q'}(\omega)} + ( u^2_k, \psi^*_{ij})_{L^q(\omega), L^{q'}(\omega)}\\
		&= \int_{\omega }
		\begin{bmatrix}
		( u_k^1)^{(1)}\\
		( u_k^1)^{(2)}\\
		\vdots \\
		( u_k^1)^{(d-1)}\\
		u_k^2
		\end{bmatrix} \cdot
		\begin{bmatrix}
		\varphi_{ij}^{*(1)}\\
		\varphi_{ij}^{*(2)}\\
		\vdots \\
		\varphi_{ij}^{*(d-1)}\\
		\psi^*_{ij}
		\end{bmatrix} \, d\omega.
		\end{align}
		The controllability Kalman/Hautus algebraic condition of the finite-dimensional projected $\mathbf{w}_N$-equation in (\ref{1.27a}) is given by \eqref{bbN4-5}
		\begin{equation}    \label{bN4-5}
		\text{rank } U_i = \text{full} = \ell_i, \quad i = 1, \dots, M.
		\end{equation}
		$M$= number of distinct unstable eigenvalues in (\ref{1.21}). Thus, given the (unstable) eigenvalues $\bar{\lambda}_i$ of $\BA_q^* \,, i = 1,2,...,M$ we need to show that the corresponding vectors $  \boldsymbol{\widehat{\Phi }}^*_{i1},..\boldsymbol{\widehat{\Phi }}^*_{il_{i}}$ (defined in (\ref{B-20}))
		
		\begin{equation}        \label{E3-7}
		\begin{tikzpicture}
		\draw (0,0) node {$\bar{\lambda}_i$};
		\draw (-0.4,-0.4) -- (-1,-1) node[anchor = north east] {$\widehat{\Phi}^*_{i1}$,};
		\draw (-0.2,-0.4) -- (-0.4,-1) node[anchor = north] {$\widehat{\Phi}^*_{i2},$};
		\draw (0.4,-1.3) node {$\ldots$};
		\draw (0.4,-0.4) -- (1,-1) node[anchor = north west] {$\widehat{\Phi}^*_{i \ell_i}$ };
		\end{tikzpicture}
		\end{equation}
		are linearly independent in $\widehat{\mathbf{L}}^{q'}(\omega)$, where $\ell_i = $ geometric multiplicity of $\lambda_i$. Thus, we seek to establish that
		\begin{equation}    \label{E3-8}
		\sum_{j=1}^{\ell_i} \alpha_j \widehat{\Phi^*}_{ij} \equiv 0 \mbox{ in } \widehat{\mathbf{L}}^{q'}(\omega) \times L^{q'}(\omega)
		\quad \Longrightarrow \quad
		\alpha_j = 0, j= 1, \ldots, \ell_i.
		\end{equation}
		To this end, define the vector $\Phi=[\boldsymbol{\varphi^*}, \psi^*]$ (we suppress dependence on $i$) by
		\begin{equation}
		\label{E3-9}
		\Phi^* = \sum_{j=1}^{\ell_i} \alpha_j \Phi^*_{ij} \quad \mbox{ in } \mathbf{L}^{q'}(\Omega) \times L^{q'}( \Omega )
		\end{equation}
		that is, $\Phi^* $ is a linear combination, with the same constants $ \alpha_j $ as in (\ref{E3-8}), of the eigenvectors $ \Phi^*_{i1},...,\Phi^*_{il_{i}}$ in (\ref{B2.4}). Then $ \boldsymbol{\Phi^*} = [\boldsymbol{\varphi}^*, \psi^* ]$ is itself an eigenvector of $\BA_q^*$ corresponding to the eigenvalue $\bar{\lambda}_i$. Thus we have
		\begin{equation}
		\label{E3-10}
		\BA_q^* \boldsymbol{\Phi}^* = \bar{\lambda}_i \boldsymbol{\Phi}^*  \mbox{ in } \widehat{\mathbf{L}}^{q'}(\Omega) \times L^{q'}( \Omega ) ; \quad \boldsymbol{\varphi^*} \in\mathbf{L}^{q'}( \Omega ) \,, \psi^* \in L^{q'}(\Omega) .
		\end{equation}
		along with
		\begin{equation} \label{B-30}
		\widehat{\boldsymbol{\Phi}}^* \equiv 0 \text{ in } \omega \text{ by } (\ref{E3-8}).
		\end{equation}
		Therefore $ \boldsymbol{\boldsymbol{\Phi}^*} = [\boldsymbol{\varphi}^*, \psi^* ],  $ satisfies the eigenvalue problem (\ref{BA2.111}a-b-c-d), along with the over-determined condition in $ \omega $, due to (\ref{B-30})
		\begin{equation}
		\boldsymbol{\varphi}^*= \{ \varphi^{*(1)} ,  \varphi^{*(2)},...,\varphi^{*(d-1)}\} \equiv 0 \text{ in } \omega \,,\qquad \psi^* \equiv 0 \text{ in } \omega.
		\end{equation}
		According to Theorem \ref{T1.4}, then
		\begin{equation}
		\label{E3-12}
		\boldsymbol{\Phi^*}=[\boldsymbol{\varphi^*}, \psi] \equiv 0 \quad \mbox{ in } \Omega.
		\end{equation}
		Indeed, in order to reach conclusion (\ref{E3-12}), one does not need the B.C. (\ref{BA2.11e}) on $ \Gamma $, as this is not needed in Theorem \ref{T1.4}. Going back to (\ref{E3-9}), we have
		\begin{equation}
		\label{E3-13}
		\Phi^* = \sum_{j=1}^{\ell_i} \alpha_j \Phi^*_{ij}  \equiv 0 \mbox{ in }  \Omega ; \mbox{ hence,  } \alpha_1 = 0, \cdots, \alpha_{\ell_i} = 0, \quad
		i=1, \ldots, M
		\end{equation}
		
		\noindent since the eigenvectors $\left \{ \Phi^*_{ij} \right \}_{j=1}^{\ell_i} $ are linearly independent in $ \mathbf{L}^{q'}(\Omega) \times L^{q'}(\Omega)$. Conclusion (\ref{E3-8}) is established. (The fact that $\bar{\lambda}_i$ is unstable, is irrelevant in the above argument.) Due to the just established property \eqref{E3-8}, we then see that the algebraic conditions \eqref{bN4-5} can be satisfied by infinitely many choices of the vectors $\mathbf{u}_1, \ldots, \mathbf{u}_K \in \widehat{\mathbf{L}}^{q'}(\Omega) \times L^{q'}( \Omega )$.\\
		
		\noindent \textbf{Acknowledgment}.
		The authors wish to thank two referees for their thorough reading of the original version and for their comments that have led to an improved exposition, with additional references. In particular they have raised stimulating questions regarding possible generalizations of the mathematical treatment. Some of their questions are addressed and answered in Section \ref{subsec_1.4} (points 3, 5, 6, 7). The proof in Appendix \ref{app-B} that the UCP of Theorem \ref{T1.4} implies the controllability rank condition \eqref{bbN4-5} was inserted upon request of one referee. In addition, the authors would like to thank the following people for providing references for Theorem \ref{Thm-1.1} for $ q \ne 2 $: Professors Andrea Giorgini, Giusy Mazzone, Hyunseok Kim.\\

	\end{appendices}

\end{document}